%% file: AH-mass-aspect-part-1.tex
\documentclass[reqno,11pt]{amsart}
\usepackage{faktor}
\usepackage{lscape}
\usepackage[usenames,dvipsnames]{xcolor}
\usepackage{enumerate,amsmath,amssymb,stmaryrd,mathtools} 
\usepackage{pstricks,pst-node,pst-coil,pst-plot}
\usepackage{tikz}
\usepackage{tikz-cd}
\usepackage{genyoungtabtikz}

\usetikzlibrary{arrows,decorations.markings,decorations.pathreplacing,patterns}
\usetikzlibrary{matrix}
\usetikzlibrary{graphs}
\usetikzlibrary{backgrounds}

\usepackage[space, compress, sort]{cite}
\usepackage[colorlinks=true, bookmarks=true, pdfstartview=FitH, pagebackref=true]{hyperref}

\usepackage{verbatim}

\input{def}

\numberwithin{equation}{section}

\def\DynkinNodeSize{2mm}
\def\DynkinArrowLength{3mm}
\def\BrauerNodeSize{1.5mm}
\tikzset{
  dnode/.style={
    circle,
    inner sep=0pt,
    minimum size=\DynkinNodeSize,
    fill=white,
    draw},
  middlearrow/.style={
    decoration={markings,
      mark=at position 0.6 with
      {\draw (0:0mm) -- +(+135:\DynkinArrowLength); \draw (0:0mm) -- +(-135:\DynkinArrowLength);},
    },
    postaction={decorate}
  },
  leftrightarrow/.style={
    decoration={markings,
      mark=at position 0.999 with
      {
      \draw (0:0mm) -- +(+135:\DynkinArrowLength); \draw (0:0mm) -- +(-135:\DynkinArrowLength);
      },
      mark=at position 0.001 with
      {
      \draw (0:0mm) -- +(+45:\DynkinArrowLength); \draw (0:0mm) -- +(-45:\DynkinArrowLength);
      },
    },
    postaction={decorate}
  },
  sedge/.style={
  },
  dedge/.style={
    middlearrow,
    double distance=1mm,
  },
  tedge/.style={
    middlearrow,
    double distance=1.0mm+\pgflinewidth,
    postaction={draw}, 
  },
  infedge/.style={
    leftrightarrow,
    double distance=0.5mm,
  },
  bnode/.style={
    circle,
    inner sep=0pt,
    minimum size=\BrauerNodeSize,
    fill=black,
    draw
  },
  bnode2/.style={
    circle,
    inner sep=0pt,
    minimum size=\BrauerNodeSize/2.5,
    fill=black,
    draw
  },
}

\allowdisplaybreaks

\begin{document}

\title[Mass-like invariants for AH metrics. Part 1: classification]
{Mass-like invariants for asymptotically hyperbolic metrics.\\~\\
  Part 1: classification}

\author{Julien Cortier}
\address{Institut Fourier, Universit\'e Grenoble-Alpes \\
  100 rue des Math\'ematiques \\
  38610 Gi\`eres \\
  France}
\email{julien-jean.cortier@ac-montpellier.fr}

\author{Mattias Dahl}
\address{Institutionen f\"or Matematik \\
  Kungliga Tekniska H\"ogskolan \\
  100 44 Stockholm \\
  Sweden} \email{dahl@kth.se}

\author{Romain Gicquaud}
\address{Institut Denis Poisson \\
  UFR Sciences et Technologie \\
  Facult\'e Fran\c cois Rabelais \\
  Parc de Grandmont \\
  37200 Tours \\
  France} \email{romain.gicquaud@univ-tours.fr}

\begin{abstract}
  In this article, we classify the set of asymptotic mass-like invariants for asymptotically hyperbolic metrics. It turns out that the standard mass is just one example among the two families of invariants we find. These invariants correspond to finite-dimensional representations of the group of isometries of hyperbolic space. We describe these invariants in terms of wave-harmonic polynomials and polynomial solutions of the linearized Einstein equations in Minkowski space.
\end{abstract}




\maketitle

\tableofcontents

\section{Introduction}\label{secIntroduction}

\input{introduction}

\section{Preliminaries and basic definitions}\label{secDefinition}

\input{definition}

\section{Action of hyperbolic isometries on mass aspect tensors}\label{secLorentzaction}

\input{lorentzaction}

\section{Classification of linear masses at infinity}\label{secClassification}

\input{classification}

\appendix
\section{On the irreducible representations that appear in the classification}\label{secRepresentation}

\input{weyl}

\bibliographystyle{amsplain}
\bibliography{biblio-AH-mass}

\end{document}

%% file: def.tex
\theoremstyle{plain}

\newtheorem{thm}{Theorem}[section]
\newtheorem{theorem}[thm]{Theorem}

\newtheorem{corollary}[thm]{Corollary}

\newtheorem{lemma}[thm]{Lemma}

\newtheorem{proposition}[thm]{Proposition}

\theoremstyle{remark}
\newtheorem{example}[thm]{Example}
\newtheorem{remark}[thm]{Remark}
\newtheorem{remarks}[thm]{Remarks}

\newtheorem{claim}{Claim}

\theoremstyle{definition}

\newtheorem{definition}[thm]{Definition}

\newcounter{mnotecount}[section]

\newcommand{\definedas}{\coloneqq}

\def\epsilon{{\varepsilon}}
\def\phi{{\varphi}}

\let\<\langle
\let\>\rangle

\newcommand{\Mbar}{\overline{M}}

\newcommand{\psibar}{\overline{\psi}}

\newcommand{\ebar}{\overline{e}}

\newcommand{\ftil}{\widetilde{f}}
\newcommand{\htil}{\widetilde{h}}

\newcommand{\bR}{\mathbb{R}}

\newcommand{\bS}{\mathbb{S}}

\newcommand{\bN}{\mathbb{N}}
\newcommand{\bH}{\mathbb{H}}
\newcommand{\bC}{\mathbb{C}}

\newcommand{\cB}{\mathcal{B}}
\newcommand{\cD}{\mathcal{D}}
\newcommand{\cI}{\mathcal{I}}
\newcommand{\cJ}{\mathcal{J}}
\newcommand{\cE}{\mathcal{E}}
\newcommand{\cC}{\mathcal{C}}

\newcommand{\cH}{\mathcal{H}}

\newcommand{\cK}{\mathcal{K}}

\newcommand{\cG}{\mathcal{G}}
\newcommand{\cN}{\mathcal{N}}

\newcommand{\cW}{\mathcal{W}}
\newcommand{\cT}{\mathcal{T}}

\newcommand{\fa}{\mathfrak{a}}

\newcommand{\fr}{\mathfrak{r}}

\newcommand{\fp}{\mathfrak{p}}

\newcommand{\Abar}{\overline{A}}

\newcommand{\Rbar}{\overline{R}}

\newcommand{\xhat}{\widehat{x}}

\newcommand{\gbar}{\overline{g}}

\renewcommand{\hbar}{\overline{h}}
\newcommand{\Gambar}{\overline{\Gamma}}
\newcommand{\Gamtil}{\widetilde{\Gamma}}
\newcommand{\gtil}{\widetilde{g}}
\newcommand{\mtil}{\widetilde{m}}
\newcommand{\Util}{\widetilde{U}}
\newcommand{\Vtil}{\widetilde{V}}

\newcommand{\bbar}{\overline{b}}
\newcommand{\zbar}{\overline{z}}

\newcommand{\omegabar}{\overline{\omega}}

\DeclareMathOperator{\tr}{tr}

\DeclareMathOperator{\divg}{div}

\newcommand{\lie}{\mathcal{L}}

\DeclareMathOperator{\Sym}{Sym}

\newcommand{\riem}{\mathcal{R}}

\DeclareMathOperator{\gradient}{grad}



%% file: introduction.tex
%
%

\subsection{Background}

Asymptotically hyperbolic manifolds are non-compact Riemannian manifolds whose
geometry, at infinity, approaches that of the hyperbolic space. Such manifolds
have attracted a lot of attention, in geometry as well as in theoretical
physics.

They naturally arise in general relativity, in the description of isolated
gravitational systems in a universe with negative cosmological constant (see,
for example~\cite{AshtekarMagnon84}), modeled on Anti-de~Sitter spacetime. One
can also see them as asymptotically umbilical hypersurfaces in asymptotically
Minkowski spacetimes. However, from the perspective of general relativity, one
must also consider the second fundamental form of the embedded hypersurface.

Properties of asymptotically hyperbolic metrics have been investigated,
motivated by the study of isolated systems in general relativity. If the
difference between such a metric $g$ and the hyperbolic metric $b$ decays at
infinity faster than some rate, part of the deviation is captured by the
so-called \emph{mass vector} $p_g \in \bR^{n+1}$. This concept is reminiscent
of the \emph{ADM mass} for asymptotically Euclidean manifolds. As in the case
of asymptotically Euclidean manifolds, a difficulty when studying the mass is
its apparent dependence on the choice of chart at infinity. This issue arises
in both Wang's approach~\cite{WangMass} and that of Chru\'sciel and
Herzlich~\cite{ChruscielHerzlich}. However, in the latter work, it is
established that a certain scalar quantity associated to the mass vector---the
mass---does not depend on the chosen chart at infinity, provided appropriate
asymptotic conditions are met.

The key result underlying this fact was first established by Chru\'sciel and
Nagy in~\cite[Theorem~3.3]{ChruscielNagy}. In particular, it states that the
transition diffeomorphism $\Psi = \phi_2 \circ \phi_1^{-1}$ between any two
charts at infinity $\phi_1$ and $\phi_2$ can be decomposed into a principal
part $A$, which is an isometry of hyperbolic space, and a correction term that
decays to zero at infinity. The decay of the correction term is sufficiently
rapid that it does not affect the components of the \emph{mass vector}, which
are computed at infinity. For this reason, the transition diffeomorphism $\Psi$
is referred to as an \emph{asymptotic isometry}. We denote by $\pi$ the map
$\Psi \mapsto A$ from the group of asymptotic isometries to the
finite-dimensional group of hyperbolic isometries.

Wang and Chru\'sciel–Herzlich subsequently showed that the mass vector $p_g$ of
any asymptotically hyperbolic metric of suitable decay satisfies the
equivariance property
\[
    p_{\Psi_* g} = \pi(\Psi) p_g ,
\]
where $\pi(\Psi)$ acts naturally as an element of the Lorentz group $O(n,1)$ on
$\bR^{n+1}$. It follows in particular that the norm in the Minkowski metric
$\eta$, $\eta(p_g,p_g)$, is invariant under the action of asymptotic
isometries.

For other types of asymptotically hyperbolic manifolds, geometric invariants
also appear, see, for example~\cite{CapGover, DahlKroenckeMcCormik,
    DjadliGuillarmouHerzlich, GrahamRenormalizedVolume}. However, these invariants
are of a different nature than those considered in this paper.

\subsection{Statement of results}

Inspired by the properties of the mass vector $p_g$, namely the facts that it
is defined solely in terms of the asymptotic geometry of the metric and that it
enjoys natural covariance properties, the goal of this paper is to define and
classify mass-like invariants for asymptotically hyperbolic manifolds whose
Riemannian metrics decay towards the standard metric of the hyperbolic space at
\emph{any} specified rate. We call them \emph{linear masses at infinity}.

Motivated by the example of the mass vector, we aim at finding invariants $g
    \mapsto \Phi(g)$, defined for asymptotically hyperbolic metrics $g$ in a
neighborhood of the infinity of the hyperbolic space $\bH^n$, taking values in
a finite-dimensional representation $V$ of the group of hyperbolic isometries.
Such a quantity is a \emph{linear mass at infinity} if the map $\Phi : g
    \mapsto \Phi(g) \in V$ is linear in $g - b$, where $b$ denotes the hyperbolic
metric, and satisfies the equivariance property
\begin{equation}\label{eqRoughMass}
    \Phi ( \Psi_* g ) = \pi( \Psi) \cdot \Phi(g)
\end{equation}
under the action of the group of asymptotic isometries.

We now describe our results. We will denote by $G_k$ the space of
asymptotically hyperbolic metrics $g$ defined outside of a compact set that
decay to order $k$ towards $b$ at infinity, meaning that the norm $|g-b|_b$
satisfies
\[
    |g-b|_b = O\left(e^{-k d(x)}\right),
\]
where $d(x) = \mathrm{dist}_b (x,x_0)$ is the $b$-distance between $x$ and a
given point $x_0 \in \bH^n$. A more precise definition will be given in
Section~\ref{secDefinition}. We will work with germs of metrics at infinity to
capture the fact that the invariants we are interested in are computed outside
of any compact sets.

It will be shown in Theorem~\ref{thmChruscielNagy} that the set $I^k(g)$ of
diffeomorphisms $\Psi$ (defined once again in a neighborhood of infinity) such
that $\Psi_* g \in G_k$ is a group for composition. This group is a semi-direct
product admitting the set of diffeomorphisms asymptotic to the identity as a
normal subgroup and the group of hyperbolic isometries as the second factor. In
particular, there is a well-defined morphism $\pi$ from $I^k(g)$ onto the
subgroup $O_\uparrow(n, 1)$ of the Lorentz group corresponding to the isometry
group of the hyperbolic space.

In order to get rid of the subgroup of diffeomorphisms tending to the identity,
we introduce the notion of transverse germs of metrics. These germs satisfy the
condition $g(P, \cdot) = b(P, \cdot)$, where $P = x^i \partial_i$ is a radial
vector field, so that the difference $g - b$ has only components tangential to
the level sets of $d$. Proposition~\ref{propTransverse} will show that, for any
given metric $g \in G_k$, there exists a unique diffeomorphism $\Theta$,
tending to the identity, called the \emph{adjustment diffeomorphism}, such that
$\Theta_* g$ is transverse. In particular, this implies that the group
$O_\uparrow(n, 1)$ acts on the set of transverse germs.
Proposition~\ref{propLinearity} shows that this action is actually linear in
$g-b$ up to some high order.

To be able to formulate a continuity assumption for $\Phi$, we introduce the
concept of \emph{jets} of metrics at infinity, that is to say we consider only
a finite number of terms in the asymptotic expansion of $g-b$. We show that all
the facts holding above for germs (transversality, adjustment, group action)
descend naturally to jets. The set $J^l_k$ of jets of fixed order has a natural
topology so we define $\Phi$ to be continuous if it factors to a continuous map
$\overline{\Phi}: J^l_k \to V$ for some given $l \geq k$.

In the present work we will assume that $\Phi(g)$ only depends on the leading
coefficient in the asymptotic expansion of $g - b$: $g - b = \rho^{k-2} m +
    O(\rho^{k-1})$. This term will be called the \emph{mass aspect tensor} of $g$.
It is a symmetric $2$-tensor on $\bS^{n-1}$. In a forthcoming paper we will
study the situation when $\Phi(g)$ depends on further terms in the jet of the
asymptotically hyperbolic metric.

So we have reduced the problem of finding operators $\Phi$ mapping
asymptotically hyperbolic metrics to elements of $V$ and satisfying the
covariance property~\eqref{eqRoughMass} to classifying operators mapping
symmetric 2-tensors to elements of $V$ that intertwine the action of
$O_\uparrow(n, 1)$.

Before stating our main results, we introduce the following relevant
representations of $O_\uparrow(n, 1)$:
\begin{itemize}
    \item For any integer $l \geq 0$, let $\cH_l$ denote the space of wave-harmonic
          polynomials on $\bR^{n, 1}$ homogeneous of degree $l$:
          \[
              \cH_l \definedas \{ P \in \bR_l [X^0, \ldots, X^n],\ \Box P = 0\},
          \]
          where $\Box$ is the d'Alembertian operator: $\Box P = \eta^{\mu\nu}
              \partial_\mu \partial_\nu P$.
    \item Let also $\cW_l$ denote the set of polynomial Weyl tensors, i.e.\ covariant
          $4$-tensors $W$ over $\bR^{n, 1}$ whose entries are homogeneous polynomials of
          degree $l$ and satisfying the following properties:
          \begin{align*}
               & W_{\mu\nu\alpha\beta} = - W_{\nu\mu\alpha\beta},\ W_{\mu\nu\alpha\beta} = W_{\alpha\beta\mu\nu},                             \\
               & W_{\mu\nu\alpha\beta} + W_{\nu\alpha\mu\beta} + W_{\mu\alpha\beta\nu} = 0,                                                   \\
               & \partial_\alpha W_{\mu\nu\beta\gamma} + \partial_\beta W_{\mu\nu\gamma\alpha} + \partial_{\gamma} W_{\mu\nu\alpha\beta} = 0, \\
               & \eta^{\mu\nu} W_{\mu\nu\alpha\beta} = 0.
          \end{align*}
\end{itemize}
It turns out that all these representations are irreducible. However, when $n=3$, the representations $\bC \otimes \cW_l$ further decompose into two non-equivalent $SO_\uparrow(n, 1)$-irreducible representations. More details on these representations can be found in Appendix~\ref{secRepresentation}.

We now state the main result of our paper.
\begin{theorem}
    The linear masses $\Phi$ with values in an irreducible representation belong to one of the following two families:
    \begin{enumerate}
        \item The harmonic masses defined for any $k \geq n-1$ as follows: Let $g \in G_k$ be
              a transverse metric and denote by $m$ its mass aspect tensor. Then $\Phi_h(g)$
              is the linear form over $\cH_l$, with $l = k - n + 1$ defined by
              \[
                  P \mapsto \int_{\bS^{n-1}} P(1, x^1, \ldots, x^n) \tr^\sigma(m) d\mu^\sigma,
              \]
              where $\bS^{n-1}$ is identified with the section ${X^0=1}$ of the light cone
              and $\sigma$ is the round metric on $\bS^{n-1}$.
        \item The Weyl masses defined for any $k \geq n+1$ as follows: Let $g \in G_k$ be a
              transverse metric and denote by $m$ its mass aspect tensor. Then $\Phi_w(g)$ is
              the linear form over $\cW_l$, with $l = k - n - 1$ defined by
              \[
                  W \mapsto \int_{\bS^{n-1}} \left\<m, W(e_+, \cdot, e_+, \cdot) \right\> d\mu^\sigma,
              \]
              where $e_+ = \partial_0 + \sum_{i=1}^n x^i \partial_i$ is the outgoing null
              vector.
    \end{enumerate}
\end{theorem}
In dimension $n=3$, the Weyl masses can be further decomposed, see Proposition~\ref{thmWeylInvariants3d}. The classical Wang-Chru\'sciel-Herzlich mass vector is obtained here as the linear mass $\Phi_h$ for $k = n$, for which $\cH_1$ is the space of linear forms on $\bR^{n, 1}$. This theorem is contained in Propositions~\ref{thmScalarInvariants3d},~\ref{thmWeylInvariants3d} and~\ref{thmClassificationso(n,1)}. 

Note that the theorem could be generalized to include non-transverse metrics
as, by means of the adjustment diffeomorphism $\Theta$, we can replace any
arbitrary asymptotically hyperbolic metric $g$ by the transverse metric
$\Theta_* g$ without changing the value of $\Phi(g)$. The relation between the
leading terms in the asymptotic expansions of $g$ and $\Theta_* g$ is given
explicitly by Proposition~\ref{propAdjMass}.

\subsection{Outline of the paper}

In Section~\ref{secDefinition}, we introduce the class of asymptotically
hyperbolic metrics we are considering. We also study the germs and jets of
metrics at infinity and the transversality condition in relation to changes of
asymptotically hyperbolic coordinates. In particular, we introduce the notion
of linear mass at infinity. In Section~\ref{secLorentzaction}, we make explicit
the action of the Lie algebra of hyperbolic isometries on the space of mass
aspect tensors. In Section~\ref{secClassification}, we use the Lie algebra
action to classify the linear masses assuming that they depend only on the mass
aspect tensor. Appendix~\ref{secRepresentation} contains further information
about the spaces $\cH_l$ of harmonic polynomials and $\cW_l$ of polynomial Weyl
tensors.

\subsection*{Acknowledgments}

The authors wish to thank Lars Andersson, Xavier Bekaert, G\'erard Besson,
Rapha\"el Beuzart-Plessis, Piotr Chru\'sciel, Charles Frances, Gary Gibbons,
Marc Herzlich, C\'ecile Huneau, Eric Larsson, C\'edric Lecouvey, Vince
Moncrief, Karim Noui, and Philippe Roche for numerous fruitful discussions
concerning this work.

JC was supported by the ERC Advanced Grant 320939, Geometry and Topology of
Open Manifolds (GETOM) and warmly thanks the Institut Fourier, Grenoble, the
MPIM, Bonn and the FIM, ETH Z\"urich, for hospitality and financial support as
well as the EPDI program during part of this work. RG is partially supported by
the grant ANR-23-CE40-0010-02 of the French National Research Agency (ANR):
Einstein constraints: past, present, and future (EINSTEIN-PPF). JC and RG
warmly thank the MSRI, Berkeley, for the invitation to the fall 2013 program on
Mathematical relativity. JC and MD equally warmly thank the IHP, Paris, for the
invitation to the fall 2015 program on Mathematical general relativity.

Most of the algebraic calculations in this paper were done with the help of the
\textsc{Sage} software.

%% file: definition.tex

In this section, we give the basic definitions that will be used throughout the
paper.

\subsection{Hyperbolic space} \label{secHyperbolic}
The standard Riemannian manifold in the context of the paper is the hyperbolic
space $\bH^n$. Unless stated otherwise, we will use the Poincar\'e ball model.
In this model, $\bH^n$ is described as the unit ball $B_1(0)$ centered at the
origin in $\bR^n$ equipped with the metric
\[
    b \definedas \rho^{-2} \delta,
\]
where $\delta$ denotes the Euclidean metric on $\bR^n$ and the function $\rho$
is defined by
\[
    \rho(x) \definedas \frac{1-\left|x\right|^2}{2} .
\]
The hyperbolic distance is given by
\begin{equation} \label{eqHyperbolicDistance}
    \cosh d^b(x, y) = 1 + \frac{|x-y|^2}{2 \rho(x) \rho(y)}
\end{equation}
for $x, y \in B_1(0)$.

We briefly recall the isometry group of the hyperbolic space. Let $\bR^{n,1}$
denote Minkowski space, that is $\bR^{n+1}$ equipped with the quadratic form
\[
    \eta \definedas -(dX^0)^2 + \sum_{k=1}^n (dX^k)^2,
\]
where $X^0, X^1, \dots, X^n$ denote the standard coordinates on $\bR^{n+1}$.
Further, we denote by $\partial_0, \partial_1, \dots, \partial_n$ the standard
basis of $\bR^{n+1}$, which has the dual basis of one-forms $dX^0, dX^1, \dots,
    dX^n$. Hyperbolic space can be embedded into Minkowski space as the unit
hyperboloid, that is
\[
    \bH^n =
    \left\{
    (X^0, X^1, \dots, X^n) \in \bR^{n, 1}
    \middle| -(X^0)^2 + \sum_{k=1}^n (X^k)^2 = -1, X^0 > 0
    \right\}.
\]

The orthogonal group $O(n,1)$ consists of linear maps preserving the quadratic
form $\eta$. The isometry group of the hyperbolic space is the subgroup of
future preserving isometries called the \emph{orthochronous Lorentz group} and
is denoted by $O_{\uparrow}(n,1)$. The connected component of the identity
consists of orientation preserving isometries. It is called the
\emph{restricted Lorentz group} and is denoted by $SO_\uparrow(n,1)$.

The ball model of hyperbolic space is identified with the hyperboloid model via
a stereographic projection $p$ to the plane $X^0 = 0$ with respect to the point
$(-1, 0, \dots, 0)$,
\[
    p(X^0, X^1, \dots, X^n) \definedas \frac{1}{1+X^0} (X^1, \dots, X^n).
\]
We identify the plane $X^0 = 0$ with $\bR^n$, and we denote the standard
coordinates on $\bR^n$ by $x^1, \dots, x^n$. The inverse of $p$ is given by
\[
    p^{-1}(x^1, \dots, x^n)
    = \left(\frac{1+|x|^2}{1-|x|^2}, \frac{2 x^1}{1-|x|^2}, \dots,
    \frac{2 x^n}{1-|x|^2}\right),
\]
where $|x|^2 = (x^1)^2 + \cdots + (x^n)^2$.

If $A \in O_\uparrow(n,1)$ is an isometry of the hyperboloid, we transfer it to
an isometry of the ball model of $\bH^n$ which acts as $\Abar = p A p^{-1}$. In
what follows, we will mainly focus on the subgroup $SO_\uparrow(n, 1)$ and
consider two particular types of such elements.
\begin{itemize}
    \item
          A \emph{rotation} by an angle $\theta$ in the $X^i X^j$-plane ($1 \leq i,j \leq
              n$) is given by
          \[
              \begin{split}
                   & R_{ij}^\theta(X^0, \dots, X^i, \dots, X^j, \dots, X^n) \\
                   & \qquad =
                  (X^0, \dots, \cos(\theta) X^i - \sin(\theta) X^j, \dots,
                  \sin(\theta) X^i + \cos(\theta) X^j, \dots, X^n).
              \end{split}
          \]
          We denote the corresponding infinitesimal generator with a script letter,
          \[
              r_{ij} \definedas \frac{d}{d\theta} R_{ij}^\theta |_{\theta = 0}
              = dX^i \partial_j - dX^j \partial_i.
          \]
          Note that rotations commute with $p$, so $\Rbar \definedas p R p^{-1}$ reads
          \[
              \begin{split}
                   & \Rbar_{ij}^\theta
                  \definedas p R_{ij}^\theta p^{-1}(x^1, \dots, x^i, \dots, x^j, \dots, x^n) \\
                   & \qquad =
                  (x^1, \dots, \cos(\theta) x^i - \sin(\theta) x^j, \dots,
                  \sin(\theta) x^i + \cos(\theta) x^j, \dots, x^n).
              \end{split} \]
          The derivative of $p R_{ij}^\theta p^{-1}$ with respect to $\theta$ at
          $\theta=0$ is the rotation vector field $\fr_{ij}$, where
          \begin{equation}\label{Lij}
              \fr_{ij} \definedas x^i \partial_j - x^j \partial_i.
          \end{equation}
    \item
          A \emph{Lorentz boost} in the direction $X^i$ ($1 \leq i \leq n$) with a
          parameter $s \in \bR$ is given by
          \[ \begin{split}
                   & A_i^s(X^0, \dots, X^i, \dots, X^n) \\
                   & \qquad =
                  (\cosh(s) X^0 + \sinh(s) X^i, \dots, \sinh(s) X^0 + \cosh(s) X^i,
                  \dots, X^n).
              \end{split} \]
          The corresponding infinitesimal generator is given by
          \[
              a_i \definedas \frac{d}{ds} A_i^s |_{s = 0}
              = dX^0 \partial_i + dX^i \partial_0,
          \]
          and the corresponding isometry of the ball model is
          \begin{equation} \label{eqBoosts}
              \begin{split}
                   & \Abar^s_i \definedas p A^s_i p^{-1}(x^1, \dots, x^i, \dots, x^n) \\
                   & \qquad = \frac{1}{D}\left(x^1, \dots,
                  \cosh(s) x^i + \sinh(s) \frac{1+|x|^2}{2},
                  \dots, x^n\right),
              \end{split}
          \end{equation}
          where
          \[
              D \definedas
              \frac{1- |x|^2}{2} + \frac{1+|x|^2}{2} \cosh(s) + x^i \sinh(s).
          \]
          The derivative of $\Abar_i^s$ with respect to $s$ at $s=0$ is the boost vector
          field $\fa_i$, where
          \begin{equation}\label{Xi}
              \fa_i \definedas \frac{1 + |x|^2}{2} \partial_i - x^i x^a \partial_a.
          \end{equation}
\end{itemize}

In this article we will use the convention that upper case latin letters denote
elements in the Lie group $O_\uparrow(n,1)$, while lower case latin letters
will be used for elements in the Lie algebra $\mathfrak{so}(n,1)$ and fraktur
letters will denote the corresponding vector fields on the ball $B_1(0) \subset
    \bR^n$, or the same vector fields restricted to the unit sphere $\bS^{n-1}
    \subset \bR^n$. These vector fields are in fact tangent to the sphere, since
$\bS^{n-1}$ is preserved by elements of $O_\uparrow(n,1)$.

\subsection{Asymptotically hyperbolic metrics}\label{secAsymptoticallyHyp}

We continue by defining asymptotically hyperbolic manifolds. Several
definitions exist in the literature; we refer the reader
to~\cite{LeeFredholm,HerzlichMassFormulae,GicquaudCompactification} for an
overview. To avoid technical complications we choose here to use the simplest
such definition.

\begin{definition}\label{defAHMetric}
    Let $\Mbar$ be a compact manifold of dimension $n$ with boundary
    $\partial M$ diffeomorphic to $\bS^{n-1}$. We choose a neighborhood $\Omega$ of
    $\partial M$ in $\Mbar$, a diffeomorphism
    $\Psi: \Omega \to \overline{B}_1(0) \setminus \overline{B}_{1-\epsilon}(0)$
    from $\Omega$ to the standard annulus and a positive integer $k$. Then
    a metric $g$ on the interior $M$ of $\Mbar$ is said to be
    \emph{asymptotically hyperbolic} of order $k$ with respect to $\Psi$ if
    the metric
    \[
        \gbar \definedas \rho^2 \Psi_* g,
    \]
    which is a priori defined only on $B_1(0) \setminus
        \overline{B}_{1-\epsilon}(0)$, extends to a smooth metric on $\overline{B}_1(0)
        \setminus \overline{B}_{1-\epsilon}(0)$ such that
    \[
        \left|\gbar - \delta\right|_{\delta} = O\left( \rho^k \right).
    \]
    The boundary $\partial M \simeq \bS^{n-1}$ is called the \emph{sphere at
        infinity} for the asymptotically hyperbolic metric $g$.

    Here, $\rho$ refers to the \emph{boundary defining function} introduced in
    paragraph~\ref{secHyperbolic}.
\end{definition}

If the metric $g$ is asymptotically hyperbolic of order $k$ with respect to two
diffeomorphisms $\Psi_i : \Omega_i \to \overline{B}_1(0) \setminus
    \overline{B}_{1-\epsilon_i}(0)$ for $i = 1,2$, then the composition $\Psi_2
    \circ \Psi_1 ^{-1}$ is a diffeomorphism between two neighborhoods of
$\bS^{n-1}$ in $\overline{B}_1(0)$. Hence we will focus on the study of
asymptotically hyperbolic metrics in neighborhoods of the sphere at infinity
$\bS^{n-1}$ in $\overline{B}_1(0)$. Our results will however apply for any
manifold $(M,g)$ which is asymptotically hyperbolic with respect to some
diffeomorphism $\Psi$ as in Definition~\ref{defAHMetric}.

Let $r = r(x)$ denote the distance from the origin in the ball model of
$\bH^n$. From Formula~\eqref{eqHyperbolicDistance} we have
\[
    \rho(x) = \frac{1}{\cosh(r)+1}.
\]
As a consequence, the estimate for the decay of the metric can be rewritten in
a more intrinsic way as
\begin{equation} \label{eqEstimateDecay}
    \left|g - b\right|_b = O\left( e^{-k r} \right),
\end{equation}
and a simple argument using the triangle inequality shows that replacing
$r$ by the distance function from any given point in $B_1(0)$ gives
an equivalent decay estimate.

Since we are interested in the asymptotic behavior of asymptotically hyperbolic
metrics, it is sufficient to restrict our attention to germs of such metrics.
An open subset $U \subset B_1(0)$ is called a \emph{neighborhood of infinity}
if $B_1(0) \setminus \overline{B}_{1-\epsilon}(0) \subset U$ for some $\epsilon
    > 0$. We denote by $\cN_\infty$ the set of neighborhoods of infinity. Given an
element $U \in \cN_\infty$ we define
\[
    \Util \definedas U \cup \bS^{n-1},
\]
that is $\Util$ is the set of points in $U$ together with the sphere at
infinity.

\begin{definition}\label{defAHGerm}
    For a positive integer $k$, we consider the set $G_k^0$ consisting of
    pairs $(U, g)$ where $U$ is a neighborhood of infinity and $g$ is a
    metric on $U$ such that $\gbar \definedas \rho^2 g$ extends to a smooth
    metric on $\Util$ satisfying
    \[
        \left|g - b\right|_b
        = \left|\gbar - \delta\right|_\delta
        = O\left(\rho^k\right).
    \]
    The relation $\sim$ is defined by $(U_1, g_1) \sim (U_2, g_2)$ if there is a
    $U_3 \in \cN_\infty$ such that $U_3 \subset U_1 \cap U_2$ and $g_1 \equiv g_2$
    on $U_3$. This defines an equivalence relation on $G_k^0$, and we define the
    \emph{stalk at infinity of asymptotically hyperbolic metrics of order $k$} as
    $G_k \definedas G^0_k / \sim$. An element $g \in G_k$ will be called an
    \emph{asymptotically hyperbolic germ at infinity}.
\end{definition}

We will, by abuse of notation, identify $g \in G_k$ with the metric of a
representative $(U, g)$. This terminology is modeled on standard terminology of
sheaf theory, see for example~\cite[Chapter~1]{BredonSheafTheory}
or~\cite[Chapter~2]{Hartshorne}. Note that the stalk at infinity of
asymptotically hyperbolic metrics of order $k$ is an affine space.

We will now introduce the stalk at infinity of asymptotic isometries of an
asymptotically hyperbolic metric.

\begin{definition}\label{defASIsom}
    Given a positive integer $k$ and $g \in G_k$, we define the set
    $I^k(g)$ of \emph{asymptotic isometries} of $g$ as the stalk at
    infinity of diffeomorphisms $\Psi: U \to V$ where $U$ and $V$ are
    neighborhoods of infinity such that $\Psi_* g \in G_k$.
\end{definition}

In Lemma~\ref{lmASIsom} it will be proven that $I^k(g)$ does not depend on the
choice of $g \in G_k$. As a corollary, we will conclude that $I^k(g)$ is a
group under composition of maps.

We also introduce a particular class of germs of asymptotically hyperbolic
metrics. This class will play an important role in what follows.

\begin{definition}\label{defTransverseGerm}
    Let $g$ be a germ of an asymptotically hyperbolic metric of order $k$.
    We say that $g$ is \emph{transverse} if there exists $(U, g)$
    representing $g$ such that
    \[
        g_{ij} x^i = b_{ij} x^i
    \]
    on $U$, or equivalently,
    \[
        g(P, \cdot) = b(P, \cdot)
    \]
    on $U$, where $P \definedas x^i \partial_i$ is the position vector field on
    $\bR^n$. We denote by $G^T_k$ the set of transverse asymptotically hyperbolic
    germs of order $k$.
\end{definition}


In this definition, and in the rest of the paper, we apply the convention that
an expression with an index appearing both as upper and lower index is summed
over the appropriate range for that index.

Note that the transversality condition can be stated in terms of the dual of
the metric as
\[
    g(d|x|^2, \cdot) = b(d|x|^2, \cdot)
\]
on $U$, or equivalently in terms of the covectors $dx^i$ as
\begin{equation} \label{eqTransversalityX}
    \gbar(d\rho, dx^i) = - x^i
\end{equation}
on $U$. In particular,
\begin{align*}
    g(d|x|^2, d|x|^2)   & = b(d|x|^2, d|x|^2),            \\
    g(d\rho, d\rho)     & = b(d\rho, d\rho),              \\
    \gbar(d\rho, d\rho) & = \delta(d\rho, d\rho) = |x|^2,
\end{align*}
and
\begin{equation} \label{eqTransversalityRho}
    |d\rho|^2_{\gbar} = 1 - 2 \rho.
\end{equation}

Note that, in the ball model, the function $\rho$ defined earlier is not a
\emph{geodesic defining function} as
in~\cite[Lemme~2.1.2]{DjadliGuillarmouHerzlich} and in Wang's approach to
defining the mass, see~\cite{WangMass}, nor an \emph{adapted defining function}
as in~\cite{CapGover}.

Next we define diffeomorphisms asymptotic to the identity.

\begin{definition}\label{defASIdent}
    Given a positive integer $k$ and a neighborhood of infinity $U$, we say
    that a diffeomorphism $\Theta: \Util \to \Theta(\Util)$ is asymptotic to
    the identity of order $k$ if
    \[
        \left| \Theta(x) - x \right|_\delta = O(\rho^k).
    \]
    We denote by $I_0^k$ the stalk at infinity of diffeomorphisms asymptotic to the
    identity. The stalk $I_0^k$ is a group under composition.
\end{definition}

The next lemma tells us that a diffeomorphism asymptotic to the identity of
order $k+1$ is an asymptotic isometry of order $k$ for any germ of metrics of
order $k$.

\begin{lemma}\label{lmI0inI}
    For any $g \in G_k$ we have $I_0^{k+1} \subset I^k(g)$.
\end{lemma}

\begin{proof}
    Since the stalk $I^{k+1}_0$ is a group under composition, it suffices to check that
    $\Theta^* g \in G_k$ for all $\Theta \in I^{k+1}_0$ since this is
    equivalent to $(\Theta^{-1})_* g \in G_k$. In component notation we have
    \[\begin{split}
            \rho^2(x) (\Theta^* g)(x)_{ij} - \delta_{ij}
             & =
            \rho^2(x) g(\Theta(x))_{kl}
            \partial_i \Theta^k \partial_j \Theta^l - \delta_{ij} \\
             & =
            \frac{\rho^2(x)}{\rho^2(\Theta(x))} \gbar(\Theta(x))_{kl}
            \partial_i \Theta^k \partial_j \Theta^l - \delta_{ij}.
        \end{split}\]
    A standard argument shows that the function $\frac{\rho \circ \Theta}{\rho}$ is
    smooth near $\bS^{n-1}$. Write
    \[ \begin{split}
            \Theta^i(x)
             & =
            \Theta^i\left(\frac{x}{|x|}\right)
            - \int_{|x|}^1 \frac{d}{ds} \Theta^i\left(s \frac{x}{|x|}\right) ds     \\
             & =
            \frac{x^i}{|x|}
            - \int_{|x|}^1 \partial_j\Theta^i\left(s \frac{x}{|x|}\right)
            \frac{x^j}{|x|} ds                                                      \\
             & =
            \frac{x^i}{|x|} - (1 - |x|) \int_{0}^1
            \partial_j\Theta^i\left((\lambda + (1-\lambda) |x|) \frac{x}{|x|}\right)
            \frac{x^j}{|x|} d\lambda                                                \\
             & =
            \frac{x^i}{|x|} - (1 - |x|) \left(\frac{x^i}{|x|}
            + \int_{0}^1 \left(\partial_j\Theta^i\left((\lambda + (1-\lambda) |x|)
            \frac{x}{|x|}\right) - \delta_j^i\right)\frac{x^j}{|x|} d\lambda\right) \\
             & =
            x^i - (1 - |x|)
            \underbrace{
                \frac{x^j}{|x|} \int_{0}^1 \left(\partial_j\Theta^i\left(\lambda \frac{x}{|x|}
                + (1-\lambda) x\right) - \delta_j^i\right) d\lambda}_{\eqqcolon E^i(x)},
        \end{split} \]
    where we set $s = \lambda + (1-\lambda) |x|$. The vector field $E^i$ is smooth
    near $\bS^{n-1}$, with $E^i(x) = O(\rho^k)$. Thus,
    \[ \begin{split}
            \frac{\rho \circ \Theta}{\rho}(x)
             & =
            \frac{1 - |\Theta(x)|^2}{1 - |x|^2}  \\
             & =
            1 + 2 \frac{(1-|x|) x^i E^i(x)}{1-|x|^2}
            - \frac{(1-|x|)^2 |E(x)|^2}{1-|x|^2} \\
             & =
            1 + 2 \frac{x^i E^i(x)}{1+|x|} - \frac{1-|x|}{1+|x|}|E(x)|^2.
        \end{split} \]
    From this we conclude that $\frac{\rho \circ \Theta}{\rho}$ is smooth and
    satisfies
    \[
        \frac{\rho \circ \Theta}{\rho} = 1 + O(\rho^k).
    \]
    A straightforward computation shows that $\gbar(\Theta(x))_{kl} \partial_i
        \Theta^k \partial_j \Theta^l - \delta_{ij} = O(\rho^k)$. The lemma follows by
    multiplying these estimates.
\end{proof}

The following proposition is a variant of~\cite[Lemma~5.2]{GrahamLee}, see
also~\cite[Lemma~5.1]{LeeSpectrum}. It states that any germ of metrics of order
$k$ can be made transverse by a unique diffeomorphism asymptotic to the
identity of order $k+1$.

\begin{proposition}\label{propTransverse}
    Given an element $g \in G_k$ there exists a unique $\Theta \in I^{k+1}_0$
    such that $\Theta_* g$ is transverse. Further, we have $\Theta_* g \in G_k$.
    In particular, if $g \in G_k$ is transverse, the only $\Theta \in I_0^{k+1}$
    such that $\Theta_* g$ is again transverse is the identity.
\end{proposition}

The diffeomorphism $\Theta$ provided by this proposition is called the
\emph{adjustment diffeomorphism} for the metric $g$.

\begin{proof}
    We construct new coordinates at infinity that satisfy the transversality
    condition for the metric $g$. These new coordinates, which
    we denote by $(x'_1, \dots, x'_n)$, will provide the required
    diffeomorphism $\Theta$ through the relation
    \[
        \Theta(x_1, \dots, x_n) = (x'_1, \dots, x'_n).
    \]
    We first set
    \[
        \rho' = \frac{1-|x'|^2}{2}.
    \]
    We seek $\rho'$ to be close to $\rho$ (in a sense to be specified below). We
    set $\rho' = \rho e^v$ for a function $v$ to be determined. We define $\gtil
        \definedas \rho'^2 g = e^{2v} \gbar$. We first impose an analogue of
    Condition~\eqref{eqTransversalityRho} for $\gtil$ and $\rho'$, that is
    \[
        |d\rho'|^2_{\gtil} = 1-2\rho'.
    \]
    We rewrite this as an equation for $v$,
    \[ \begin{split}
            1-2\rho e^v
             & =
            e^{-2v} |\rho e^v dv + e^v d\rho|^2_{\gbar} \\
             & =
            |\rho dv + d\rho|^2_{\gbar}                 \\
             & =
            \rho^2 |dv|^2 + 2 \rho \gbar(dv, d\rho) + |d\rho|_{\gbar}^2,
        \end{split} \]
    or
    \begin{equation} \label{eqTransversalityV}
        2 \gbar(dv, d\rho) + \rho |dv|^2
        = \frac{1 - 2\rho - |d\rho|_{\gbar}^2}{\rho} + 2 - 2 e^v.
    \end{equation}
    Since $|\gbar-\delta|_{\delta} = O(\rho^k)$, we have
    \[
        |d\rho|_{\gbar}^2 = 1 - 2 \rho + O(\rho^k),
    \]
    meaning that
    \begin{equation} \label{eqErrorTermRho}
        \frac{1 - 2\rho - |d\rho|_{\gbar}^2}{\rho} = O(\rho^{k-1}).
    \end{equation}
    Equation~\eqref{eqTransversalityV} is a first order partial differential
    equation for $v$. The relevant theory for such equations can be found in
    ~\cite[Chapter~2]{CourantHilbert2} or in
    ~\cite[Theorem~22.39]{LeeSmoothManifolds}.
    In particular, the condition $|d\rho|^2_{\gbar} \equiv 1$ on $\bS^{n-1}$ ensures
    that there exists a unique solution $v$ in a neighborhood of $\bS^{n-1}$ such
    that $v = 0$ on $ \bS^{n-1}$. From Estimate~\eqref{eqErrorTermRho} it follows that $v = O(\rho^k)$. We have now determined the function $\rho' = e^v \rho$.

    We introduce the analog of Equation~\eqref{eqTransversalityX} for the
    coordinates $x'^i$,
    \[
        \gtil(d\rho', dx'^i) = - x'^i,
    \]
    with the boundary condition $x'^i = x^i$ on $\bS^{n-1}$. As for the previous
    equation, existence and uniqueness of a smooth solution in a neighborhood of
    $\bS^{n-1}$ is guaranteed by classical results. A straightforward calculation
    yields
    \[
        \gtil(d\rho', dx^i) = - x^i + O(\rho'^k).
    \]
    This implies that $x'^i-x^i = O(\rho'^{k+1})$ which means that $\Theta \in
        I^{k+1}_0$. It remains to verify that $\rho' = \frac{1 - |x'|^2}{2}$. Letting
    $\rho'' = \frac{1 - |x'|^2}{2}$ we have
    \[ \begin{split}
            \gtil(d\rho', d\rho'' )
             & = -x'^i \gtil(d\rho', dx'^i) \\
             & = \sum_i (x'^i)^2            \\
             & = 1 - 2 \rho''.
        \end{split} \]
    This is a non-characteristic first order PDE for the function $\rho''$. Since
    $x'^i$ coincides with $x^i$ on $\bS^{n-1}$ we have that $\rho'' = 0 = \rho'$ on
    $\bS^{n-1}$. It follows that $\rho'' = \rho'$ in a neighborhood of $\bS^{n-1}$.

    Uniqueness is easy to prove. Assume that we have two such diffeomorphisms
    $\Theta$ and $\Theta'$. Then, considering $\Psi = \Theta' \circ \Theta^{-1}$,
    we have to prove that for a transverse metric $g$ the only element $\Psi \in
        I^{k+1}_0$ such that $\Psi_* g$ is also transverse is the identity. Assume that
    the diffeomorphism $\Psi$ is given in coordinates by
    \[
        \Psi(x^1, \dots, x^n) = (x'^1, \dots, x'^n),
    \]
    where $x'^i - x^i = o(\rho)$. As before we set
    \[
        \rho = \frac{1-|x|^2}{2}
        \quad \text{and} \quad
        \rho' = \frac{1-|x'|^2}{2}.
    \]
    From the assumption on the coordinates we have that $\rho' = \rho + o(\rho)$.
    Computing as for Equation~\eqref{eqTransversalityV}, we get that $v \definedas
        \log \frac{\rho'}{\rho}$ vanishes on $\bS^{n-1}$ and satisfies
    \[
        2 \<dv, d\rho\>_{\gbar} + \rho |dv|^2 = 2 - 2 e^v,
    \]
    where we have used the fact that $|d\rho|_{\gbar}^2 = 1 - 2\rho$ since $g$ is
    transverse. Since the solution to this equation is unique, we must have $v
        \equiv 0$, or equivalently, $\rho \equiv \rho'$. From
    Condition~\eqref{eqTransversalityX} we deduce that the coordinates $x^i$ and
    $x'^i$ both satisfy the transport equation
    \[
        \left\<d\rho, dx^i\right\>_{\gbar} = - x^i.
    \]
    Since they coincide on $\bS^{n-1}$ they coincide in a neighborhood of
    $\bS^{n-1}$.

    The proof that $\Psi_* g \in G_k$ is contained in Lemma~\ref{lmI0inI} above.
\end{proof}

\begin{remarks}
    \begin{itemize}
        \item
              We note that in all previous proofs of the existence of geodesic normal coordinates
              (see for example~\cite[Lemma~5.1]{LeeSpectrum}) there is a loss of
              regularity of one derivative. Avoiding this is one of the reasons why we chose
              to restrict our study to smooth conformally compact manifolds. However, the
              proof of Proposition~\ref{propAdjustmentJet} indicates that there should be no
              such loss.
        \item The dependence of the adjustment diffeomorphism $\Theta$ on the germ $g$ will
              be studied in Subsection~\ref{secGeomMasses}. The first non-trivial term in the
              asymptotic expansion of $\Theta$ will be computed in
              Proposition~\ref{propAdjMass}.
    \end{itemize}
\end{remarks}

Before stating the main result of this section, we prove two important lemmas.
The first lemma is a rephrasing of~\cite[Theorem~6.1]{ChruscielHerzlich}. It
states that the germ of an asymptotic isometry can be extended as a
diffeomorphism to the sphere at infinity. The proof we give is based on the
theory developed in~\cite{BahuaudMarsh,Bahuaud,BahuaudGicquaud,
    GicquaudCompactification}. In order to keep this section reasonably short, we
refer the reader to~\cite{BahuaudGicquaud} for the definition of essential sets
and the relevant results concerning them.

\begin{lemma}\label{lmExtension}
    For $k$ a positive integer let $g \in G_k$ be an asymptotically
    hyperbolic metric and let $\Psi \in I^k(g)$ be an asymptotic isometry.
    If $\psi: U \to V$ represents the germ $\Psi$, then $\psi$ extends to
    a smooth diffeomorphism $\psibar: \Util \to \Vtil$.
\end{lemma}

\begin{proof}
    By Proposition~\ref{propTransverse}, we may pull back the metrics $g$ and
    $\Psi^*g$ via elements in $I_0^{k+1}$ so that they satisfy the transversality
    condition. Since this corresponds to composing the diffeomorphism $\psi$ on the
    left and on the right with elements of $I_0^{k+1}$ which are smooth up to the
    boundary, we can assume without loss of generality that $g$ and $\Psi^*g$ are
    transverse.

    From~\cite[Lemma~2.5.11]{GicquaudThesis} we know that the set $K = \{\rho \geq
        \epsilon\}$ is an essential set for both $g$ and $\psi^*g$ provided $\epsilon >
        0$ is small enough. \footnote{The notion of an essential set is a priori
        defined only for complete manifolds. Here we can simply extend $g$ (resp.
        $\Psi^* g$) to all of $B_1(0)$ in an arbitrary way. The argument
        of~\cite[Lemma~2.5.11]{GicquaudThesis} depends only on the metric outside some
        compact set.} Equivalently, $K_1 \definedas K$ and $K_2 \definedas \psi(K)$ are
    essential sets for $g$. Further, since the metrics $g$ and $\psi^* g$ are
    smoothly conformally compactifiable, their sectional curvatures satisfy
    \begin{equation} \label{eqAsymptoticCurvature}
        \sec^g = -1 + O(\rho)
        \quad \text{and} \quad
        \sec^{\psi^*g} = -1 + O(\rho),
    \end{equation}
    see~\cite[Section 1.2.1]{GicquaudThesis}. Using the conformal transformation
    law for the curvature tensor $\riem$, see~\cite[Proposition 1.1.2]{GicquaudStability},
    together with~\cite[Lemma 2.7]{LeeFredholm}, one can prove the following estimates:
    \begin{equation} \label{eqAsymptoticCurvature2}
        \left|\nabla^g \riem^g\right|_g = O(\rho)
        \quad \text{and} \quad
        \left|\nabla^{\psi^*g} \riem^{\psi^*g}\right|_{\psi^*g} = O(\rho).
    \end{equation}

    The transversality condition imposes that the distance from $K_1$ (resp. $K_2$)
    with respect to the background hyperbolic metric and $g$ (resp. $\psi^* g$)
    agree. For points $y$ in the boundary of $\{\rho \geq \epsilon\}$, that is such
    that $\rho(y) = \epsilon$, we have $|y| = \sqrt{1-2\epsilon}$. We also remark
    that the closest point projections $\pi(x)$ of a point $x$ onto $\Sigma =
        \rho^{-1}(\{\epsilon \})$ with respect to the metrics $b$, $g$ and $\psi^* g$
    all coincide with the Euclidean closest point projection onto $\Sigma$ because
    $\Sigma$ is a round sphere centered at the origin. Hence,
    \[
        |x - \pi(x)|^2 = \left(|x| - \sqrt{1-2\epsilon}\right)^2.
    \]
    From Equation~\eqref{eqHyperbolicDistance}, the distance from $K$ to any point
    $x$ lying outside $K$ is given by
    \[
        \cosh d(x, K) =
        1 + \frac{\left(|x| - \sqrt{1-2\epsilon}\right)^2}{2\epsilon \rho(x)}.
    \]
    A straightforward calculation shows that $e^{-d(x, K)}$ can be expressed as an
    analytic function of $\rho$ such that
    \begin{equation}\label{eqAsymptoticDistance}
        e^{-d(x, K)} \sim
        \frac{\epsilon}{\left(1 - \sqrt{1-2\epsilon}\right)^2} \rho(x).
    \end{equation}
    Thus we can replace the conformal factor $e^{-d(x, K)}$ in
    ~\cite{BahuaudGicquaud} by $\rho$ and the results of this article
    still apply. It follows from Equations~\eqref{eqAsymptoticCurvature},
    ~\eqref{eqAsymptoticCurvature2}, and~\eqref{eqAsymptoticDistance} that
    both $g$ and $\Psi^* g$ fulfill the conditions of
    ~\cite[Theorem~A]{BahuaudGicquaud} with $a = 1$. In particular, since the
    $C^{1, \alpha}$-structure ($0 < \alpha < 1$) of the manifold with boundary
    obtained by geodesic conformal compactification is unique, we conclude
    that $\psi$ extends to a $C^{1, \alpha}$ diffeomorphism
    $\psibar: \Util \to \Vtil$.

    The proof is completed using a standard bootstrapping argument involving
    Christoffel symbols. Let $\Gambar$ denote the Christoffel symbols of the metric
    $\gbar = \rho^2 g$ in the coordinate system $(x^1, \dots, x^n)$. Similarly, let
    $\Gamtil$ denote the Christoffel symbols of the metric $\gtil \definedas \rho^2
        \psi^* g$. By assumption, all the components of $\Gambar$ and $\Gamtil$ are
    smooth functions. The transformation law for the Christoffel symbols reads
    \[
        \frac{\partial^2 \psi^k}{\partial x^i \partial x^j}
        = \frac{\partial \psi^l}{\partial x^i}
        \frac{\partial \psi^m}{\partial x^j} \Gambar^k_{lm}
        - \Gamtil^l_{ij} \frac{\partial \psi^k}{\partial x^l}.
    \]
    Since $\psi$ is a $C^1$ diffeomorphism up to the boundary, the previous formula
    immediately shows that $\psi$ is actually a $C^2$ diffeomorphism since the
    right hand side only involves first order derivatives of $\psi$ (hence $C^0$
    functions) together with $C^\infty$ functions (the Christoffel symbols). The
    process can be iterated to conclude that $\psi$ is actually smooth up to the
    boundary.
\end{proof}

The second lemma states that any asymptotic isometry can be written as a
composition of a true isometry of the hyperbolic metric and a diffeomorphism
asymptotic to the identity. This is a variant of a result by Chru\'sciel and
Nagy~\cite{ChruscielNagy}, and Chru\'sciel and
Herzlich~\cite{ChruscielHerzlich}. A simplified proof (of a weaker result) can
be found in~\cite{HerzlichMassFormulae}.

\begin{lemma}\label{lmASIsom}
    Given a germ of metrics $g \in G_k$, any element $\Psi \in I^k(g)$
    decomposes uniquely as $\Psi = \Abar \circ \Psi_0$ (resp.
    $\Psi = \Psi'_0 \circ \Abar$), where $A \in O_{\uparrow}(n,1)$ and
    $\Psi_0 \in I^{k+1}_0$ (resp. $\Psi'_0 \in I^{k+1}_0$). Conversely,
    any element of the form $\Abar \circ \Psi_0$
    (resp. $\Psi'_0 \circ \Abar$), where $\Psi_0 \in I^{k+1}_0$
    (resp. $\Psi'_0 \in I^{k+1}_0$) and $A \in O_{\uparrow}(n,1)$, belongs
    to $I^k(g)$. In particular, the set $I^k(g)$ does not depend on the
    choice of $g \in G_k$.
\end{lemma}

\begin{proof}
    We first prove that if $A$ is an isometry of the metric $b$ then $A$
    belongs to $I^k(g)$ for any metric $g \in G_k$. Any such element
    can be written as the composition of an element of $O(n)$ and a Lorentz
    boost $\Abar_i^s = p A_i^s p^{-1}$. Hence, from the formulas given in
    Subsection~\ref{secHyperbolic}, we immediately see that $A$ is actually
    a smooth diffeomorphism of $\overline{B}_1(0)$. The metric
    $\rho^2 \Abar_* g$ can be rewritten as
    \[
        \gtil =
        \frac{\rho^2}{\left(\rho \circ \Abar^{-1}\right)^2} \Abar_* \gbar.
    \]
    Arguing as in the proof of Lemma~\ref{lmI0inI}, we have that the function
    $\frac{\rho^2}{\left(\rho \circ \Abar^{-1}\right)^2}$ is smooth. Since the
    metric $\Abar_* \gbar$ is smooth we conclude that $\gtil$ is a smooth metric on
    some $\Util$, $U \in \cN_\infty$. The decay estimate~\eqref{eqEstimateDecay}
    for the metric $\Abar_* g$, that is $\left|\Abar_* g - b\right|_b = O(\rho^k)$,
    is readily verified since $r(\Abar(x)) = d^b(\Abar(x), 0) = d^b(x,
        \Abar^{-1}(0))$ (see the remark following Equation~\eqref{eqEstimateDecay}).

    Thus, for any $g \in G_k$, given $A \in O_{\uparrow}(n, 1)$ and $\Psi_0 \in
        I_0^{k+1}$, we have $(\Abar \circ \Psi_0)_* g \in G_k$. A similar argument
    shows that $(\Psi_0 \circ \Abar)_* g \in G_k$.

    Next, given $\Psi \in I^k(g)$, we seek the corresponding isometry $\Abar$. It
    follows from Lemma~\ref{lmExtension} that $\Psi$ extends to a smooth
    diffeomorphism up to the boundary $\bS^{n-1}$. We claim that $\Psi$ induces a
    conformal diffeomorphism on $\bS^{n-1}$. The metric $\gtil \definedas \rho^2
        \Psi_* g$ can be rewritten as
    \[
        \gtil = \frac{\rho^2}{\left(\rho \circ \Psi^{-1}\right)^2} \Psi_* \gbar.
    \]
    Arguing once again as in the proof of Lemma~\ref{lmI0inI}, we have that the
    function $\frac{\rho^2}{\left(\rho \circ \Psi^{-1}\right)^2}$ is smooth on some
    $\Util$ where $U \in \cN_\infty$. Restricting to $\bS^{n-1}$, since both
    $\gbar$ and $\gtil$ restrict to the round metric $\sigma$, $\Psi$ induces a
    conformal isometry. The restriction of $\Psi$ to $\bS^{n-1}$ coincides with the
    restriction of a unique isometry $\Abar$ of the ball model of hyperbolic space.

    Considering $\Abar^{-1} \circ \Psi$ (resp. $\Psi \circ \Abar^{-1}$) we are left
    with proving that an element $\Psi \in I^k(g)$ such that $\Psi$ induces the
    identity on $\bS^{n-1}$ belongs to $I_0^{k+1}$. By composing $\Psi$ on the left
    and right with elements of $I_0^{k+1}$, we may assume that both $g$ and $\Psi_*
        g$ are transverse. The lemma will follow if, under these assumptions, $\Psi$ is
    the identity.

    Following the proof of Proposition~\ref{propTransverse}, it suffices to show
    that the ratio $\frac{\rho \circ \Psi}{\rho}$ tends to $1$ at infinity. We set
    $\Xi \definedas \Psi^{-1}$ and write
    \[
        \Xi^j(x) = x^j + \rho \xi^j(x)
    \]
    together with $\gbar = \rho^2 g = \delta + \ebar$ where $|\ebar|_\delta =
        O(\rho^k)$. Note that
    \[
        \partial_i \Xi^k = \delta_i^k + (\partial_i \rho) \xi^k
        + \rho \partial_i \xi^k
        = \delta_i^k - x^i \xi^k + O(\rho).
    \]
    The condition $\Psi_* g = \Xi^* g \in G_k$ can be written in coordinates as
    \begin{align*}
        \rho^{-2}(\Xi(x))\left(\delta_{kl}
        + \ebar_{kl}(\Xi(x))\right) \partial_i \xi^k \partial_j \xi^l
         & = \rho^{-2} \delta_{ij} + O(\rho^{k-2}),               \\
        \delta_{kl} \left(\delta^k_i - x^i \xi^k\right)
        \left(\delta^l_j - x^j \xi^l\right)
         & = \frac{\rho^2(\Xi(x))}{\rho^2} \delta_{ij} + O(\rho)  \\
        \delta_{ij} - x^i \xi^j - x^j \xi^i + x^i x^j |\xi|^2
         & = \frac{\rho^2(\Xi(x))}{\rho^2} \delta_{ij} + O(\rho).
    \end{align*}
    The limit of $\frac{\rho \circ \Psi}{\rho}$ on $\bS^{n-1}$ is then
    obtained by considering the last equality contracted twice with any
    vector $V$ orthogonal to the position vector field $P$, that is
    $V^i x^i = 0$,
    \[
        \delta_{ij}V^i V^j
        = \frac{\rho^2(\Xi(x))}{\rho^2} \delta_{ij} V^i V^j + O(\rho).
    \]
    Hence,
    \[
        \left.\frac{\rho \circ \Psi}{\rho}\right\vert_{\bS^{n-1}} = 1.
    \]
    This concludes the proof of the lemma.
\end{proof}

One of the consequences of Lemma~\ref{lmASIsom} is that the set $I^k(g)$ is
independent of $g \in G_k$ since any of its elements can be written as $\Abar
    \circ \Psi_0$ where $\Abar$ is an isometry of $b$ and $\Psi_0 \in I^{k+1}_0$.
It can also be seen that $I^k(g)$ is actually a group under composition.
Indeed, since $I^k(g)$ is independent of $g \in G_k$, given $\Psi_1, \Psi_2 \in
    I^k(g)$, we have $(\Psi_2)_* g \in G_k$ and
\[
    (\Psi_1 \circ \Psi_2)_* g = (\Psi_1)_* ((\Psi_2)_* g) \in G_k,
\]
thus proving that $I^k(g)$ is closed under composition.

The relationship between the groups of isometries, asymptotic isometries, and
diffeomorphisms asymptotic to the identity is summarized in the following
theorem.

\begin{theorem}\label{thmChruscielNagy}
    There is a short exact sequence of groups
    \begin{equation}\label{eqChruscielNagy}
        1 \longrightarrow I^{k+1}_0
        \overset{i}\longrightarrow I^k(g)
        \overset{\pi}{\longrightarrow} O_\uparrow(n, 1)
        \longrightarrow 1.
    \end{equation}
\end{theorem}

\begin{proof}
    It follows from Lemma~\ref{lmI0inI} that we have an inclusion
    $i$ of $I_0^{k+1}$ in $I^k(g)$. Given any element $\Theta \in I^{k+1}_0$
    and any $\Psi \in I^k(g)$, Lemma~\ref{lmASIsom} gives us a decomposition
    $\Psi = \Abar \circ \Psi_0$ where $\Psi_0 \in I^{k+1}_0$ and
    $A \in O_\uparrow(n,1)$. Then
    \[
        \Psi \circ \Theta \circ \Psi^{-1}
        = \Abar \circ \Psi_0 \circ \Theta \circ \Psi_0^{-1} \circ \Abar^{-1}.
    \]
    Note that $\Psi_0 \circ \Theta \circ \Psi_0^{-1} \in I^{k+1}_0$ so there exists
    an element $\Theta_1 \in I^{k+1}_0$ such that $\Abar \circ \Psi_0 \circ \Theta
        \circ \Psi_0^{-1} = \Theta_1 \circ \Abar$. Consequently,
    \[
        \Psi \circ \Theta \circ \Psi^{-1}
        = \Theta_1 \circ \Abar \circ \Abar^{-1}
        = \Theta_1 \in I^{k+1}_0.
    \]
    This proves that $I^{k+1}_0$ is a normal subgroup of $I^{k}(g)$, and the
    quotient $I^k(g)/I^{k+1}_0$ is therefore identified with $O_\uparrow(n,1)$.
\end{proof}

The sequence~\eqref{eqChruscielNagy} actually splits which can be seen by
taking $g = b$ (and recalling that $I^k(g)$ does not depend on $g$). The
natural action of $O_\uparrow(n, 1)$ on $\bH^n$ gives an embedding of
$O_\uparrow(n, 1)$ into $I^k(b)$ which is a right inverse for $\pi$.

Theorem~\ref{thmChruscielNagy} allows us to construct a natural action of
$O_\uparrow(n, 1)$ on the set $G^T_k$ of transverse germs of metrics.

\begin{proposition}\label{propActionLorentz}
    There exists a unique action of $O_\uparrow(n, 1)$ on the set $G^T_k$ such
    that, for any $A \in O_\uparrow(n, 1)$ and $g \in G^T_k$, we have
    \[
        A \cdot g = \widetilde{A}_* g,
    \]
    where $\widetilde{A}$ is the unique element in $I^k(g)$ such that
    $\widetilde{A}_* g$ is transverse and $\pi(\widetilde{A}) = A$.
\end{proposition}

\begin{proof}
    Let $g \in G_k^T$ and $A \in O_\uparrow(n, 1)$. From Proposition
    \ref{propTransverse}, there exists a unique $\Theta \in I_0^{k+1}$ such
    that the metric $\Theta_* \Abar_* g$ is transverse. We set
    $A \cdot g \definedas \Theta_* \Abar_* g$. Note that
    $A_g \definedas \Theta\circ \Abar$ is the unique element of $I^k(g)$
    such that $\pi(A_g) = A$ and $(A_g)_* g \in G^T_k$.

    We check that this defines an action of $O_\uparrow(n, 1)$ on $G^T_k$. The
    property $\mathrm{Id} \cdot g = g$ is immediate. Let $A, B \in O_\uparrow(n,
        1)$ and $g \in G^T_k$ be given. Note that $\pi(A_{B\cdot g} B_g) =
        \pi(A_{B\cdot g}) \pi(B_g) = AB$ and
    \[
        (A_{B\cdot g})_* (B_g)_* g = (A_{B\cdot g})_* (B\cdot g) \in G^T_k.
    \]
    From the discussion above, we conclude that $A_{B\cdot g} \circ B_g = (AB)_g$
    and hence
    \[
        A \cdot(B \cdot g)
        = (A_{B\cdot g})_* (B_g)_*
        = ((AB)_g)_* g = (AB) \cdot g.
    \]
\end{proof}

\subsection{Jets of asymptotically hyperbolic metrics}\label{secJets}

The linear masses we define should be functions of the germ at infinity of an
asymptotically hyperbolic metric. There is however an important caveat
preventing us from using such a definition. Namely, we want to impose some
continuity assumption for the invariants, and the problem is then that the
stalk of asymptotically hyperbolic metrics has no natural topology (we refer
the reader to the notes~\cite[Chapter 2]{TaylorTVS}, see also~\cite[Chapter 2,
    Section 6]{SchaeferWolff}). For this reason we introduce the set of $l$-jets of
asymptotically hyperbolic metrics.

\begin{definition}\label{defJets}
    For positive integers $k,l$, we define the set of \emph{$l$-jets of asymptotically
        hyperbolic metrics} of order $k$ as
    \[
        J^l_k \definedas G_k / \sim_{l+1}
    \]
    where the equivalence relation $\sim_{l+1}$ is defined by $g_1 \sim_{l+1} g_2$
    if and only if $\left|g_1 - g_2\right|_b = O(\rho^{l+1})$. We denote the
    projection from $G_k$ to $J^l_k$ by $\Pi^l_k$. A jet $j \in J^l_k$ is called
    \emph{transverse} if there exists a germ $g \in G_k$ representing $j$ which is
    transverse. We denote the set of transverse jets in $J^l_k$ by $T^l_k$.
\end{definition}
\noindent
Note that if $l < k$, the sets $J^l_k$ and $T^l_k$ each consist of a single element,
namely the jet of the hyperbolic metric $b$. For general $k$ and $l$, the sets
$J^l_k$ and $T^l_k$ are naturally affine spaces, with $b$ as a distinguished
base point. It is therefore natural to consider the associated vector spaces, obtained by identifying each jet with its difference from the reference hyperbolic jet:

\begin{definition}\label{defVectorizedJets}
    We define the \emph{vector space of $l$-jets of asymptotically hyperbolic metrics of order $k$} as
    \[
        \cJ^l_k \definedas J^l_k - [b],
    \]
    where $[b] \in J^l_k$ denotes the jet of the hyperbolic metric. Similarly, we
    define the \emph{vector space of transverse $l$-jets} as
    \[
        \cT^l_k \definedas T^l_k - [b].
    \]
    Finally, we let $\vec{\Pi}^l_k: G_k \to \cJ^l_k$ be defined as $\vec{\Pi}^l_k
        (g) = \Pi^l_k(g) - [b]$.
\end{definition}
Now remark that there exist (linear) isomorphisms
\[
    \cJ^l_k \simeq \bigoplus_{j=k}^l \Gamma(\Sym^2(T^* \bR^n \vert_{\bS^{n-1}})),
    \quad \cT^l_k \simeq \bigoplus_{j=k}^l \Gamma(\Sym^2(T^*\bS^{n-1})),
\]
where, for any given jet $[g] \in J^l_k$ (resp.\ $[g] \in T^l_k$),
\[
    g = \rho^{-2} (\delta + \rho^k g^{(k)} + \cdots + \rho^l g^{(l)} + O(\rho^{l+1})),
\]
$[g] - [b]$ is mapped to
\[
    (g^{(k)}, \ldots, g^{(l)}).
\]
These isomorphisms naturally endow topologies on $J^l_k$ and $T^l_k$ as
follows. The Levi-Civita connection of $\bR^n$ induces a connection on
$\operatorname{Sym}^2(T^* \bR^n \vert_{\bS^{n-1}})$ which allows us to define
the standard $C^{\infty}$ Fr\'echet space topology on the space of smooth
sections $\Gamma(\operatorname{Sym}^2(T^* \bR^n \vert_{\bS^{n-1}}))$. The
relevant theory of Fr\'echet spaces can be found in~\cite{HamiltonNashMoser,
    RudinFunctionalAnalysis}.

In this terminology, an asymptotically hyperbolic germ can be thought of as a
germ at $\rho = 0$ of curves $\rho \mapsto g(\rho)$ defined on an interval of
the form $(0, \epsilon)$ such that $\gbar(\rho) \definedas \rho^2 g(\rho)$
extends smoothly to the interval $[0, \epsilon)$. In polar coordinates the
hyperbolic metric reads
\[
    b = \rho^{-2}\left(\frac{d\rho^2}{1-2\rho} + (1-2\rho) \sigma\right),
\]
where $\sigma$ stands for the round metric on the unit sphere.

A metric $g$ is then transverse if $g - b$ is a $1$-parameter family of
sections of $\operatorname{Sym}^2(T^*\bS^{n-1})$ extended trivially in the
$\rho$-direction. The set of $l$-jets of asymptotically hyperbolic metrics is
thus identified with $l+1$ copies of $\Gamma(\operatorname{Sym}^2(T^* \bR^n
    \vert_{\bS^{n-1}}))$ via
\[
    g \mapsto
    \left(
    \gbar(0), \partial_\rho\gbar(0), \dots, \partial^{l}_\rho\gbar(0)
    \right).
\]
We define the topology on $\cJ^l_k$ (resp.\ on $\cT^l_k$) and, hence, on
$J^l_k$ (resp.\ on $T^l_k$) as the product topology through this
identification.

We now relate jets of asymptotically hyperbolic metrics to the theory developed
in Subsection~\ref{secAsymptoticallyHyp}. The next proposition shows that the
Taylor expansion of the diffeomorphism $\Theta$ can be obtained formally from
the Taylor expansion of the metric $\gbar = \rho^2 g$.

\begin{proposition}\label{propAdjustmentJet}
    Given $g \in G_k$, we let $\Theta$ be the adjustment diffeomorphism
    corresponding to $g$ (see Proposition~\ref{propTransverse}). For any
    $l > k$, the $(l+1)$-jet of $\Theta$ is fully determined by the $l$-jet
    of $g$ and depends smoothly on it.
    In particular, if $\Pi^l_k(g) \in T^l_k$, the $(l+1)$-jet of $\Theta$ coincides
    with the identity: $\Theta \in I_0^{l+2}$.
\end{proposition}

\begin{proof}
    Set $\Psi = \Theta^{-1}$. We will show that the $(l+1)$-jet of $\Psi$
    is determined by the $l$-jet of $g$. We assume that
    $\Psi \in I^{k+1}_0$ is such that $\Psi^* g$ is transverse, that is
    \[
        \rho^2 \Psi^* (\rho^{-2} \gbar)_{ij} x^i = x^j.
    \]
    Introducing ``polar'' coordinates $(\rho, \phi^A)$, where $(\phi^A)$ are
    coordinates on an open subset of the sphere $\bS^{n-1}$, the transversality
    condition reads
    \[
        \rho^2 \Psi^* (\rho^{-2} \gbar)(\partial_\rho, \cdot)
        = \rho^2 b(\partial_\rho, \cdot) = \frac{d\rho}{1-2\rho}.
    \]
    We use the index $0$ to denote the $\rho$ direction and capital letter indices
    ranging from $1$ to $n-1$ to denote directions tangent to the sphere. Using
    this convention, the transversality condition can be rephrased as
    \[
        \gbar_{ab}(\Psi(x)) \partial_0 \Psi^a(x) \partial_j \Psi^b(x)
        = \left(\frac{\rho(\Psi(x))}{\rho(x)}\right)^2 \frac{\delta^0_j}{1-2\rho}.
    \]
    From Proposition~\ref{propTransverse}, we have $|\Psi(x) - x| = O(\rho^{k+1})$,
    with $k \geq 1$. Hence,
    \[
        \Psi^0(x) =
        \rho(x) + \psi^0(x),\qquad \Psi^A(x) = \phi^A(x) + \psi^A(x),
    \]
    with $\psi^0(x), \psi^A(x) = O(\rho^{k+1})$. Introducing this into the
    transversality condition, we get
    \[
        \left(1+\frac{\psi^0}{\rho}\right)^2 \frac{\delta^0_j}{1-2\rho}
        = \bbar_{ab} \partial_0 \Psi^a(x) \partial_j \Psi^b(x)
        + (\gbar_{ab}(\Psi(x)) - \bbar_{ab}) \partial_0 \Psi^a(x) \partial_j \Psi^b(x)
    \]
    In the cases $j = 0$ and $j = C \neq 0$ this gives
    \begin{subequations}
        \begin{align}
            \left(1+\frac{\psi^0}{\rho}\right)^2 \frac{1}{1-2\rho}
             & = \frac{1}{1-2\rho}\left(1+\partial_0 \psi^0\right)^2
            + (1-2\rho) \sigma_{AB} \partial_0 \psi^A \partial_0 \psi^B\nonumber               \\
             & \qquad + (\gbar_{ab}(\Psi(x)) - \bbar_{ab})
            \partial_0 \Psi^a(x) \partial_0 \Psi^b(x),\label{eqRadialRadial}                   \\
            0
             & = \frac{1}{1-2\rho}\left(1+\partial_0 \psi^0\right) \partial_C \psi^0 \nonumber \\
             & \qquad + (1-2\rho) \sigma_{AB} \partial_0 \psi^A \left(\delta^B_C
            + \partial_C \psi^B\right)\nonumber                                                \\
             & \qquad + (\gbar_{ab}(\Psi(x)) - \bbar_{ab})
            \partial_0 \Psi^a(x) \partial_C \Psi^b(x)\label{eqRadialTangential}.
        \end{align}
    \end{subequations}
    The first equation can be rewritten as
    \begin{equation}\label{eqRadialRadial2}\tag{\ref*{eqRadialRadial}'}
        \begin{split}
             & \left(2+\frac{\psi^0}{\rho} + \partial_0 \psi^0\right)
            \left(\frac{\psi^0}{\rho} - \partial_0 \psi^0\right)        \\
             & \qquad =
            (1-2\rho)^2 \sigma_{AB} \partial_0 \psi^A \partial_0 \psi^B \\
             & \qquad\qquad
            + (1-2\rho)(\gbar_{ab}(\Psi(x)) - \bbar_{ab})
            \partial_0 \Psi^a(x) \partial_0 \Psi^b(x).
        \end{split}
    \end{equation}
    We next make Taylor expansions of $\psi^0$, $\psi^A$ in $\rho$,
    \begin{equation}\label{eqTaylorPsi}
        \left\lbrace
        \begin{aligned}
            \psi^0
             & = \psi^0_2 \rho^2 + \cdots + \psi^0_{l+1} \rho^{l+1} + O(\rho^{l+2}), \\
            \psi^A
             & = \psi^A_2 \rho^2 + \cdots + \psi^A_{l+1} \rho^{l+1} + O(\rho^{l+2}).
        \end{aligned}
        \right.
    \end{equation}
    Note that the first two terms in the Taylor expansion disappear since we assumed
    that $\psi^i = O(\rho^2)$. The proof now goes by induction on $l$. The idea is that,
    having determined $\psi^0$ and $\psi^A$ up to order $\rho^l$, we can determine the
    coefficient $\psi^0_{l+1}$ by expanding Equation~\eqref{eqRadialRadial2} up to order
    $\rho^l$. The only place where $\psi^0_{l+1}$ shows up when Taylor
    expanding~\eqref{eqRadialRadial2} up to order $O(\rho^{l+2})$ is in
    \[
        \frac{\psi^0}{\rho} - \partial_0 \psi^0
        = - \rho \psi^0_2 - 2 \rho^2 \psi^0_3 - \cdots - l \rho^l \psi^0_{l+1}
        + O(\rho^{l+1}).
    \]
    We can then determine $\psi^A_{l+1}$ by looking at the Taylor expansion of
    Equation~\eqref{eqRadialTangential} up to order $O(\rho^{l+2})$. Also here,
    $\psi^A_{l+1}$ shows up only at one place, namely in the term $(1-2\rho)
        \sigma_{AB} \partial_0 \psi^A \delta^B_A$.

    It is important at each step to notice that in order to determine the
    coefficients $\psi^0_{l+1}$ and $\psi^A_{l+1}$, we only need the Taylor
    expansion of the coefficients $\gbar_{ab}$ up to order $\rho^l$, and that the
    coefficients of the Taylor expansions~\eqref{eqTaylorPsi} are actually
    polynomials in $\gbar_{ab}$ and its derivatives at $\rho=0$. Since the details
    of the proof are rather involved, we leave them to the interested reader.
\end{proof}

Proposition~\ref{propAdjustmentJet} also implies that the
``transversalization'' operation via the adjustment diffeomorphism descends to
jets. This is the content of the next proposition.

\begin{proposition}\label{propLinearity}
    There exists a map
    \[
        \theta: \cJ^l_k \to \cT^l_k
    \]
    with the property that the diagram
    \begin{center}
        \begin{tikzcd}
            G_k \arrow{r}{\Theta_*} \arrow{d}{\vec{\Pi}_k^l} & G_k^T \arrow{d}{\vec{\Pi}_k^l}\\
            \cJ^l_k \arrow[dotted]{r}{\theta} & \cT^l_k
        \end{tikzcd}
    \end{center}
    commutes. Here $\Theta_*$ is (with some abuse of notation) the map associating to a given
    element $g \in G_k$ the element $\Theta_* g \in G^T_k$, where $\Theta \in I_0^{k+1}$ is
    the adjustment diffeomorphism defined in Proposition~\ref{propTransverse}. Assuming that
    $l < 2k$, the map $\theta$ is linear in $[g]-[b]$. More generally, for any $l'$ such that
    $0 \leq l' \leq l$, $\theta$ is linear in $[g] - [b]$ when restricted to jets $[g] \in J^l_k$
    which are transverse up to order $l'$ (i.e.\ which project to elements of $T_k^{l'}$) as long
    as $l \leq k + l'$.
\end{proposition}

\begin{proof}
    The proof of the first part of the proposition follows directly from
    Proposition~\ref{propAdjustmentJet}. Left to prove is that $\theta$ is linear
    in $g-b$ when $l \leq k + l'$ and $g$ projects to an element of $T_k^{l'}$ (note
    that the linearity of $\theta$ for $l < 2k$ follows from the general case by
    taking any $l' = k-1$).

    Let $g_0$ and $g_1$ be two elements of $G_k$ which project to elements of
    $T_k^{l'}$. For $\lambda \in [0,1]$, we set
    \[
        g_\lambda \definedas (1-\lambda) g_0 + \lambda g_1
    \]
    and denote by $\Psi_\lambda$ the unique element in $I^{k+1}_0$ for which
    $\Psi_\lambda^* g_\lambda$ is transverse. Since $g_\lambda$ projects to an
    element of $T_k^{l'}$, it follows from Proposition~\ref{propAdjustmentJet}
    applied with $l = l'$ that $\Psi_\lambda \in I_0^{l'+2}$. We wish to show that
    \[
        \Psi^*_\lambda g_\lambda = (1-\lambda) \Psi^*_0 g_0 + \lambda \Psi^*_1 g_1,
    \]
    at the level of $l$-jets. Since the right hand side is affine in $\lambda$ it
    is enough to show that
    \[
        \rho^2 \frac{d^2}{d\lambda^2} \Psi^*_\lambda g_\lambda = 0,
    \]
    once again at the level of $l$-jets. Note that we have not proven that
    $\Psi_\lambda$ depends in a $C^2$-manner on $\lambda$ but from the previous
    proposition, we know that the $l$-jet of $\Psi_\lambda$ depends smoothly on
    $\lambda$ and this is all that will be needed. We will compute this second
    order derivative at some $\lambda_0 \in [0, 1]$. To simplify calculations, we
    can replace the metrics $g_0$ and $g_1$ by $\Psi_{\lambda_0}^* g_0$ and
    $\Psi_{\lambda_0}^* g_1$ and hence assume that $\Psi_{\lambda_0}$ is the
    identity and that $g_{\lambda_0}$ is transverse.

    Using component notation we compute at $\lambda= \lambda_0$,
    \[\begin{split}
             & \rho^2 \frac{d^2}{d\lambda^2} \left(\Psi^*_\lambda g_\lambda\right)_{ij}  \\
             & \qquad=
            6 \frac{\partial_a \rho}{\rho} \frac{\partial_b \rho}{\rho}
            \Psi_\lambda'^a \Psi_\lambda'^b \gbar_{\lambda ij}(x)
            -2 \frac{\partial_a \partial_b \rho}{\rho}
            \Psi_\lambda'^a \Psi_\lambda'^b \gbar_{\lambda ij}(x)
            - 2 \frac{\partial_a \rho}{\rho} \Psi_\lambda''^a \gbar_{\lambda ij}(x)      \\
             & \qquad\qquad
            - 4 \frac{\partial_a \rho}{\rho} \Psi_\lambda'^a \gbar'_{\lambda ij}
            - 4 \frac{\partial_a \rho}{\rho} \Psi_\lambda'^a \partial_b \gbar_{\lambda ij}
            \Psi_\lambda'^b                                                              \\
             & \qquad\qquad
            - 4 \frac{\partial_a \rho}{\rho} \Psi_\lambda'^a \gbar_{\lambda kj}(x)
            \partial_i\Psi_\lambda'^k
            - 4 \frac{\partial_a \rho}{\rho} \Psi_\lambda'^a \gbar_{\lambda il}(x)
            \partial_j\Psi_\lambda'^l                                                    \\
             & \qquad\qquad
            + 2\partial_a \gbar_{\lambda ij}' \Psi'^a
            + 2 \gbar'_{\lambda kj} \partial_i \Psi'^k_\lambda
            + 2 \gbar'_{\lambda il} \partial_j \Psi'^l_\lambda                           \\
             & \qquad\qquad
            + \partial_a \partial_b \gbar_{\lambda ij} \Psi_\lambda'^a \Psi_\lambda'^b
            + \partial_a \gbar_{\lambda ij} \Psi_\lambda''^a                             \\
             & \qquad\qquad
            + 2 \partial_a \gbar_{\lambda kj} \Psi'^a_\lambda \partial_i \Psi'^k_\lambda
            + 2 \partial_a \gbar_{\lambda il} \Psi'^a_\lambda \partial_j \Psi'^l_\lambda \\
             & \qquad\qquad
            + \gbar_{\lambda kj}(x) \partial_i \Psi''^k_\lambda
            + \gbar_{\lambda il}(x) \partial_j \Psi''^l_\lambda
            + 2 \gbar_{\lambda kl}(x) \partial_i \Psi'^k_\lambda \partial_j \Psi'^l_\lambda .
        \end{split}\]
    Since $\Psi', \Psi'' = O(\rho^{l'+2})$, it follows by inspection of the decay
    order of each term that
    \begin{equation}\label{eqSecondDerivative}
        \begin{split}
            \rho^2 \frac{d^2}{d\lambda^2} \left(\Psi^*_\lambda g_\lambda\right)_{ij}
             & = - 2 \frac{\partial_a \rho}{\rho} \Psi_\lambda''^a \gbar_{\lambda ij}(x) + \gbar_{\lambda kj}(x) \partial_i \Psi''^k_\lambda \\
             & \qquad \qquad + \gbar_{\lambda il}(x) \partial_j \Psi''^l_\lambda + O(\rho^{k+l'+1}).
        \end{split}
    \end{equation}
    To prove that the remaining three terms are also $O(\rho^{k+l'+1})$, we need to take a closer look at Equations~\eqref{eqRadialRadial2} and~\eqref{eqRadialTangential}. Let $\psi^0_\lambda$ and $\psi^A_\lambda$ be such that
    \[
        \Psi_\lambda^0(x) = \rho(x) + \psi^0_\lambda(x), \qquad \Psi_\lambda^A(x) = \phi^A(x) + \psi^A_\lambda(x).
    \]
    Since $\Psi_\lambda \in I_0^{l'+2}$, the Taylor expansions of $\psi^0_\lambda$
    and of $\psi^A_\lambda$ now start at the order $O(\rho^{l'+2})$:
    \[
        \left\lbrace
        \begin{aligned}
            \psi^0_\lambda
             & = \psi^0_{\lambda, l'+2} \rho^{l'+2} + \cdots + \psi^0_{\lambda, l+1} \rho^{l+1} + O(\rho^{l+2}), \\
            \psi^A_\lambda
             & = \psi^A_{\lambda, l'+2} \rho^{l'+2} + \cdots + \psi^A_{\lambda, l+1} \rho^{l+1} + O(\rho^{l+2}).
        \end{aligned}
        \right.
    \]
    From the fact that $\psi^0_\lambda$ and $\psi^A_\lambda$ are $O(\rho^{l'+2})$,
    Equations~\eqref{eqRadialRadial2} and~\eqref{eqRadialTangential} reduce to
    \begin{align*}
        2 \left(\frac{\psi^0_\lambda}{\rho} - \partial_0 \psi^0_\lambda\right)
         & =
        (1-2\rho)(\gbar_{00}(\Psi_\lambda(x)) - \bbar_{00}) + O(\rho^{k+l'+1})         \\
         & =
        (1-2\rho)(\gbar_{00}(x) - \bbar_{00} + \partial_a \bbar_{00} \psi^a_\lambda(x) \\
         & \qquad
        + \partial_a (\gbar_{00}-\bbar_{00}) \psi^a_\lambda(x) ) + O(\rho^{k+l'+1}),   \\
         & =
        (1-2\rho)(\gbar_{00}(x) - \bbar_{00} + \partial_a \bbar_{00} \psi^a_\lambda(x))
        + O(\rho^{k+l'+1}),                                                            \\
        0
         & =
        \frac{1}{1-2\rho} \partial_C \psi^0_\lambda
        + (1-2\rho) \sigma_{AC} \partial_0 \psi^A_\lambda                              \\
         & \qquad + (\gbar_{0C}(\Psi_\lambda(x)) - \bbar_{0C}) + O(\rho^{k+l'+1})      \\
         & =
        \frac{1}{1-2\rho} \partial_C \psi^0_\lambda
        + (1-2\rho) \sigma_{AC} \partial_0 \psi^A_\lambda                              \\
         & \qquad + \gbar_{0C}(x) + O(\rho^{k+l'+1}),
    \end{align*}
    where we used the fact that $\bbar_{0C} = 0$ in the second calculation. It follows by
    induction that $\psi$ depends linearly on $\gbar - \bbar$ up to terms of order
    $O(\rho^{k+l'+1})$. This implies that $\Psi_\lambda'' = O(\rho^{k+l'+1})$. So, from
    Equation~\eqref{eqSecondDerivative}, we deduce that
    \[
        \vec{\Pi}_k^l \left(\frac{d^2}{d\lambda^2} \Psi^*_\lambda g_\lambda\right) = 0
    \]
    provided $l < k + l' + 1$.
\end{proof}

In the next proposition we compute the first non-trivial term in the asymptotic
expansion of $\Theta_* g$.

\begin{proposition}\label{propAdjMass}
    Let $g \in G_k$. We denote by $m$ the first non-trivial term in the asymptotic expansion of $g$, that is
    \[
        g = b + \rho^{k-2} m + O(\rho^{k-1})
    \]
    where $m$ is a section of $\operatorname{Sym}^2(T^*\bR^n\vert_{\bS^{n-1}})$.
    Let $\Theta$ be the adjustment diffeomorphism of $g$ from
    Proposition~\ref{propTransverse}. The metric $\gtil \definedas \Theta_* g$ has
    the asymptotic expansion
    \[
        \gtil = b + \rho^{k-2} \mtil + O(\rho^{k-1}),
    \]
    where the section $\mtil$ of $\operatorname{Sym}^2(T^*\bS^{n-1})$ is given by
    \begin{equation}\label{eqAdjustedMassAspect}
        \mtil_{ij}
        = m_{ij} - m_{aj} x^a x_i - m_{ia} x^a x_j
        + m_{ab} \frac{x^a x^b}{k} \left( (k-1)x_i x_j + \delta_{ij}\right).
    \end{equation}
\end{proposition}

\begin{proof}
    Calculations starting from Equations~\eqref{eqRadialRadial} and~\eqref{eqRadialTangential} are fairly straightforward. Indeed, the first non-trivial terms in the asymptotic expansion of $\Psi = \Theta^{-1}$ are the ones of order $\rho^{k+1}$. They can be obtained from the terms of order $O(\rho^k)$ in~\eqref{eqRadialRadial}-\eqref{eqRadialTangential}, that is
    \[
        \psi^0_{k+1} = - \frac{1}{2k} m_{00}, \quad
        \sigma_{AB}\psi^B = -\frac{m_{0B}}{k+1}.
    \]
    Hence, from the identity
    \[
        \rho^2 \gtil = \left(\frac{\rho}{\rho\circ\Psi}\right)^2 \Psi^* \gbar
    \]
    we find that
    \[
        \mtil_{00} = 0, \quad
        \mtil_{0A} = 0, \quad
        \mtil_{AB} = m_{AB} + \frac{m_{00}}{k} \sigma_{AB}.
    \]
    The relation between $m$ and $\mtil$ can be condensed into
    \[
        \mtil =
        m - d\rho \otimes m(\cdot, \partial_\rho)
        - m(\partial_\rho, \cdot) \otimes d\rho
        + \frac{m(\partial_\rho, \partial_\rho)}{k}
        \left(\delta + (k-1) d\rho \otimes d\rho\right).
    \]
    We now wish to return to Cartesian coordinates. To do this, it suffices to note
    that they are related to the $(\rho, \phi^A)$-coordinates via
    \[
        x^i = \sqrt{1-2\rho} F^i(\phi),
    \]
    where $(F^i)$ is a set of $n$ given functions such that $\sum (F^i)^2 = 1$,
    meaning that $F^i = \frac{x^i}{|x|}$. As a consequence,
    \[
        \frac{\partial}{\partial \rho}
        = \frac{\partial x^i}{\partial \rho} \frac{\partial}{\partial x^i}
        = - \frac{1}{\sqrt{1-2 \rho}} F^i \frac{\partial}{\partial x^i}
        = - \frac{1}{|x|^2} x^i \frac{\partial}{\partial x^i}.
    \]
    Formula~\eqref{eqAdjustedMassAspect} follows.
\end{proof}

The following theorem will allow us to define asymptotic invariants.

\begin{theorem}\label{thmActionJets}
    For any $l \geq k$, there exists a unique action of the group $O_\uparrow(n, 1)$ on $\cT^l_k$ such that the projection $\vec{\Pi}_k^l: G^T_k \to \cT^l_k$ is a $O_\uparrow(n,1)$-equivariant map. Namely,
    \[
        A \cdot \left(\vec{\Pi}_k^l(g)\right) = \vec{\Pi}_k^l(A \cdot g),
    \]
    for all $A \in O_\uparrow(n, 1)$ and all $g \in G_k^T$, where the action of
    $O_\uparrow(n, 1)$ on $G^T_k$ was defined in
    Proposition~\ref{propActionLorentz}. Further, the action is linear and smooth
    as long as $l \leq 2k$ and reduces in the case $l=k$ to the pushforward action
    \[
        A \cdot g = \Abar_* g,
    \]
    where $A \in O_\uparrow(n, 1)$ is any hyperbolic isometry.
\end{theorem}

\begin{proof}
    The only non-trivial point in the proof is the fact that the action is linear for all $l \leq 2k$. We shall prove that we have $\Abar_* g \in T^k_k$ for $g \in T^l_k$ and $A\in O_\uparrow(n, 1)$, which means that transversality of the first non-trivial term in the asymptotic expansion of $g$ is preserved under the (non-adjusted) action of $O_\uparrow(n, 1)$. From Proposition~\ref{propLinearity}, it then follows that $g \mapsto \theta \Abar_* g$ is linear in $g-b$ for $l \leq 2k$.

    As before we denote by $P = x^i \partial_i$ the position vector field and we
    set $g = b + e$. Since $A$ is a hyperbolic isometry, we have
    \[
        \Abar_* g = b + \Abar_* e.
    \]
    Transversality of $g$ reads $e(P, \cdot) = 0$, so we need to check that
    \[
        \left|\rho^2(\Abar_* e)(P, \cdot)\right|_{\delta}
        = \left|\rho^2 \Abar_*\left(e(\Abar^* P, \cdot)\right)\right|_{\delta}
        = O(\rho^{k+1}).
    \]
    This will follow by showing that $\Abar^* P = \lambda(x, A) P + O(\rho)$ for
    some function $\lambda(x, A)$. Before starting calculations, we note that this
    fact is natural from the point of view of hyperbolic geometry. Indeed, the
    action of a hyperbolic isometry can be understood as a change of origin of the
    hyperbolic space and $P$ is a vector field pointing in the direction of
    geodesics emanating from the origin. All hyperbolic geodesics intersect the
    boundary $\bS^{n-1}$ orthogonally. Hence, a tensor that is transverse with
    respect to some choice of an origin remains transverse up to correction terms
    with respect to the new origin.

    The condition
    \begin{equation}\label{eqTranversalityP}
        \Abar^* P = \lambda(x, A) P + O(\rho)
    \end{equation}
    only needs to be checked for generators of the Lorentz group. It is
    obvious for rotations since $\Rbar^* P = P$ for any rotation $R$. The case
    of Lorentz boosts requires some calculations.

    With $A$ replaced by $A^{-1}$, Condition~\eqref{eqTranversalityP} is equivalent
    to
    \[
        \Abar_* P = \lambda(x, A^{-1}) P + O(\rho),
    \]
    or to
    \[
        d\Abar(P)(x) = \lambda(\Abar(x), A^{-1}) P(\Abar(x)) + O(\rho).
    \]
    We use Formula~\eqref{eqBoosts} to write
    \[\begin{split}
            d\Abar^s_i(x)
             & =
            \frac{1}{D}
            \left(dx^1, \dots, \cosh(s) dx^i + \sinh(s) x^j dx^j, \dots, dx^n\right)     \\
             & \qquad
            - \frac{1}{D^2}
            \left(x^1, \dots, \cosh(s) x^i + \sinh(s) \frac{1+|x|^2}{2},\dots,x^n\right) \\
             & \qquad\qquad\qquad
            \cdot
            \left((\cosh(s)-1) x^j dx^j + \sinh(s) dx^i\right),                          \\
        \end{split}\]
    so
    \[\begin{split}
            d\Abar^s(P)
             & =
            \left(\frac{1}{D} - \frac{(\cosh(s)-1) |x|^2 + \sinh(s) x^i}{D^2}\right) \\
             & \qquad\qquad
            \cdot
            \left(x^1, \dots, \cosh(s) + \sinh(s) \frac{1+|x|^2}{2},\dots,x^n\right)
            + O(\rho)                                                                \\
             & =
            \frac{1}{\cosh(s) + \sinh(s) x^i} P(\Abar^s_i(x)) + O(\rho),
        \end{split}\]
    and the claim follows.
\end{proof}

\subsection{Linear masses at infinity}\label{secGeomMasses}

We now come to the central definition of this article.

\begin{definition}\label{defGeometricMass}
    Let $V$ be a finite-dimensional representation of the group
    $O_\uparrow(n, 1)$.  A \emph{linear mass at infinity}
    for the set of asymptotically hyperbolic metrics of order $k$ is a map
    $\Phi: G_k \to V$ satisfying the following properties.
    \begin{enumerate}
        \item
              $\Phi(b) = 0$.
        \item
              For any $g \in G_k$ and any $\Psi \in I^k(g)$ it holds that
              \[
                  \Phi(\Psi_* g) = \pi(\Psi) \cdot \Phi(g).
              \]
              That is, $\Phi$ is an intertwining map.
        \item The map $\Phi$ factors through $\cT^k_k$,
              \begin{center}
                  \begin{tikzcd}
                      G_k \arrow{rr}{\Phi}\arrow{rd}[swap]{\widetilde{\Pi}^k_k} & & V\\
                      & \cT^k_k \arrow[dotted]{ru}[swap]{\phi}
                  \end{tikzcd}
              \end{center}
              where $\widetilde{\Pi}^k_k$ denotes the composition
              $\vec{\Pi}^k_k \circ \Theta_*$, and $\Theta_*$ is (with some abuse of notation) the map which takes an element $g \in G_k$ to the unique transverse element $\Theta_* g \in G^T_k$ where $\Theta \in I^{k+1}_0$.
        \item
              The map $\phi$ is continuous and linear.
    \end{enumerate}
\end{definition}

Let us make a few comments on this definition.

\begin{remarks}
    \begin{itemize}
        \item
              At this point we could be more general assuming only that the
              space $V$ is a $O_\uparrow(n, 1)$-set. The restriction to vector space
              representations is in a sense irrelevant, since we can embed
              any reasonable $O_\uparrow(n, 1)$-set into a (linear) representation of
              $O_\uparrow(n, 1)$, see~\cite[Chapter~7,~Section~1.3]{Procesi}. The only
              non-trivial point is that the map $\Phi: G_k \to V$ has to be linear.
              Invariants depending polynomially on the leading term in the asymptotic
              expansion of $g - b$ (the so-called mass aspect) have been
              studied for example in~\cite{LiNguyen}.
        \item
              We call these objects ``masses'' instead of ``invariants'' since we
              strictly speaking do not get something independent of the chosen chart
              at infinity. The situation is the same for the classical asymptotically
              hyperbolic mass. If one insists on having a true invariant, one has to
              look at the $O_\uparrow(n,1)$-orbit to which $\Phi(g)$ belongs. This
              relies on the classical invariant theory described in
              ~\cite[Chapter~11]{Procesi} or~\cite[Chapter~5]{GoodmanWallach}
              and will be addressed in future work.
        \item
              It is possible to extend the definition to linear masses factoring through jets
              $T^l_k$ for $l>k$. This will also be studied in forthcoming work. Preparing for
              this is the reason for studying general jets in the present paper. \item The
              condition $\Phi(b) = 0$ is imposed since we want our invariants to measure the
              difference between some given metric and the hyperbolic metric. This condition
              will be immediately fulfilled if we assume that $V$ has no $1$-dimensional
              trivial subrepresentation, in particular if $V$ is an irreducible
              representation with $\dim V > 1$. Indeed, the hyperbolic metric is a fixed
              point for the action of $O_\uparrow(n, 1)$ and under the condition stated
              above, the only such fixed point in $V$ is $0$.

    \end{itemize}
\end{remarks}

In the next section we will start the classification of the linear masses at
infinity.

%% file: lorentzaction.tex
%
%

\newcommand{\Fhi}{F}
\newcommand{\fhi}{{\bf F}}

We are going to classify all linear masses at infinity as defined in
Definition~\ref{defGeometricMass}. Our way of finding such maps $\Phi: G_k \to
    V$ is to begin by looking for continuous linear maps
\[
    \phi : \cT_k^k \to V,
\]
which are intertwining for the action of the group $O_{\uparrow}(n,1)$, where
$V$ is a finite-dimensional irreducible representation of the group
$O_{\uparrow}(n,1)$.

Note that whenever such a map $\phi$ is found, one immediately gets from
Theorem~\ref{thmActionJets} that the map $\Phi : G_k \to V$ defined as
\[
    \Phi \definedas \phi \circ \widetilde{\Pi}^k_k
\]
satisfies the requirements of Definition~\ref{defGeometricMass}.

On the way to describe the action of the orthochronous Lorentz group
$O_{\uparrow}(n,1)$ on such ($k$-jets of) metrics, we notice that it is in fact
equivalent to describe the action on the \emph{mass aspect tensors}, that is
the first non-trivial term $\tilde m$ in the asymptotic expansion of the metric
$\tilde g$ in Proposition~\ref{propAdjMass}, see
Lemma~\ref{lmGroupActionMassAspect}.

One motivation for considering directly the group action (and intertwining
maps) on the set of mass aspect tensors is that it allows one to read the
linear masses at infinity as expressions computed directly at the conformal
boundary of the asymptotically hyperbolic manifold (here the sphere with its
standard conformal class), and not through a limit process. This is precisely
the idea of the definition of the mass vector by Wang in~\cite{WangMass}.

Let us denote by $S^2(\bS^{n-1}) := \Gamma
    \big(\operatorname{Sym}^2(T^*\bS^{n-1})\big)$ the set of symmetric
$(0,2)$-tensors on $\bS^{n-1}$, also called \emph{mass aspect tensors}. After
equipping it with an $O_{\uparrow}(n,1)$-action, it follows from
Proposition~\ref{propMassaspect} that our quest for linear masses at infinity
reduces to looking for $O_{\uparrow}(n,1)$-intertwining maps
\[
    \Fhi : S^2(\bS^{n-1}) \to V
\]
sending a mass aspect tensor to an element of a finite-dimensional
representation $V$ of $O_{\uparrow}(n,1)$.

For the sake of simplicity, we shall first find such maps $\Fhi$ which are
intertwining with respect to the Lorentz Lie algebra $\mathfrak{so}(n,1)$. Note
that any finite-dimensional $O_{\uparrow}(n,1)$-representation $V$ is naturally
a $\mathfrak{so}(n,1)$-representation as well.

As $O_{\uparrow}(n,1)$ is not connected, the set of
$\mathfrak{so}(n,1)$-intertwining maps from $S^2(\bS^{n-1})$ to $V$ is larger
than the set of such maps which are intertwining for the action of the group
$O_{\uparrow}(n,1)$. We will however argue that all the maps we find in
Section~\ref{secClassification} are genuine $O_{\uparrow}(n,1)$-intertwining,
hence give linear masses at infinity.

\subsection{Action of \texorpdfstring{$O_{\uparrow}(n,1)$}{TEXT} on
    mass aspect tensors}
Let $g$ be a $k$-jet of metrics in $T_k^k$. As before, we will abuse notation
and also write $g$ for a representative metric in this class. Such a metric $g$
is then identified with its principal part of the asymptotic expansion,
\[
    g = b + \rho^{k-2} m \mod \ker(\widetilde{\Pi}_k^k),
\]
where $m \in S^2(\bS^{n-1})$. We say that $m$ is the \emph{mass aspect tensor
    of $g$}.

To emphasize the role of the decay parameter $k$ on the action of
$O_{\uparrow}(n,1)$ on the set $S^2(\bS^{n-1})$ of mass aspect tensors as
described in Lemma~\ref{lmGroupActionMassAspect} and
Corollary~\ref{remGroupAction}, we will use the formalism of \emph{conformal
    weights}, following~\cite{DjadliGuillarmouHerzlich} and~\cite{CurryGover}. In
the subsequent sections however, since the conformal class $[\sigma]$ of the
unit round metric $\sigma$ on $\bS^{n-1}$ is fixed, we will only consider mass
aspect tensors as ordinary tensors on $\bS^{n-1}$, for the sake of clarity.

In order to define the notion of conformal weights, we first introduce the
trivial line bundle $\cG = \bR_+^* \times \bS^{n-1}$ over $\bS^{n-1}$. Any
element $(t,x)$ of $\cG$ may then be identified with the quadratic form $t^2
    \sigma_x$. It follows that any metric $h$ in the conformal class $[\sigma]$
identifies with a section of $\cG$, since in this case, there is a smooth
positive function $\Omega$ on $\bS^{n-1}$ such that $h = \Omega^2 \sigma$.

\begin{definition}\label{defConfWeight}
    Given any real parameter $w$, a smooth map $f : \cG \to \cC^{\infty}(\bS^{n-1})$ is a \emph{conformal density of weight $w$} if, for all $x \in \bS^{n-1}$, the equality
    \[
        f(\varphi^*\sigma)_x = \Omega(x)^w f(\sigma)_{\varphi^{-1}(x)}
    \]
    holds for all conformal diffeomorphisms $\varphi$ such that $\varphi^*\sigma =
        \Omega^2 \sigma$, where $\Omega$ is a smooth positive function on $\bS^{n-1}$.
    We denote $\cE[w]$ the space of such conformal weight densities. Similarly, for
    a vector bundle $\cB$ over $\bS^{n-1}$, we denote by $\cB[w]$ the \emph{bundle
        of elements of $\cB$ with conformal weight $w$}, namely the bundle of smooth
    maps $f : \cG \to \cB$ satisfying
    \[
        f(\varphi^*\sigma)_x = \Omega(x)^w \varphi_* \left(f(\sigma)\right)_x\ ,
    \]
    for all $x \in \bS^{n-1}$, where the conformal diffeomorphism $\varphi$ and the
    smooth positive function $\Omega$ are again related by $\varphi^* \sigma =
        \Omega^2 \sigma$.
\end{definition}

Now, following the identification introduced after
Definition~\ref{defVectorizedJets}, we define an action of $O_{\uparrow}(n,1)$
on $S^2(\bS^{n-1})$. The following lemma describes the behaviour of the mass
aspect tensor of an asymptotically hyperbolic metric under the action of the
group of isometries of the hyperbolic space.

\begin{lemma}\label{lmGroupActionMassAspect}
    The action of $O_{\uparrow}(n,1)$ on $\cT_k^k$ defines a unique action
    \[
        (A,m) \mapsto A \cdot m
    \]
    on $S^2(\bS^{n-1})$. Moreover, if $m$ is the mass aspect tensor of a metric $g$
    in $\cT_k^k$, then for every $A \in O_{\uparrow}(n,1)$, there is a smooth,
    positive function $u[A]$ defined on $\bS^{n-1}$ such that
    \[
        A \cdot m = u[A]^{k-2} \Abar_* m .
    \]
\end{lemma}

\begin{proof}
    We define $A\cdot m$ as the image of $A \cdot g$ in $S^2(\bS^{n-1})$ where $g
        \in G_k$ is any germ of metric whose image in $\cT_k^k$ is $m$. One easily sees
    that this defines an $O_{\uparrow}(n,1)$-action on $S^2(\bS^{n-1})$. To compute
    this action, we use Proposition~\ref{propActionLorentz} which tells us that $A
        \cdot g = \widetilde{A}_* g$ with $\widetilde{A} = \Theta \circ \Abar$, where
    $\Theta$ is the unique adjustment diffeomorphism in $I_0^{k+1}$ of
    Proposition~\ref{propTransverse} associated to $\Abar_* g$. Both metrics
    $\widetilde{A}_* g$ and $\Abar_* g$ then lie in the same equivalence class
    modulo $\cT_k^k$. Next, we write
    \[
        \Abar_* g
        = b + \left(\rho \circ \Abar^{-1}\right)^{k-2} \Abar_* m
        + O\left( (\rho \circ \Abar^{-1})^{k-1} \right).
    \]
    For $x \neq 0$ we set $\hat x \definedas \frac{x}{|x|} \in \bS^{n-1}$, and
    define
    \begin{equation} \label{def_u[A]}
        u[A](\hat x) \definedas \lim_{\lambda \to 1}
        \frac{\rho (\Abar^{-1}(\lambda \hat x))}{\rho(\lambda \hat x)}.
    \end{equation}
    This defines $u[A]$ as a positive smooth function on $\bS^{n-1}$. We obtain
    \[
        \Abar_* g
        = b + \rho^{k-2} \big(u[A](\hat x)\big)^{k-2} \Abar_* m
        + O\left( \rho^{k-1} \right).
    \]
    Finally, we get the equality in $\cT_k^k$ (hence we omit the remainder),
    \[
        A \cdot g = b + \rho^{k-2} \big(u[A](\hat x)\big)^{k-2} \Abar_* m .
    \]
    Thus, we identify the expression for $A \cdot m$ as
    \[
        A \cdot m = \big(u[A](\hat x)\big)^{k-2} \Abar_* m.
    \]
\end{proof}

\begin{remark} Since $A$ induces a conformal diffeomorphism of the $\bS^{n-1}$ sphere, with the equality $\Abar_*\sigma = \left(\Abar^{-1}\right)^*\sigma = u[A]^2 \sigma$, we can rephrase this result by stating that the map $\cG \to S^2(\bS^{n-1})$, which sends any metric $\Abar_* \sigma$ to the mass aspect tensor $A \cdot m$, is a conformal density of weight $2-k$, along with Definition~\ref{defConfWeight}. We may then say that Lemma~\ref{lmGroupActionMassAspect} gives the action of $O_{\uparrow}(n,1)$ on $S^2(\bS^{n-1})[2-k]$ induced from the action on $\cT_k^k$.
\end{remark}


\begin{proposition}\label{propMassaspect}
    All linear masses at infinity are obtained as intertwining maps
    \[
        \Fhi : S^2(\bS^{n-1}) \to V
    \]
    between the set of mass aspect tensors of order $k$ and a finite-dimensional
    representation $V$ of $O_{\uparrow}(n,1)$.
\end{proposition}

Other consequences of Lemma~\ref{lmGroupActionMassAspect} are expressions for
the corresponding action on the \emph{mass aspect function}, $\tr ^{\sigma} m$,
as well as on the product $\tr^{\sigma} m \, d\mu ^{\sigma}$ of the mass aspect
function and the volume form. Here again, $\sigma$ stands for the round metric
on the unit sphere $\bS^{n-1}$.

\begin{corollary}\label{remGroupAction}
    We have
    \[
        A \cdot (\tr ^\sigma m) = (u[A])^{k} \tr ^\sigma m \circ \Abar^{-1}
    \]
    and
    \[
        A \cdot (\tr^{\sigma} m \, d\mu^{\sigma})
        = (u[A])^{k+1-n} \Abar_* \left(\tr^{\sigma} m \, d\mu^{\sigma}\right).
    \]
\end{corollary}

\begin{proof}
    The fact that $\Abar_* b = b$ translates into
    \[
        (\rho \circ \Abar^{-1})^{-2} \Abar_* \delta = \rho ^{-2} \delta,
    \]
    which we can rewrite as
    \[
        \Abar_* \sigma
        = \Abar_* \delta |_{|x| = 1}
        = \lim_{|x| \rightarrow 1}
        \left(\frac{\rho \circ \Abar^{-1}(x)}{\rho(x)}\right)^2 \sigma.
    \]
    This tells us that
    \[
        \Abar_* \sigma = u[A]^2 \sigma.
    \]
    To get the expression for the action on $\tr^{\sigma} m$ we write
    \[
        A \cdot g
        =
        \rho^{-2}
        \left[\sigma + \rho^k u[A]^{k-2}(\hat x) \Abar_* m
            + O(\rho^{k+1}) \right],
    \]
    so the mass aspect function of the metric $A \cdot g$ is
    \[
        A \cdot (\tr ^{\sigma} m)
        = u[A]^{k-2}(\hat x) \tr ^{\sigma} (\Abar_* m) .
    \]
    Using the identity
    \[
        \Abar_* (\tr^{\sigma} m)
        = \tr^{\Abar_* \sigma} (\Abar_* m)
        = u[A]^{-2} \tr ^{\sigma} (\Abar_* m),
    \]
    we find that $\tr^{\sigma}(\Abar_* m) = u[A]^2 \tr^{\sigma} m \circ \Abar^{-1}$
    and the result follows. For the action on the product $\tr^{\sigma}m \,
        d\mu^{\sigma}$, we write
    \[
        A_*(\tr^{\sigma} m \, d\mu^{\sigma})
        = (\tr^{\sigma} m \circ A^{-1}) \, d\mu^{A_*\sigma}
        = (\tr^{\sigma} m \circ A^{-1}) \, u[A]^{n-1} d\mu^{\sigma} ,
    \]
    and, again using the expression above for $A \cdot g$, we conclude that
    \[
        A \cdot (\tr^{\sigma} m \, d\mu^{\sigma})
        = \tr^{\sigma} (\Abar_* m) u[A]^{k-2} d\mu^{\sigma}
    \]
    which can be rewritten as
    \[
        u[A]^k \tr^{\sigma} m \circ \Abar^{-1} \, d\mu^{\sigma}
        = u[A]^{k+1-n} \Abar_*(\tr^{\sigma} m \, d\mu^{\sigma}) ,
    \]
    as desired.
\end{proof}

In particular, we have that $\tr^{\sigma} m \, d\mu^{\sigma}$ is an invariant
under the action of the group $\mathrm{Conf}(\bS^{n-1})$ of conformal
diffeomorphisms of the sphere for $k = n-1$.

\begin{remark} Following Definition~\ref{defConfWeight}, this result can be stated again in terms of \emph{conformal weights}. Precisely, the map sending the conformal metric $\Abar_*\sigma = u[A]^{-2} \sigma$ to the mass aspect function $A\cdot (\tr^{\sigma} m)$ is a conformal density of weight $-k$, meaning that we get the induced action of $O_{\uparrow}(n,1)$ on $\cC^{\infty}(\bS^{n-1})[-k]$. Similarly, the map sending the same conformal metric $\Abar_*\sigma$ to the $n$-form $A\cdot (\tr^{\sigma}m d\mu^{\sigma})$ has weight $n - 1 - k$.
\end{remark}



\subsection{Action of the Lorentz algebra on mass aspect tensors}

We now define the associated Lie algebra action of $\mathfrak{so}(n,1)$ on mass
aspect tensors in $S^2(\bS^{n-1})$. It is given by
\[
    X \cdot m \definedas
    \frac{\partial}{\partial s} \left(\Abar^s \cdot m\right)|_{s=0},
\]
for $X \in \mathfrak{so}(n,1)$ and $m \in S^2(\bS^{n-1})$, where $(A^s)_{s \in
            \bR}$ is the one-parameter subgroup of $O_{\uparrow}(n,1)$ generated by $X$. In
the next proposition we compute this action for Lorentz boosts $a_i$ and for
rotations $r_{ij}$, whose expressions were given in
Subsection~\ref{secHyperbolic}.

\begin{proposition}\label{propBoostRotationActonMassAspect}
    Let $\fa_i$ be the Lorentz boost vector field as defined in~\eqref{Xi}
    and let $\fr_{ij}$ be the rotation vector field defined
    in~\eqref{Lij}. Then
    \begin{align}
        a_i \cdot m
         & = -\nabla_{\fa_i}^{\sigma} m  + k x^i m, \label{X_icdotm} \\
        r_{ij} \cdot m
         & = - \nabla^{\sigma}_{\fr_{ij}} m
        - \left( m(r_{ij}(\cdot),\cdot) +  m(\cdot,r_{ij}(\cdot)) \right),
        \label{L_ijcdotm}
    \end{align}
    where $r_{ij}$ acts on vector fields tangent to $\bS^{n-1}$ by
    \[
        r_{ij} (U)
        = U^i \partial_j - U^j \partial_i
        = U^i \fa_j - U^j \fa_i \mod (x^a \partial_a).
    \]
\end{proposition}

\begin{proof}
    From the expression for the boost hyperbolic isometries $A^s_i$ defined
    in Section~\ref{secHyperbolic} we compute
    \[
        \rho \circ \Abar_i^{-s}
        =
        \frac{1}{\cosh s \frac{1+|x|^2}{1-|x|^2} - \sinh s \frac{2 x^i}{1-|x|^2} + 1}
        = \frac{\rho}{\cosh s - x^i \sinh s} + O(\rho^2)
    \]
    which gives us
    \[
        u[A_i^s](\hat x) = \frac{1}{\cosh s - x^i \sinh s}.
    \]
    Hence by Lemma~\ref{lmGroupActionMassAspect} we have the expression
    \[
        A^s _i \cdot m = \frac{1}{(\cosh s - x^i \sinh s)^{k-2}}
        \left(\Abar^s_i\right)_* m
    \]
    for the action of $A^s_i$. We compute the derivative of this expression with
    respect to $s$, and get
    \begin{equation}\label{LieXi}
        \begin{split}
            a_i \cdot m
             & = \frac{\partial}{\partial s} \left(
            \frac{1}{(\cosh s - x^i \sinh s)^{k-2}}
            \left(\Abar^s_i \right)_* m \right)|_{s=0} \\
             & = -\lie_{\fa_i} m + (k-2)x^i m .
        \end{split}
    \end{equation}
    Next, we want to rewrite this using covariant derivatives instead of
    Lie derivatives. Since elements of $O_\uparrow(n, 1)$ preserve the sphere
    at infinity $\bS^{n-1}$, the vector field $\fa_i$ is tangent to this
    sphere. We have
    \[ \begin{split}
            (\lie_{\fa_i} m)(\fa_a, \fa_b)
             & =
            \fa_i (m(\fa_a, \fa_b)) - m([\fa_i, \fa_a], \fa_b)
            - m(\fa_a,[\fa_i, \fa_b]) \\
             & =
            (\nabla^{\sigma}_{\fa_i} m)(\fa_a, \fa_b)
            + m(\nabla^{\sigma}_{\fa_a} \fa_i, \fa_b)
            + m(\fa_a,\nabla^{\sigma}_{\fa_b} \fa_i).
        \end{split}\]
    We then write
    \[
        \nabla^{\sigma} _{\fa_a} \fa_i
        = \nabla^{\delta}_{\fa_a} \fa_i + \delta(\fa_a, \fa_i)\nu
    \]
    where $\nu = x/|x|$ is the unit outward pointing normal of $\bS^{n-1}$ in
    $\bR^n$. The transversality property states that $\iota_{\nu} m = 0$. We need
    to compute $\nabla ^{\delta}_{\fa_a} \fa_i$ at $|x|=1$, that is
    \[
        \nabla^{\delta}_{\fa_a} \fa_i
        = -x^i \partial_a + 2 x^a x^i x^c \partial_c - \delta_a^i x^c \partial _c.
    \]
    From the transversality property, this yields
    \[
        m(\nabla ^{\sigma}_{\fa_a} \fa_i, \fa_b)
        = -x^i m_{ab}
        = -x^i m(\fa_a, \fa_b)
    \]
    since we restrict to the sphere $\bS^{n-1}$. Thus, we obtain
    \[
        (\lie_{\fa_i} m)(\fa_a, \fa_b)
        = (\nabla^{\sigma}_{\fa_i} m)(\fa_a, \fa_b)
        - 2 x^i m(\fa_a, \fa_b),
    \]
    for all indices $a,b$, or
    \[
        \lie_{\fa_i} m = \nabla^{\sigma}_{\fa_i} m - 2 x^i m .
    \]
    The right-hand side of~\eqref{LieXi} therefore is
    \[
        -\nabla ^{\sigma}_{\fa_i} m + k x^i m .
    \]
    which finally yields
    \[
        a_i \cdot m = -\nabla_{\fa_i}^{\sigma} m  + k x^i m ,
    \]
    as desired.

    We now derive the infinitesimal action for the one-parameter group of rotations
    $R^{\theta}_{ij}$. Since $r_{ij} = - [a_i, a_j]$, we have
    \[
        r_{ij} \cdot m = - [a_i, a_j] \cdot m.
    \]
    A straightforward computation yields
    \[ \begin{split}
            [a_i, a_j] \cdot m
             & = [ -\nabla^{\sigma}_{\fa_i} + kx^i , -\nabla^{\sigma}_{\fa_j} + kx^j] m \\
             & = \nabla^{\sigma}_{[\fa_i, \fa_j]} m + \riem^{\sigma} (\fa_i, \fa_j) m.
        \end{split} \]
    For a tangent vector $V$ to $\bS^{n-1}$ we have
    \[ \begin{split}
            \riem^{\sigma} (\fa_i , \fa_j) V
             & = \sigma(\fa_j , V) \fa_i - \sigma(\fa_i , V ) \fa_j \\
             & = dx^j (V) \fa_i -  dx^i (V) \fa_j                   \\
             & = - r_{ij} (V),
        \end{split} \]
    so
    \[
        \riem^{\sigma} (\fa_i, \fa_j) m (U,V) =  m(r_{ij}(U), V) + m(U, r_{ij}V).
    \]
    Thus we find
    \[
        (r_{ij} \cdot m)(U,V)
        = - \left(\nabla^{\sigma}_{\fr_{ij}} m (U,V)
        + m(r_{ij}(U), V) +  m(U, r_{ij}(V)) \right),
    \]
    which concludes the proof of the proposition.
\end{proof}

%

\subsection{Lie algebra intertwining operators}

We denote by $V$ an arbitrary irreducible finite-dimensional representation of
the group of isometries $O_{\uparrow}(n,1)$ of $\bH^n$. In
Lemma~\ref{lmGroupActionMassAspect}, we computed the action of the group
$O_{\uparrow}(n,1)$ on the mass aspect tensor corresponding to some jet of
metric in $\cT_k^k$. The identification we made between the space of jets
$\cT_k^k$ and the mass aspect tensors $S^2(\bS^{n-1})$ in
Lemma~\ref{lmGroupActionMassAspect} leads us to seek continuous intertwining
maps
\[
    \Fhi : S^2(\bS^{n-1}) \to V,
\]
that is, continuous maps commuting with the respective
$O_{\uparrow}(n,1)$-actions,
\[
    \forall (A,m) \in O_{\uparrow}(n,1) \times S^2(\bS^{n-1})\ ,\
    \Fhi(A \cdot m) = A \cdot (\Fhi(m)).
\]
For the associated Lie algebra representations this definition implies that
\[
    \Fhi(a \cdot m) = a \cdot (\Fhi(m)).
\]
for any element $a \in \mathfrak{so}(n,1)$ and any $m \in S^2(\bS^{n-1})$.

We define an action of the Lie algebra $\mathfrak{so}(n,1)$ on linear maps
$S^2(\bS^{n-1}) \to V$ by
\[
    (a \cdot \Fhi)(m) \definedas a \cdot (\Fhi(m)) - \Fhi(a \cdot m)
\]
for $a \in \mathfrak{so}(n,1)$ and $\Fhi : S^2(\bS^{n-1}) \to V$. Intertwining
operators are then the fixed points for this action.

Let $(v_\mu)$ be a basis of $V$ and write $\Fhi = \sum_\mu \Fhi^\mu v_\mu$.
Then the components $\Fhi^\mu$ are linear forms on $S^2(\bS^{n-1})$ which are
continuous with respect to the standard topology on $S^2(\bS^{n-1})$ from
Definition~\ref{defGeometricMass}, that is the $\Fhi^\mu$ are distributions.

As usual in distribution theory, we use the same notation for the distribution
$\Fhi^\mu$ itself and for the distribution density with values in
$S^2(\bS^{n-1})$, so that we may write
\begin{equation}\label{eqAnsatz}
    \Fhi (m) =
    \sum_{\mu} \int_{\bS^{n-1}} \langle \Fhi^\mu , m \rangle d\mu^{\sigma} v_\mu,
\end{equation}
where the symbol $\langle \cdot ,\cdot \rangle$ denotes the inner
product induced by the metric $\sigma$ on $S^2(\bS^{n-1})$.
The notation $\Fhi \definedas \sum_{\mu} \Fhi ^\mu v_\mu$ will apply for
both the distribution and the $S^2(\bS^{n-1})$-valued density.
In the following we hope that it will be clear at any place which object
we are referring to with this symbol.

By dualizing the action of $\mathfrak{so}(n,1)$ we next find conditions that
$\Fhi$ must satisfy to be an intertwining operator.

\begin{proposition}\label{propIntertwiningMap}
    If $\Fhi : S^2(\bS^{n-1}) \to V$ is intertwining for the action of
    $\mathfrak{so}(n,1)$ then
    \begin{equation}\label{eqInterBoost}
        \nabla_{\fa_i} \Fhi + (k + 1 - n) x^i \Fhi
        - \sum_\mu \Fhi^\mu a_i \cdot v_\mu = 0
    \end{equation}
    and
    \begin{equation} \label{eqInterRotation}
        \nabla_{\fr_{ij}} \Fhi + \Fhi(r_{ij}(\cdot), \cdot) + \Fhi(\cdot, r_{ij}(\cdot))
        - \sum_\mu \Fhi^\mu r_{ij} \cdot v_\mu = 0.
    \end{equation}
    Conversely, if~\eqref{eqInterBoost} and~\eqref{eqInterRotation} hold
    then the map $\Fhi : S^2(\bS^{n-1}) \to V$ is intertwining for the
    Lie algebra action.
\end{proposition}

From
\[ \begin{split}
        0 & =
        \int_{\bS^{n-1}} \divg^{\sigma} \left( \<\Fhi^\mu, m\> X\right) \, d\mu^\sigma \\
          & =
        \int_{\bS^{n-1}} \left(
        \<\nabla_X^{\sigma} \Fhi^\mu, m\> + \<\Fhi^\mu, \nabla^{\sigma}_X m\>
        + \<\Fhi^\mu, m\> \divg^{\sigma} X
        \right) \, d\mu^\sigma
    \end{split}\]
it follows that
\begin{equation} \label{eqIntByParts}
    \int_{\bS^{n-1}} \<\Fhi^\mu, \nabla^{\sigma}_X m\> \, d\mu^\sigma
    =
    - \int_{\bS^{n-1}}
    \<\nabla_X^{\sigma} \Fhi^\mu + (\divg^{\sigma} X) \Fhi^\mu,  m\>
    \, d\mu^\sigma,
\end{equation}
for any vector field $X$ on $\bS^{n-1}$. Note that
$\fa_i = \gradient^\sigma x^i$ so
$\divg^{\sigma} \fa_i = \Delta^{\sigma} x^i = -(n-1) x^i$, while
$\fr_{ij}$ is a Killing vector field of $(\bS^{n-1},\sigma)$,
and hence $\divg^{\sigma} \fr_{ij} = 0$.

\begin{proof}
    Assume $\Fhi$ is an intertwining operator and
    $\fa_i \in \mathfrak{so}(n,1)$ is a boost. Using~\eqref{X_icdotm}
    and~\eqref{eqIntByParts} we have
    \[ \begin{split}
            0
             & =
            (a_i \cdot \Fhi)(m)                                               \\
             & =
            \sum_\mu \int_{\bS^{n-1}} \<\Fhi^\mu, m\> \, d\mu^\sigma a_i \cdot v_\mu
            - \int_{\bS^{n-1}} \<\Fhi^\mu, a_i \cdot m\> \, d\mu^\sigma v_\mu \\
             & =
            \sum_\mu \int_{\bS^{n-1}} \<\Fhi^\mu, m\> \, d\mu^\sigma a_i \cdot v_\mu
            + \int_{\bS^{n-1}} \<\Fhi^\mu, \nabla^{\sigma}_{\fa_i} m - k x^i m\>
            \, d\mu^\sigma v_\mu                                              \\
             & =
            \sum_\mu \int_{\bS^{n-1}} \<\Fhi^\mu, m\> \, d\mu^\sigma a_i \cdot v_\mu
            - \int_{\bS^{n-1}}
            \<\nabla^{\sigma}_{\fa_i} \Fhi^\mu + (\divg^{\sigma} \fa_i + kx^i) \Fhi^\mu, m \>
            \, d\mu^\sigma v_\mu                                              \\
             & =
            \sum_\mu \int_{\bS^{n-1}} \<\Fhi^\mu, m\> \, d\mu^\sigma a_i \cdot v_\mu
            - \int_{\bS^{n-1}} \<\nabla^{\sigma}_{\fa_i} \Fhi^\mu - (n-1-k) x^i \Fhi^\mu, m \>
            \, d\mu^\sigma v_\mu
        \end{split} \]
    for $i = 1, \dots, n$. Since this holds for all $m \in S^2(\bS^{n-1})$ we find
    that
    \[
        \nabla_{\fa_i} \Fhi - (n - 1 - k) x^i \Fhi
        - \sum_\mu \Fhi^\mu a_i \cdot v_\mu = 0,
    \]
    where at this stage the derivatives are in the sense of distributions.

    To compute the action of a rotation $\fr_{ij} \in \mathfrak{so}(n,1)$ we need
    the following formula where $\epsilon_A$ denotes an orthonormal frame on
    $\bS^{n-1}$,
    \[ \begin{split}
            \<\Fhi, m(r_{ij} (\cdot), \cdot)\>
             & =
            \sum_{A,B} \Fhi(\epsilon_A, \epsilon_B) m (r_{ij} (\epsilon_A), \epsilon_B)   \\
             & =
            \sum_{A,B,C} \Fhi(\epsilon_A, \epsilon_B)
            m (\<r_{ij} (\epsilon_A), \epsilon_C\> \epsilon_C, \epsilon_B)                \\
             & =
            \sum_{A,B,C} \<\epsilon_A, r_{ij}^* (\epsilon_C) \>
            \Fhi(\epsilon_A, \epsilon_B) m (\epsilon_C, \epsilon_B)                       \\
             & =
            -\sum_{A,B,C} \Fhi(\<\epsilon_A, r_{ij} (\epsilon_C)\>\epsilon_A,\epsilon_B)
            m(\epsilon_C,\epsilon_B)                                                      \\
             & =
            - \sum_{B,C} \Fhi(r_{ij} (\epsilon_C), \epsilon_B) m (\epsilon_C, \epsilon_B) \\
             & =
            - \<\Fhi(r_{ij} (\cdot), \cdot), m\>,
        \end{split}\]
    and similarly,
    \[
        \<\Fhi, m(\cdot, r_{ij}(\cdot) )\> = - \<\Fhi(\cdot, r_{ij}(\cdot) ), m\>.
    \]
    Hence~\eqref{L_ijcdotm} together with~\eqref{eqIntByParts} tells us that
    \[ \begin{split}
            0 & =
            (r_{ij} \cdot \Fhi)(m)                                                      \\
              & =
            \sum_\mu \int_{\bS^{n-1}} \<\Fhi^\mu, m\> \, d\mu^\sigma r_{ij} \cdot v_\mu
            - \int_{\bS^{n-1}} \<\Fhi^\mu, r_{ij} \cdot m \> \, d\mu^\sigma v_\mu       \\
              & =
            \sum_\mu \int_{\bS^{n-1}} \<\Fhi^\mu, m\> \, d\mu^\sigma r_{ij} \cdot v_\mu \\
              & \qquad
            + \int_{\bS^{n-1}} \<\Fhi^\mu, \nabla^{\sigma}_{\fr_{ij}} m
            + m(r_{ij} (\cdot), \cdot) + m(\cdot, r_{ij} (\cdot) )
            \> \, d\mu^\sigma v_\mu                                                     \\
              & =
            \sum_\mu \int_{\bS^{n-1}} \<\Fhi^\mu, m\> \, d\mu^\sigma r_{ij} \cdot v_\mu \\
              & \qquad
            - \int_{\bS^{n-1}} \<\nabla^{\sigma}_{\fr_{ij}} \Fhi^\mu
            + \Fhi^\mu(r_{ij}(\cdot), \cdot) + \Fhi^\mu(\cdot, r_{ij}(\cdot) ), m\>
            \, d\mu^\sigma v_\mu,
        \end{split} \]
    and we conclude that
    \[
        \nabla_{\fr_{ij}} \Fhi + \Fhi(r_{ij}(\cdot), \cdot) + \Fhi(\cdot, r_{ij}(\cdot))
        - \sum_\mu \Fhi^\mu r_{ij} \cdot v_\mu = 0,
    \]
    again in the sense of distributions. Since $r_{ij} = - [a_i, a_j]$ we see
    that~\eqref{eqInterRotation} follows from~\eqref{eqInterBoost}. Since boosts
    $a_i$ and rotations $r_{ij}$ form a basis of $\mathfrak{so}(n,1)$ it is
    sufficient that~\eqref{eqInterBoost} and~\eqref{eqInterRotation} hold to
    conclude that $\Fhi$ is an intertwining map.
\end{proof}

An important remark to make at this point is the following proposition.

\begin{proposition}\label{propAnalytic}
    The distributions $\Fhi^{\mu}$ are analytic functions.
\end{proposition}

This follows from the fact that $\Fhi$ is a solution to
Equation~\eqref{eqInterBoost} which has analytic coefficients, together with
the fact that the vector fields $\fa_i$ span the tangent space $T_p \bS^{n-1}$
at each point $p$.

We end this section by looking at two simple examples of Lie algebra
intertwining maps. Both turn out to be linear masses at infinity in the sense
of Definition~\ref{defGeometricMass}, and appear in the classification
established in Section~\ref{secClassification}.

\begin{example}\label{example_trivialmass}
    As a first example, we take for $V$ the trivial $1$-dimensional
    representation of $\mathfrak{so}(n, 1)$. A basis consists of a single
    vector $v_0 \neq 0$, and the action is $a \cdot v_0 = 0$ for any
    $a \in \mathfrak{so}(n, 1)$. Set $\Fhi = \Fhi_0 v_0$. From
    Equation~\eqref{eqInterBoost} evaluated at the point where
    $\fa_i = 0$ we get that $k = n-1$. The rotations $\fr_{AB}$,
    for $2 \leq A < B \leq n$, all vanish at the south pole $p_0 = (-1,0,\dots,0)$,
    so Equation~\eqref{eqInterRotation} evaluated at this point yields
    \[
        \Fhi_0(r_{AB} (\cdot), \cdot) + \Fhi_0(\cdot, r_{AB}(\cdot)) = 0.
    \]
    This means that $\Fhi_0$ is a bilinear form which is invariant under the action
    of $\mathfrak{so}(n-1)$. The only such forms are proportional to the metric
    $\sigma$. It follows that $\Fhi = \sigma v$ for some $v \in V$ at the south
    pole. By rotational symmetry this extends to $\Fhi = \sigma v$ on all of
    $\bS^{n-1}$, where $v : \bS^{n-1} \to V$ is smooth. From~\eqref{eqInterBoost}
    it follows that $v$ is constant, and $\Fhi = \sigma v$ for some $v \in V$.

    Conversely it is clear that $\Fhi = \sigma v$ satisfies
    Equations~\eqref{eqInterBoost} and~\eqref{eqInterRotation} with $k = n-1$.
    Going back to the expression of the $\mathfrak{so}(n,1)$-intertwining map in
    this case, we obtain
    \[
        \Fhi : m \mapsto \int_{\bS^{n-1}} \tr^{\sigma} m \, d\mu ^{\sigma} \, v .
    \]
    The equivariance of such a map under the full group of hyperbolic isometries
    $O_{\uparrow}(n,1)$ holds due to Corollary~\ref{remGroupAction} for $k = n-1$.
    From Proposition~\ref{propMassaspect}, we conclude that the number
    $\int_{\bS^{n-1}} \tr ^{\sigma} m \, d\mu^{\sigma}$ is a linear mass at
    infinity. This extends to the non-Einstein case (with conformal infinity
    $(\bS^{n-1},[\sigma])$) the notion of \emph{conformal anomaly}, see for
    example~\cite[Chapter~3]{DjadliGuillarmouHerzlich}, which is defined for
    general Poincar\'e-Einstein manifolds.
\end{example}

\begin{example}\label{example_AHmassWang}
    As a second example, we choose $V = \bR^{n, 1}$ with the standard
    action of $O_{\uparrow}(n,1)$. Let
    $\partial_0, \partial_1, \dots, \partial_n$ be the
    standard orthonormal basis of $V$. For $k = n$, Wang defines in~\cite{WangMass}
    the \emph{energy-momentum vector} of $g$ as
    \begin{equation}
        \Fhi (m) =
        \int_{\bS^{n-1}} \tr^{\sigma} m \, d\mu^{\sigma} \, \partial_0 +
        \sum_i \int_{\bS^{n-1}} x^i \tr^{\sigma} m \, d\mu^{\sigma} \, \partial_i .
    \end{equation}
    Here the $x^i$'s are the coordinate functions on $\bR^n$ restricted to
    $\bS^{n-1}$. In the work of Wang as well as in the work of
    Chru\'sciel-Herzlich~\cite{WangMass,ChruscielHerzlich}, it is proven
    that the group $O_{\uparrow}(n,1)$ acts on such vectors so that $\Fhi$ is
    actually an $O_{\uparrow}(n,1)$-intertwining map, and hence a
    linear mass at infinity by Proposition~\ref{propMassaspect}.
    Further, the Minkowski norm square of this vector does not depend on
    the choice of the asymptotically hyperbolic chart. Whenever it is
    non-negative, it is the square of the so-called
    \emph{asymptotically hyperbolic mass}.
\end{example}

%% file: classification.tex
%
%

In this section we will classify all maps $\Fhi : S^2(\bS^{n-1}) \to V$ which
are intertwining for the action of the Lie group $O_{\uparrow}(n, 1)$. In
particular, such maps are $\mathfrak{so}(n,1)$-intertwining.

We first reduce the classification to a representation theoretic problem
involving only finite-dimensional representations. This is done in the
following subsections using the following steps.
\begin{itemize}
      \item
            The Lie algebra $\mathfrak{so}(n,1)$ has a parabolic subalgebra
            $\fp$ consisting of all elements for which the associated vector field
            vanishes at the south pole $p_0 = (-1,0,\ldots,0)$. This corresponds to
            the parabolic subgroup $P$ of $SO_\uparrow(n,1)$ fixing the south pole
            of the sphere at infinity. In subsection~\ref{subsec_fixingthesouthpole} we
            study the necessary conditions on $\Fhi$ coming from Equations~\eqref{eqInterBoost} and~\eqref{eqInterRotation} with vector fields fixing the south pole.
            In this case, the derivative terms in those equations vanish at the
            south pole, and we obtain a set of algebraic equations for $\Fhi(p_0)$.
            These equations are stated in Proposition~\ref{propSouthPole}.
      \item
            In Subsections~\ref{subsec_the3-dimensionalcase} and~\ref{subsecGeneral}
            we use methods from representation theory to classify all solutions to the
            equations at the south pole.
      \item
            Using the $O_{\uparrow}(n,1)$-intertwining property of $\Fhi$, in particular
            for $SO(n)$ as described in Subsection~\ref{subsec_rotations}, we will be able
            to deduce the expression of $\Fhi(x)$ at any point $x \in \bS^{n-1}$.
\end{itemize}

The conclusions of this part of the argument are collected in
Proposition~\ref{thmClassificationso(n,1)}. It remains to check that these
$\mathfrak{so}(n, 1)$-intertwining operators do yield linear masses at infinity
according to Definition~\ref{defGeometricMass}. The first step is in
Propositions~\ref{thmScalarInvariants3d},~\ref{thmWeylInvariants3d},~\ref{thmClassificationso(n,1)},
where we show that the $\mathfrak{so}(n, 1)$-intertwining operators we have
found lift to $O_{\uparrow}(n,1)$-intertwining operators on the set of mass
aspect tensors. Proposition~\ref{propMassaspect} will then ensure that we
obtain the final classification of linear masses at infinity.

\subsection{Elements fixing the south pole}\label{subsec_fixingthesouthpole}

The subalgebra $\fp$ of $\mathfrak{so}(n,1)$ fixing the south pole has a basis
consisting of
\begin{itemize}
      \item
            the infinitesimal boost $a_1$ in the $x^1$-direction,
      \item
            infinitesimal translations fixing the south pole $s_A \definedas a_A + r_{1A}$, $2 \leq A \leq n$,
      \item
            infinitesimal rotations $r_{AB}$, $2 \leq A < B \leq n$.
\end{itemize}

Let $(v_\mu)$ be a basis of the $\mathfrak{so}(n,1)$-representation $V$.

\begin{proposition}\label{propSouthPole}
      Suppose $\Fhi : S^2(\bS^{n-1}) \to V$ is an
      $\mathfrak{so}(n,1)$-intertwining map. Then
      $\Fhi(p_0) = \sum_\mu \Fhi^\mu(p_0) v_\mu$ satisfies
      \begin{equation} \label{BoostC}
            \sum_\mu \Fhi^\mu(p_0) (a_1 \cdot v_\mu)  = (n-1-k) \Fhi(p_0),
      \end{equation}
      \begin{equation} \label{TranslationC}
            \sum_\mu \Fhi^\mu(p_0) s_A \cdot v_\mu = 0,
      \end{equation}
      \begin{equation} \label{RotationC}
            \sum_\mu \Fhi^\mu(p_0) r_{AB} \cdot v_\mu
            =  \Fhi(p_0)(r_{AB}(\cdot), \cdot) + \Fhi(p_0)(\cdot, r_{AB}(\cdot)) .
      \end{equation}
\end{proposition}

\begin{proof}
      Equation~\eqref{eqInterBoost} for $a_1$ evaluated at the south pole
      $p_0 = (-1,0,\ldots,0)$ gives us
      \[ \begin{split}
                  0 & =
                  \nabla_{\fa_1} \Fhi(p_0) - (n-1-k) x^1 \Fhi(p_0)
                  - \sum_\mu \Fhi^\mu(p_0) a_1 \cdot v_\mu \\
                    & =
                  (n-1-k) \Fhi(p_0) - \sum_\mu \Fhi^\mu(p_0) a_1 \cdot v_\mu .
            \end{split} \]
      For translations at the north pole, $s_A = a_A + r_{1A}$, we get
      \[ \begin{split}
                  0 & =
                  \nabla_{\fa_A} \Fhi - (n - 1 - k) x^A \Fhi
                  - \sum_\mu \Fhi^\mu a_A \cdot v_\mu    \\
                    & \qquad
                  +\nabla_{\fr_{1A}} \Fhi + \Fhi(r_{1A}(\cdot), \cdot) + \Fhi(\cdot, r_{1A}(\cdot))
                  - \sum_\mu \Fhi^\mu r_{1A} \cdot v_\mu \\
                    & =
                  \nabla_{\mathfrak{s}_A} \Fhi - (n - 1 - k) x^A \Fhi
                  - \sum_\mu \Fhi^\mu s_A \cdot v_\mu    \\
                    & \qquad
                  + \Fhi(r_{1A}(\cdot), \cdot) + \Fhi(\cdot, r_{1A}(\cdot)) .
            \end{split} \]
      from~\eqref{eqInterBoost} and~\eqref{eqInterRotation}. At the south pole, the
      vector field $\mathfrak{s}_A \definedas \fa_A + \fr_{1A}$ vanishes and $x^A=0$.
      Further, the tangent space at the south pole is spanned by $\partial_C$, $2\leq
            C \leq n$, and
      \[
            r_{1A}(\partial_C)
            = dx^1 (\partial_C)  \partial_A - dx^A(\partial_C) \partial_1
            = -\delta_C^A \partial_1,
      \]
      so $\Fhi(r_{1A}(\cdot), \cdot) = 0$. Together we find that
      \[
            \sum_\mu \Fhi^\mu(p_0) s_A \cdot v_\mu = 0.
      \]
      Finally, Equation~\eqref{eqInterRotation} evaluated at the south pole tells us
      that
      \[ \begin{split}
                  0 & =
                  \nabla_{\fr_{AB}} \Fhi(p_0) + \Fhi(p_0)(r_{AB}(\cdot), \cdot)
                  + \Fhi(p_0)(\cdot, r_{AB}(\cdot))
                  - \sum_\mu \Fhi^\mu (p_0)r_{AB} \cdot v_\mu \\
                    & =
                  \Fhi(p_0)(r_{AB}(\cdot), \cdot) + \Fhi(p_0)(\cdot, r_{AB}(\cdot))
                  - \sum_\mu \Fhi^\mu(p_0) r_{AB} \cdot v_\mu.
            \end{split} \]
\end{proof}

\subsection{Rotations of the sphere}\label{subsec_rotations}

For $x \in \bS^{n-1}$, let $B_x$ be the rotation in the plane spanned by $x$
and $p_0$ with $B_x(p_0) = x$. Since the map $B_x$ is a rotation, the
corresponding function $u[B]$ of Lemma~\ref{lmGroupActionMassAspect} is equal
to $1$, and it therefore acts on the space of mass aspect tensors
$S^2(\bS^{n-1})$ by $B_x \cdot m = (B_x)_* m$. As a result,
\[
      \Fhi(B_x \cdot m)
      =
      \sum_{\mu} \int_{\bS^{n-1}} \< \Fhi^\mu, (B_x)_* m \> \, d\mu^{\sigma} v_\mu  .
\]
Since $B_x$ is an isometry of the sphere $\bS^{n-1}$, we have
\[ \begin{split}
            B_x \cdot \Fhi(m)
             & =
            B_x \cdot
            \left(\sum_{\mu} \int_{\bS^{n-1}} \< \Fhi^\mu, m \> \, d\mu^{\sigma} v_\mu \right) \\
             & =
            \sum_{\mu} \int_{\bS^{n-1}} \< (B_x)_*\Fhi^\mu, (B_x)_*m \> \, d\mu^{\sigma}
            \left(B_x \cdot v_\mu \right).
      \end{split} \]
The $O_{\uparrow}(n,1)$-intertwining property for $\Fhi$ tells us that
\[
      \sum_{\mu} \Fhi^\mu v_\mu
      = \sum_{\mu} (B_x)_*\Fhi^\mu \left( B_x \cdot v_\mu \right)
\]
at all points. If we in particular evaluate at the point $x$ we get
\[
      \Fhi(x) (U,V) =
      \sum_{\mu} \Fhi^\mu (p_0) ((B_x)^* U , (B_x)^* V) \left( B_x \cdot v_\mu \right)
\]
for $U,V \in T_x \bS^{n-1}$. This last equation gives a way to transport the
expression of $\Fhi$ we will obtain at the south pole $p_0$ to any other point
of $\bS^{n-1}$, thanks to the $O_{\uparrow}(n,1)$-intertwining property of
$\Fhi$.

\subsection{The 3-dimensional case}\label{subsec_the3-dimensionalcase}

Here we specialize to the case where the dimension is $n=3$, and assume that
$\Fhi$ is an $O_{\uparrow}(3,1)$-intertwining map. To classify all
possibilities for $\Fhi$ we start by finding all solutions to the equations at
the south pole, this is Proposition~\ref{propSouthPole}.

The approach we use is specific to the 3-dimensional case, and the resulting
classification is also partly specific, compare
Proposition~\ref{thmWeylInvariants3d} with the higher-dimensional version
stated in Proposition~\ref{thmClassificationso(n,1)}. The methods used in this
section are rather elementary and might serve as a warm-up for the case where
the dimension $n \geq 4$, which is treated in Subsection~\ref{subsecGeneral}.

It is a standard fact that the complexified Lie algebra $\mathfrak{so}(3, 1)
      \otimes \bC$ can be written as the sum of two copies of $\mathfrak{sl}(2,\bC)$.
Several such splittings exist since $\mathfrak{sl}(2,\bC)$ has many
automorphisms. We implement a splitting by choosing generators
\[
      \left\lbrace
      \begin{array}{rl}
            h_1 & = -a_1 - i r_{23}                                          \\
            e_1 & = \frac{1}{2}\left(-a_2 + i a_3 - r_{12} + i r_{13}\right) \\
            f_1 & = \frac{1}{2}\left(-a_2 - i a_3 + r_{12} + i r_{13}\right)
      \end{array}
      \right. ,
      \quad
      \left\lbrace
      \begin{array}{rl}
            h_2 & = -a_1 + i r_{23}                                          \\
            e_2 & = \frac{1}{2}\left(-a_2 - i a_3 - r_{12} - i r_{13}\right) \\
            f_2 & = \frac{1}{2}\left(-a_2 + i a_3 + r_{12} - i r_{13}\right)
      \end{array}
      \right. .
\]
The families $\{h_1, e_1, f_1\}$ and $\{h_2, e_2, f_2\}$ each satisfy the
commutation relations of $\mathfrak{sl}(2,\bC)$, and they commute with each
other. Complex irreducible representations of the Lie algebra
$\mathfrak{so}(3,1)$ are therefore in bijection with tensor product
representations $V^1 \otimes V^2$ where $V^1$ (resp. $V^2$) is an irreducible
representation of $(\mathfrak{sl}_2)_1 = \mathrm{vect}(h_1, e_1, f_1)$ (resp.
$(\mathfrak{sl}_2)_2 = \mathrm{vect}(h_2, e_2, f_2)$),
see~\cite[Section~4.2.2]{GoodmanWallach}. Let $\mathfrak{h}$ denote the Cartan
subalgebra of $\mathfrak{so}(3, 1)$ generated by $h_1$ and $h_2$ and let
$(\omegabar_1, \omegabar_2)$ be the basis of $\mathfrak{h}^*$ dual to
$(h_1,h_2)$. The highest-weight theory tells us that irreducible
finite-dimensional representations of $\mathfrak{so}(3, 1)$ are in bijection
with the set $\bN \omegabar_1 + \bN \omegabar_2$.

Suppose that the irreducible representation $V$ has highest weight $w_1
      \omegabar_1 + w_2 \omegabar_2$. In order to characterize $\Fhi(p_0)$ when $\Fhi
      : S^2(\bS^{2}) \to V$ is a $\mathfrak{so}(3,1)$-intertwining map, we write
\[
      s_2 = a_2 + r_{12} = -e_1 - e_2, \quad s_3 = a_3 + r_{13} = -i\left(e_1 - e_2\right).
\]
From Equation~\eqref{TranslationC} we find
\[
      e_1 \cdot (\Fhi(p_0)) = e_2 \cdot (\Fhi(p_0)) = 0,
\]
which means that
\begin{equation}\label{PhiHighestWeight}
      \Fhi(p_0) = v_{w_1 \omegabar_1 + w_2 \omegabar_2} \fhi,
\end{equation}
where $v_{w_1 \omegabar_1 + w_2 \omegabar_2}$ is a highest weight vector of the representation $V = V_{w_1 \omegabar_1 + w_2 \omegabar_2}$ and $\fhi$ is a symmetric 2-tensor over $T_{p_0}\bS^2$. Further, Equation~\eqref{BoostC} tells us that
\[
      a_1 \cdot (\Fhi(p_0)) = (2-k) \Fhi(p_0),
\]
where $a_1$ is equal to $-\frac12 (h_1 + h_2)$. Therefore we have $a_1 \cdot
      (\Fhi(p_0)) = -\frac{w_1 + w_2}{2} \Fhi(p_0)$, which gives the expression
\[
      k = 2 + \frac{w_1 + w_2}{2}
\]
for the decay order $k$.

Next we use Equation~\eqref{RotationC} together with the relation $r_{23} =
      \frac{i}{2}(h_1 - h_2)$. At the south pole where $y=0$ we have
$r_{23}(\partial_2) = \partial_3$ and $r_{23}(\partial _3) = -\partial_2$.
Evaluating~\eqref{RotationC} on all pairs of vectors tangent to the sphere
$\bS^2$ at the south pole, we obtain the system of equations
\[
      \left\lbrace
      \begin{aligned}
            \frac{i}{2} (w_1 - w_2) \fhi_{22} & = 2 \fhi_{23},           \\
            \frac{i}{2} (w_1 - w_2) \fhi_{23} & = \fhi_{33} - \fhi_{22}, \\
            \frac{i}{2} (w_1 - w_2) \fhi_{33} & = -2 \fhi_{23}.
      \end{aligned}
      \right.
\]
It is a simple exercise to solve this eigenvalue problem. Up to a
multiplicative constant, we find three solutions,
\begin{itemize}
      \item $w_1 = w_2$, with
            $\fhi = dx^2 \otimes dx^2 + dx^3 \otimes dx^3
                  = \frac{1}{2}(dz \otimes d\bar{z} + d\bar{z} \otimes dz)$,
      \item $w_1 = w_2 + 4$, with
            $\fhi = dx^2 \otimes dx^2 - dx^3 \otimes dx^3
                  + i (dx^2 \otimes dx^3 + dx^3 \otimes dx^2)
                  = dz \otimes dz$, and
      \item $w_1 = w_2 - 4$, with
            $\fhi = dx^2 \otimes dx^2 - dx^3 \otimes dx^3
                  - i (dx^2 \otimes dx^3 + dx^3 \otimes dx^2)
                  = d\bar{z} \otimes d\bar{z}$,
\end{itemize}
where we have set $z \definedas x^2 + i x^3$.

For $n_1 \geq 0$ let $\cH_{n_1}$ denote the space of wave-harmonic homogeneous
polynomials of degree $n_1$ on $\bR^{3, 1}$. This is a representation of
$O_\uparrow(3,1)$ under
\[
      \begin{array}{ccc}
            O_\uparrow(3,1) \times \cH_{n_1} & \to     & \cH_{n_1}       \\
            (A,P)                            & \mapsto & P \circ A^{-1}.
      \end{array}
\]

In the first of these three cases, we have $w_1 = w_2 = n_1 = k-2$. We get the
following family of $O_{\uparrow}(3,1)$-intertwining maps and corresponding
linear masses at infinity.

\begin{proposition}[Harmonic masses]\label{thmScalarInvariants3d}
      Let $\Fhi : S^2(\bS^2) \to V$ be an $O_{\uparrow}(3,1)$-intertwining map, such that $V$ as a $\mathfrak{so}(3,1)$-representation has highest weight $n_1 \omegabar_1 + n_1 \omegabar_2$, with $n_1 \geq 0$. The representation $V$ can then be identified with the dual of $\cH_{n_1}$, and the map $\Fhi$ is a multiple of the map
      \[
            \begin{array}{rccl}
                  \Fhi_h (m) : & \cH_{n_1} & \to     & \bR \\
                               & P         & \mapsto &
                  \displaystyle{\int_{\bS^2} P(1, x^1,x^2,x^3) \tr^\sigma m
                        \, d\mu^\sigma},
            \end{array}
      \]
      for $m \in S^2(\bS^2)$. Further, this map is $O_{\uparrow}(3,1)$-intertwining,
      hence a linear mass at infinity for asymptotically hyperbolic metrics of order
      $k = n_1 + 2$.
\end{proposition}

Here the linear masses from Examples~\ref{example_trivialmass}
and~\ref{example_AHmassWang} appear as the cases $n_1=0$ and $n_1=1$.

\begin{proof}
      The representation with highest weight $n_1 (\omegabar_1 + \omegabar_2)$ corresponding to the case $w_1 = w_2 = n_1 = k-2$ can be realized as the subrepresentation
      $\mathring{\Sym}^{n_1}(\bR^{3,1})$ of $(\bR^{3,1})^{\otimes n_1}$ consisting of harmonic (that is, trace-free) symmetric $n_1$-tensors. The highest weight vector $v_{n_1 (\omegabar_1 + \omegabar_2)}$ is given by $(\partial_0 - \partial_1)^{\otimes n_1}$. From~\eqref{PhiHighestWeight}, we know that $\Fhi(p_0)$ can be written up to a constant as $(\partial_0 - \partial_1)^{\otimes n_1} \sigma(p_0)$.

      The vector $\partial_0 - \partial_1$ is the position vector of the south pole
      when we view the sphere $\bS^2$ as the submanifold
      \begin{equation}\label{SphereMinkowski}
            \left\{ (1, x^1, x^2, x^3) \mid
            (x^1)^2 + (x^2)^2 + (x^3)^2 = 1 \right\} \subset \bR^{3,1}
      \end{equation}
      Using the subgroup $SO(3) \subset SO_\uparrow(3, 1)$ stabilizing the vector $\partial_0$ to translate the formula for $\Fhi(p_0)$ to any given point on $\bS^2$ (see Subsection~\ref{subsec_rotations}), we obtain
      \[
            \Fhi(x^1, x^2, x^3) =
            (\partial_0 + x^1 \partial_1 + x^2 \partial_2 + x^3 \partial_3)^{\otimes n_1}
            \otimes \sigma (x^1,x^2,x^3).
      \]

      Let $\bR_{n_1}[X^0, X^1, X^2, X^3]$ be the space of homogeneous polynomials of
      degree $n_1$ in $X^0, X^1, X^2, X^3$. Elements of the basis $(\partial_\mu)$ of
      $\bR^{3, 1}$ act as derivations
      \[
            \partial_\mu:  \bR_{n_1}[X^0, X^1, X^2, X^3]
            \to \bR_{n_1 - 1}[X^0, X^1, X^2, X^3].
      \]
      The operator $x^\mu \partial_\mu$ is a polarization operator on
      $\bR_{n_1}[X^0,X^1, X^2, X^3]$, see~\cite[Section~2]{Procesi}. It follows that
      for any polynomial $P \in \bR_{n_1}[X^0, X^1, X^2, X^3]$, we have
      \[
            \left(
            \partial_0 + x^1 \partial_1 + x^2 \partial_2 + x^3 \partial_3
            \right)^{n_1} P
            = (n_1)! P(1, x^1, x^2, x^3).
      \]
      This way $\Fhi(x^1, x^2, x^3)$ appears as an element of $(\bR_{n_1}[X^0,
                  X^1,X^2, X^3])^* \otimes \sigma$, namely
      \[
            \Fhi(m)(P) = (n_1)!
            \int_{\bS^2} P(1, x^1, x^2, x^3) \tr^{\sigma} m \, d\mu^{\sigma}.
      \]
      Notice that since $(x^1)^2 + (x^2)^2 + (x^3)^2 = 1$ on $\bS^2$, if
      \[
            P = \left((X^0)^2 - (X^1)^2 - (X^2)^2 - (X^3)^2\right) Q,
      \]
      we have
      \[ \begin{split}
                  \Fhi(m)(P)
                   & = (n_1)!
                  \int_{\bS^2} \left(1 - (x^1)^2 - (x^2)^2 - (x^3)^2\right)
                  Q(1, x^1, x^2, x^3) \tr^{\sigma} m \, d\mu^{\sigma} \\
                   & = 0.
            \end{split} \]
      From the decomposition
      \[ \begin{split}
                   & \bR_{n_1}[X^0, X^1, X^2, X^3] \\
                   & \qquad
                  = \cH_{n_1} \oplus \left((X^0)^2 - (X^1)^2 - (X^2)^2 - (X^3)^2\right)
                  \bR_{n_1 - 2}[X^0, X^1, X^2, X^3]
            \end{split} \]
      it follows that $\Fhi(m)$ is an element of $(\cH_{n_1})^*$.

      The argument to check that the maps $\Fhi_h : S^2(\bS^2) \to (\cH_{n_1})^*$ are
      $O_{\uparrow}(3,1)$-intertwining is the same as in the higher-dimensional case,
      see the proof of Proposition~\ref{thmClassificationso(n,1)}.
\end{proof}

For the cases $w_1 = w_2 \pm 4$, the situation is more complicated. Each
representation $V_{(n_1+4) \omegabar_1 + n_1 \omegabar_2}$ and $V_{n_1
                  \omegabar_1 + (n_1+4) \omegabar_2}$ is a complex representation. But their sum
$V_{(n_1+4) \omegabar_1 + n_1 \omegabar_2} \oplus V_{n_1 \omegabar_1 + (n_1+4)
                  \omegabar_2}$ is a real representation. From the solution of the Clebsch-Gordan
problem (see e.g.~\cite[Chapter X]{CohenTannoudji}) it follows that this
combination of representations appears in the decomposition of
\begin{equation}\label{eqDecompositionRep}
      V_{n_1 (\omegabar_1 + \omegabar_2)} \otimes \left(V_{4 \omegabar_1} \oplus V_{4 \omegabar_2}\right).
\end{equation}

The representation $V_{2 \omegabar_1} \oplus V_{2 \omegabar_2}$ can be
identified with the representation $\Lambda^2 \bR^{3, 1}$ which has a basis
given by $\partial_\mu \wedge \partial_\nu$, where $0 \leq \mu < \nu \leq 3$.
The highest weight vectors of this representation can be computed explicitly.
They are $v_{2 \omegabar_1} = (\partial_0 - \partial_1) \wedge (\partial_2 -
      i\partial_3)$ for $V_{2 \omegabar_1}$ and $v_{2 \omegabar_2} = (\partial_0 -
      \partial_1) \wedge (\partial_2 + i \partial_3)$ for $V_{2 \omegabar_2}$.

Elements in the representation $V_{4 \omegabar_1}$ (resp. $V_{4 \omegabar_2}$)
can be understood as symmetric tensor products of elements in
$V_{2\omegabar_1}$ (resp. $V_{2 \omegabar_2}$). In particular, the highest
weight vector of the representation $V_{4 \omegabar_1}$ (resp. $V_{4
                  \omegabar_2}$) is given by
\[
      v_{4 \omegabar_1} = \left((\partial_0 - \partial_1) \wedge
      (\partial_2 - i \partial_3)\right)^{\otimes 2}
      \quad \text{resp.} \quad
      v_{4 \omegabar_2} = \left((\partial_0 - \partial_1) \wedge
      (\partial_2 + i \partial_3)\right)^{\otimes 2}.
\]
It should be noted at this point that these tensors belong to $\Sym^2(\Lambda^2
      \bC^4)$ and hence can be written in component notation as
\[
      W^{\mu\nu\alpha\beta}
      \partial_\mu \otimes \partial_\nu \otimes \partial_\alpha \otimes \partial_\beta.
\]
It is straightforward to check that
\[
      \begin{aligned}
            0 & = W^{\mu\nu\alpha\beta} \eta_{\mu\alpha},                                \\
            0 & = W^{\mu\nu\alpha\beta} + W^{\nu\alpha\mu\beta} + W^{\alpha\mu\nu\beta},
      \end{aligned}
\]
which means that these vectors are trace-free with respect to the Minkowski
metric $\eta$ and satisfy the contravariant version of the first Bianchi
identity.

Using~\eqref{PhiHighestWeight} we will now look for the form of $\Fhi(p_0)$.
Note that, introducing the complex structure $J$ on $\bS^2$, we have
\[
      \begin{split}
            (\partial_2 - i \partial_3) dz
             & = \partial_2 dx^2 + \partial_3 dx^3 - i (\partial_3 dx^2 - \partial_2 dx^3)
            = \mathrm{Id} - i J,                                                           \\
            (\partial_2 + i \partial_3) d\zbar
             & = \partial_2 dx^2 + \partial_3 dx^3 + i (\partial_3 dx^2 - \partial_2 dx^3)
            = \mathrm{Id} + i J,
      \end{split}
\]
at the south pole. As a consequence, the element $v_{4 \omegabar_1} \otimes
      (dz\otimes dz)$ (resp. $v_{4 \omegabar_2} \otimes (d\zbar\otimes d\zbar)$) can
be understood as a map
\[
      \begin{array}{rccl}
            \Psi_+: & T_{p_0} \bS^2 \otimes T_{p_0} \bS^2 & \to     &
            \Sym^2(\Lambda^2 \bC^4)                                   \\
                    & X \otimes Y                         & \mapsto &
            \left((\partial_0-\partial_1) \wedge (X - i J(X))\right) \otimes
            \left((\partial_0-\partial_1) \wedge (Y - i J(Y))\right),
      \end{array}
\]
resp. a map
\[
      \begin{array}{rccl}
            \Psi_-: & T_{p_0} \bS^2 \otimes T_{p_0} \bS^2 & \to     &
            \Sym^2(\Lambda^2 \bC^4)                                                                                            \\
                    & X \otimes Y                         & \mapsto & \left((\partial_0-\partial_1) \wedge (X + i J(X))\right)
            \otimes \left((\partial_0-\partial_1) \wedge (Y + i J(Y))\right),
      \end{array}
\]
where $T_{p_0} \bS^2$ denotes the tangent space at the south pole of $\bS^2$,
when $\bS^2$ is seen as the $2$-sphere embedded in $\bR^{3,1}$ as
in~\eqref{SphereMinkowski}. This formula can be written using the Hodge star
operator $\star$ acting on $\Lambda^2\bR^{3, 1} \otimes \bC$, see for
example~\cite[Definition~1.51]{Besse}. Straightforward computations lead to
\[
      \star (e_+ \wedge X) = e_+ \wedge J(X)
\]
for any vector $X \in T_{p_0} \bS^2$, where $e_+$ is the future pointing null
vector at $p_0$ defined by $e_+ \definedas \partial_0 - \partial_1$. Further,
$\star^2 = - \mathrm{Id}$ so the eigenvalues of $\star$ on $\Lambda^2\bR^{3,
            1}$ are $\pm i$.

\begin{definition}\label{defEigenspaces}
      We denote by $\left(\Lambda^2\bR^{3, 1}\right)^\pm$ the eigenspaces of the Hodge star operator on the complexified space
      $\Lambda^2\bR^{3, 1} \otimes \bC$, that is
      \[
            \left(\Lambda^2\bR^{3, 1}\right)^\pm
            =
            \left\{ \omega \in \Lambda^2\bR^{3, 1} \otimes \bC \mid
            \star \omega = \pm i \omega  \right\}.
      \]
      We define $p_\pm$ as the projection operator onto $\left(\Lambda^2\bR^{3,
                  1}\right)^\pm$ with kernel $\left(\Lambda^2\bR^{3, 1}\right)^\mp$, that is
      $p_{\pm}(\omega) = \frac{1}{2}(\omega \mp i \star \omega)$.
\end{definition}

The decomposition
\[
      \Lambda^2\bR^{3, 1} \otimes \bC =
      \left(\Lambda^2\bR^{3, 1}\right)^+
      \oplus \left(\Lambda^2\bR^{3, 1}\right)^-
\]
corresponds to the decomposition into irreducible representations under the
action of $\mathfrak{so}(3, 1)$, where $\left(\Lambda^2\bR^{3, 1}\right)^+$
corresponds to the representation $V_{2 \omegabar_1}$ and $\left(\Lambda^2
      \bR^{3, 1} \right)^-$ corresponds to $V_{2 \omegabar_2}$.

The operators $\Psi_\pm$ are then better understood as
\[
      \left\lbrace
      \begin{aligned}
            \Psi_+ (X, Y) & = p_+ (e_+ \wedge X) \otimes p_+ (e_+ \wedge Y), \\
            \Psi_- (X, Y) & = p_- (e_+ \wedge X) \otimes p_- (e_+ \wedge Y),
      \end{aligned}
      \right.
\]
where we have dropped an irrelevant factor $4$ in the definition of $\Psi_\pm$.

Let $\cW_0$ denote the space of constant Weyl tensors, namely the set of
covariant 4-tensors $W \in \Sym^2(\Lambda^2 \bR^{3, 1})$ satisfying
\[
      \eta^{\mu\alpha} W_{\mu\nu\alpha\beta} = 0,
      \quad W_{\mu\nu\alpha\beta} + W_{\nu\alpha\mu\beta} + W_{\alpha\mu\nu\beta} = 0.
\]
Elements of $\cW_0$ can be seen as endomorphisms of $\Lambda^2\bR^{3,1}$
commuting with $\star$, see~\cite[Paragraphs 3.19 and 3.20]{Besse}. The space
$\cW_0$ naturally pairs with $\Sym^2(\Lambda^2 \bR^{3, 1})$, and since the
Hodge star operator is self-adjoint we get
\[ \begin{split}
            \left\< W, \alpha \otimes \star \beta\right\>
             & = \left\< \alpha, W(\star \beta) \right\>       \\
             & = \left\< \alpha, \star W(\beta) \right\>       \\
             & = \left\< \star \alpha, W(\beta) \right\>       \\
             & = \left\< W, \star \alpha \otimes \beta\right\>
      \end{split} \]
for any $W \in \cW_0$ and any $\alpha, \beta \in \Lambda^2\bR^{3,1}$. Pairing
the elements of $\cW_0$ with $\Psi_\pm$, we get
\[
      \<W, \Psi_\pm(X, Y)\> = \<W, p_\pm(e_+ \wedge X) \otimes (e_+ \wedge Y)\>,
\]
where one projection $p_\pm$ disappears due to the previous calculations.

Returning to the representations $V_{(n_1+4) \omegabar_1 + n_1 \omegabar_2}$
(resp. $V_{n_1 \omegabar_1 + (n_1+4) \omegabar_2}$), the element
$v_{(n_1+4)\omegabar_1 + n_1 \omegabar_2} dz \otimes dz$, resp. $v_{n_1
                  \omegabar_1 + (n_1+4) \omegabar_2} d\overline{z} \otimes d\overline{z}$, can be
written as
\[
      \Fhi_\pm(p_0) =
      (e_+)^{\otimes n_1} \otimes
      \big(p_\pm (e_+ \wedge \cdot ) \otimes p_\pm (e_+ \wedge \cdot) \big),
\]
where the factor $(e_+)^{\otimes n_1}$ comes from the
$V_{n_1(\omegabar_1+\omegabar_2)}$-part of~\eqref{eqDecompositionRep}. Let
$\cW_{n_1}$ denote the set of \emph{polynomial Weyl tensors} of degree $n_1$,
that is
\[
      \cW_{n_1} \definedas
      \left\{ W \in \bR_{n_1}[X^0,X^1,X^2,X^3] \otimes \cW_0 \mid
      \partial_\mu W_{\nu\sigma\alpha\beta} + \partial_\nu W_{\sigma\mu\alpha\beta}
      + \partial_\sigma W_{\mu\nu\alpha\beta} = 0 \right\}.
\]
We extend $e_+$ by rotations to the null vector field
\[
      e_+ \definedas x^{\mu} \partial_{\mu}
      = \partial_0 + x^1 \partial_1 + x^2 \partial_2 + x^3 \partial_3
\]
on $\bS^2$. By invariance under rotations, and arguing as in the proof of
Proposition~\ref{thmScalarInvariants3d}, we have for all $W \in \cW_{n_1}$,
\[
      \< W,\Fhi_{\pm}\>  =
      (n_1!)\ W(e_+,(\mathrm{Id} \mp i J)(\cdot), e_+, \cdot)(1,x^1,x^2,x^3)
\]
where the right-hand side is an element of $\bR_{n_1}[X^0,X^1,X^2,X^3] \otimes
      S^2(\bS^{2})$ evaluated at the point $(X^0,X^1,X^2,X^3) = (1,x^1,x^2,x^3)$. In
the proof of Proposition~\ref{thmClassificationso(n,1)}, we will see that the
maps $\Fhi_{w, \pm}$ are $O_\uparrow(n, 1)$-intertwining. Hence, we have the
following proposition.

\begin{proposition}[Weyl masses]\label{thmWeylInvariants3d}
      Let $\Fhi : S^2(\bS^2) \to V$ be a $O_{\uparrow}(3,1)$-intertwining map, such that $V$, as a $\mathfrak{so}(3,1)$-representation has highest weight $(n_1+4) \omegabar_1 + n_1 \omegabar_2$ (resp. $n_1 \omegabar_1 +(n_1+4) \omegabar_2$), for $n_1 \geq 0$. The representation $V$ can then be identified with a subspace of the dual of $\bC \otimes \cW_{n_1}$, and the map $\Fhi$ is a multiple of the map defined by
      \[
            \begin{array}{rccl}
                  \Fhi_{w,+}(m): & \cW_{n_1} & \to     & \bC                                                                                                   \\
                                 & W         & \mapsto & \int_{\bS^{2}} \left\<m, W(e_+,(\mathrm{Id} - i J)(\cdot), e_+, \cdot)\right\>_\sigma \, d\mu^\sigma,
            \end{array}
      \]
      resp.
      \[
            \begin{array}{rccl}
                  \Fhi_{w,-}(m): & \cW_{n_1} & \to     & \bC \\
                                 & W         & \mapsto &
                  \int_{\bS^{2}} \left\<m, W(e_+,(\mathrm{Id} + i J)(\cdot), e_+, \cdot)\right\>_\sigma \, d\mu^\sigma
            \end{array}
      \]
      for $m \in S^2(\bS^2)$. Further, these maps are
      $O_{\uparrow}(3,1)$-intertwining, hence linear masses at infinity for
      asymptotically hyperbolic metrics of order $k = n_1 + 4$.
\end{proposition}

Note that the masses at infinity found in this proposition depend only on the
trace-free part of the mass aspect tensor $m$.

\subsection{The case of dimension \texorpdfstring{$n \geq 4$}{TEXT}}\label{subsecGeneral}

In this section we assume $n \geq 4$. For the representation theory we follow
the conventions of~\cite{GoodmanWallach}.

Recall that $(\partial_\mu)_{\mu=0,\ldots,n}$ denotes the standard basis of
$\bR^{n, 1}$. We introduce a basis $(e_i)$ of $\bR^{n, 1} \otimes \bC$ as
follows.
\begin{itemize}
      \item
            If $n+1 = 2l$ is even, we set $e_{+1} \definedas \partial_0 - \partial_1$ and $e_{-1} \definedas - \frac{\partial_0 + \partial_1}{2}$. For $k = 2, \ldots, l$, we set $e_{+k} \definedas \partial_{2k-2} + i\partial_{2k-1}$ and $e_{-k} \definedas \frac{1}{2} (\partial_{2k-2} - i\partial_{2k-1})$.
      \item
            If $n+1 = 2l+1$ is odd, we set $e_0 \definedas \partial_n$ and define $e_{\pm k}$, $k = 1, \ldots, l$ as in the previous case.
\end{itemize}
This basis is chosen in such a way that $\eta(e_i,e_j)=\delta_{i,-j}$. We also define the dual basis $(e^i)$ with respect to the Minkowski metric, so that $e^i=\eta(e_i,\cdot)$. We consider the maximal torus $H$ of $SO_{\bC}(n, 1)$ which is the subgroup of matrices which are diagonal in the basis $(e_i)$. We denote its Lie algebra by $\mathfrak{h}$.

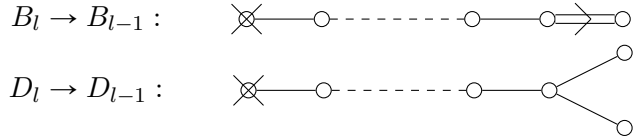
\begin{figure}[ht]
      \begin{tikzpicture}
            \draw (-1,0) node[anchor=east]  {$B_l \to B_{l-1}:$};

            \node[dnode] (1) at (0,0) {};
            \node[dnode] (2) at (1,0) {};
            \node[dnode] (3) at (3,0) {};
            \node[dnode] (4) at (4,0) {};
            \node[dnode] (5) at (5,0) {};

            \path (1) edge[sedge] (2)
            (2) edge[sedge,dashed] (3)
            (3) edge[sedge] (4)
            (4) edge[dedge] (5)
            ;

            \draw (-0.2, -0.2) -- (0.2,  0.2);
            \draw (-0.2,  0.2) -- (0.2, -0.2);
      \end{tikzpicture}

      \begin{tikzpicture}
            \draw (-1,0) node[anchor=east]  {$D_l \to D_{l-1}:$};

            \node[dnode] (1) at (0,0) {};
            \node[dnode] (2) at (1,0) {};
            \node[dnode] (3) at (3,0) {};
            \node[dnode] (4) at (4,0) {};
            \node[dnode] (5) at (5,0.5) {};
            \node[dnode] (6) at (5,-0.5) {};

            \path (1) edge[sedge] (2)
            (2) edge[sedge,dashed] (3)
            (3) edge[sedge] (4)
            (4) edge[sedge] (5)
            edge[sedge] (6)
            ;

            \draw (-0.2, -0.2) -- (0.2,  0.2);
            \draw (-0.2,  0.2) -- (0.2, -0.2);
      \end{tikzpicture}
      \caption{Reduction from $\mathfrak{so}(n, 1)$ to $\mathfrak{so}(n-1)$
            in terms of Dynkin diagrams.}\label{figDynkin}
\end{figure}

This convention is convenient since reduction from $\mathfrak{so}(n,1)$ to the
subalgebra $\mathfrak{so}(n-1)$ of orientation preserving isometries at the
south pole corresponds to deleting the leftmost node in the Dynkin diagrams
$B_l$ (if $n+1=2l+1$ is odd) or $D_l$ (if $n+1=2l$ is even). See
Figure~\ref{figDynkin}. The number of nodes in these diagrams correspond to
$l$, the rank of the Lie algebra $\mathfrak{so}(n,1)$.

Notice that for the diagram $B_l$ with $l=2$ and for $D_l$ with $l=3$ (so when
$n=4$ or $5$), the shape of the Dynkin diagram changes. This means that these
cases will require some special attention.

The elements $s_A$ introduced right before Proposition~\ref{propSouthPole} in
the 3-dimensional case generalize to the ladder operator
\begin{equation}\label{eqTranslationSouthPole}
      \left\lbrace
      \begin{aligned}
            X_{\epsilon_1-\epsilon_k} & = \frac{1}{2} \left(s_{2k-2} + i s_{2k-1}\right), \\
            X_{\epsilon_1+\epsilon_k} & = s_{2k-2} - i s_{2k-1},                          \\
            X_{\epsilon_1}            & = s_n \qquad \text{(for odd $n+1$),}
      \end{aligned}
      \right.
\end{equation}
for $k=-l, \ldots, -1, 1, \ldots, l$, see~\cite[Section~2.4.1]{GoodmanWallach}. Here $\epsilon_k \in \mathfrak{h}^*$ is the linear functional whose value on a matrix is the $k$-th term on its diagonal. The notation $X_{f}$ where $f \in \mathfrak{h}^*$ means that $\mathrm{ad}(h)(X_f) = f(h) X_f$, that is $X_f$ belongs to the root space associated to $f$.

The fundamental weights $\omegabar_i$ of $\mathfrak{so}(n, 1)$ are given by
\begin{itemize}
      \item
            ($n+1$ odd) $\omegabar_i = \epsilon_1 + \cdots + \epsilon_i$ for
            $1 \leq i \leq l-1$ and
            $\omegabar_l = \frac{1}{2} \left(\epsilon_1 + \cdots + \epsilon_l\right)$,
      \item
            ($n+1$ even) $\omegabar_i = \epsilon_1 + \cdots + \epsilon_i$ for $1 \leq i \leq l-2$ and
            $\omegabar_{l-1} = \frac{1}{2} \left(\epsilon_1 + \cdots + \epsilon_{l-1} - \epsilon_l\right)$,
            $\omegabar_l = \frac{1}{2} \left(\epsilon_1 + \cdots + \epsilon_{l-1} + \epsilon_l\right)$,
\end{itemize}
see~\cite[Section~3.1.3]{GoodmanWallach}.

We let $E_i$ be the matrix having coefficient $1$ on the $i$-th diagonal
position and zeros elsewhere, where the diagonal positions are enumerated by
$-l \leq i \leq l$. The coroots associated to the fundamental weights are given
by
\begin{itemize}
      \item
            ($n+1$ odd) $H_i = E_i - E_{i+1} + E_{-i-1} - E_{-i}$ for
            $1 \leq i \leq l-1$ and $H_l = 2(E_l - E_{-l})$,
      \item
            ($n+1$ even) $H_i = E_i - E_{i+1} + E_{-i-1} - E_{-i}$ for
            $1 \leq i \leq l-1$ and $H_l = E_{l-1} + E_l - E_{-l} - E_{-(l-1)}$.
\end{itemize}

The subalgebra $\mathfrak{so}(n-1)$ then has a Cartan subalgebra generated by
$H_2, \ldots, H_l$. Hence, the fundamental weights of $\mathfrak{so}(n-1)$ are
$\omegabar_i$ for $i \geq 2$ except in the case $n+1=5$ for which the
fundamental weight is $\omegabar = \epsilon_l$ which is twice the restriction
of $\omegabar_2$.

We now assume that $\Fhi : S^2(\bS^{n-1}) \to V$ is an
$O_{\uparrow}(n,1)$-intertwining map. Formula~\eqref{RotationC} tells us that
$\Fhi(p_0)$ intertwines the actions of $\mathfrak{so}(n-1)$ on
$\Sym^2(T_{p_0}^* \bS^{n-1})$ and on $V$. Hence, the preimage of a highest
weight vector $v$ of $V$ is a highest weight vector in $\Sym^2 (T_{p_0}^*
      \bS^{n-1})$ for the action of $\mathfrak{so}(n-1)$.

We assume for the moment that $n \geq 6$, we will later indicate the
modifications needed for $n=4,5$. The space $\Sym^2 (T_{p_0}^* \bS^{n-1})$
decomposes into two irreducible representations for the action of
$\mathfrak{so}(n-1)$,
\begin{itemize}
      \item
            the trivial one-dimensional representation generated by $\sum_{A = 2}^n dx^A \otimes dx^A = \sum_{j \neq \pm 1} e^j \otimes e^{-j}$,
      \item
            the set $\Sym^2_0(T_{p_0}^* \bS^{n-1})$ of symmetric trace-free $2$-tensors with highest weight vector $e^2 \otimes e^2$.
\end{itemize}
The corresponding highest weights are $0$ for the trivial representation and $2 \omegabar_2$ for $\Sym^2_0(T_{p_0}^* \bS^{n-1})$.

From the formulas in Equation~\eqref{eqTranslationSouthPole},
Condition~\eqref{TranslationC} implies that if $v_0$ is one of the highest
weight vectors given above, then $\Fhi(p_0)(v_0)$ is annihilated by the full
nilpotent algebra $\mathfrak{n}_+$ consisting of positive roots of
$\mathfrak{so}(n, 1)$. In particular, up to a normalization constant, $v =
      \Fhi(p_0)(v_0)$.

Let $\omega$ denote the highest weight of $V$, that is $h \cdot v = \omega(h)
      v$ for all $h \in \mathfrak{h}$. By classical highest weight theory, we know
that $\omega = \sum_{i=1}^l n_i \omegabar_i$ with non-negative integers $n_i$.
From the discussion above, restriction to the Cartan subalgebra of
$\mathfrak{so}(n-1)$ imposes that $n_2 = 0$ or $2$, $n_i = 0$ for all $i > 2$,
while $n_1$ is arbitrary.

Condition~\eqref{BoostC} gives the decay order $k$ of the invariant. Since
$\epsilon_1(a_1) = -1$ and $\epsilon_2(a_1) = 0$ we have
\[\begin{split}
            (n-1-k) \Fhi(p_0)(v_0)
             & = a_1 \cdot \left(\sum_\mu \Fhi^\mu(p_0)(v_0) v_\mu\right)     \\
             & = \omega(a_1) \Fhi(p_0)(v_0)                                   \\
             & = (n_1 \omegabar_1(a_1) + n_2 \omegabar_2(a_1)) \Fhi(p_0)(v_0) \\
             & = -(n_1+n_2) \Fhi(p_0)(v_0),
      \end{split} \]
and hence $k = n-1 + n_1 + n_2$.

The representation with highest weight $\omegabar_2$ is $\Lambda^2 \bR^{n, 1}$.
As a consequence the representation with highest weight $\omega =
      n_1\omegabar_1 + n_2 \omegabar_2$ is obtained as the submodule of highest
weight $n_1 \omegabar_1 + n_2 \omegabar_2$ of
\[
      \underbrace{\bR^{n,1} \otimes \cdots \otimes \bR^{n, 1}}_{n_1\text{ times}}
      \otimes
      \underbrace{
            \Lambda^2\bR^{n,1} \otimes \cdots \otimes \Lambda^2\bR^{n,1}
      }_{n_2\text{ times}} ,
\]
namely the so-called Cartan product, see~\cite[Section~5.5.3]{GoodmanWallach}.

We now briefly discuss the modifications of the argument for the cases $n = 4,
      5$.
\begin{itemize}
      \item
            If $n=4$, the standard representation of $\mathfrak{so}(3)$ has highest weight $2 \omegabar_2$ and hence $\Sym^2_0(T^*_{p_0} \bS^{n-1})$ has highest weight $4\omegabar_2$ for $\mathfrak{so}(n-1)$. When returning to $\mathfrak{so}(4, 1)$ this gives again the representation $\Lambda^2\bR^{4, 1}$ so the previous discussion applies.
      \item
            If $n=5$, the highest weight $\omega$ corresponding to the standard representation of $SO(4)$ has $\omega(H_{l-1}) = \omega(H_l) = 1$. This is the $(1, 1)$-representation according to the previous. This lifts to the representation $\Lambda^2 \bR^{5,1}$ so the previous discussion still applies.
\end{itemize}

Since the longest element in the Weyl group acts on these highest weights by
changing them to their opposite we conclude that the representations $V$ that
we found are isomorphic to their duals,
see~\cite[Section~3.2.3]{GoodmanWallach}.

First, we study the representation with highest weight $n_1 \omegabar_1$. This
corresponds to the space $V = \Sym^{n_1}_0 (\bR^{n, 1})$ of symmetric
trace-free contravariant tensors of rank $n_1$. The highest weight vector of
this representation is $e_{+1} \otimes \cdots \otimes e_{+1}$ so the operator
$\Fhi(p_0)$ reads
\[
      \Fhi(p_0)=
      \underbrace{e_{+1} \otimes \cdots \otimes e_{+1}}_{n_1 \text{ times}} \sigma(p_0).
\]
The vector $e_{+1}$ extends by rotations to $\bS^{n-1}$ as the future pointing
null vector $e_+ = \partial_0 + x^1 \partial_1 + \cdots + x^n \partial_n $. And
hence, by invariance under rotations, the operator $\Fhi$ reads
\[
      \Fhi =
      \underbrace{e_+ \otimes \cdots \otimes e_+}_{n_1 \text{ times}} \sigma.
\]
As in the proof of Proposition~\ref{thmScalarInvariants3d}, the space
$\Sym^{n_1}_0 (\bR^{n, 1})$ can be identified with the dual of $\cH_{n_1}$ of
homogeneous polynomials of degree $n_1$ in $X^0,\ldots,X^n$. Hence, the formula
for these invariants takes a form similar to the one we obtained in the
previous section,
\[
      \begin{array}{rccl}
            \Fhi(m): & \cH_{n_1} & \to     & \bC                                                   \\
                     & P         & \mapsto & \displaystyle{\int_{\bS^{n-1}} P(1, x^1, \ldots, x^n)
                  \tr^\sigma m \, d\mu^\sigma}.
      \end{array}
\]

Second, we look at the representation with highest weight $n_1 \omegabar_1 +
      2\omegabar_2$. This can be realized as the subspace of $\left(\bR^{n,
            1}\right)^{\otimes (n_1 + 4)}$ corresponding to the Young tableau
\begin{center}
      $Y_{n_1} = $ \gyoung(125:\cdots_3<n_1+4>,34).
\end{center}
Namely, let $\mathbf{r}(Y_{n_1})$ denote the row symmetrizer
\[
      \mathbf{r}(Y_{n_1}) =
      (1 + (34)) \sum_{\sigma \in S(\{1, 2, 5, \ldots, n_1+4\})} \sigma,
\]
and let $\mathbf{c}$ be the column skew symmetrizer
\[
      \mathbf{c}(Y_{n_1}) = (1 - (13))(1 - (24)).
\]
These elements are elements of the group algebra $\bC[S_{n_1+4}]$ and act on
$\left(\bR^{n, 1}\right)^{\otimes (n_1 + 4)}$ by permuting the positions of the
elements in a simple tensor.

The Young symmetrizer $\mathbf{s}(Y_{n_1}) =
      \mathbf{c}(Y_{n_1})\mathbf{r}(Y_{n_1})$ is, up to a normalization constant, a
projection from the subspace of harmonic tensors in $(\bR^{n,
            1})^{\otimes(n_1+4)}$ onto an irreducible representation $V$ of highest weight
$n_1 \omegabar_1 + 2\omegabar_2$. Here, as before, harmonic means that the
trace with respect to any pair of indices vanishes.

It can be checked that
\[
      (1 + (123) + (132)) \mathbf{s}(Y_{n_1}) = 0,
\]
and that
\[
      (1 + (125) + (152)) \mathbf{s}(Y_{n_1}) = 0.
\]
These two identities are the (contravariant) analogs of the first and the
second Bianchi identities.

The map $\Fhi(p_0)$ sends the highest weight vector $e^{+2} \otimes e^{+2}$ to
a highest weight vector $v$ in $V$. Up to the action of an element of
$\mathfrak{so}(n,1)$, it is equal to $(e_{+1})^{\otimes n_1} \otimes (e_{+1}
      \wedge e_{+2}) \otimes (e_{+1} \wedge e_{+2})$, which is a highest weight
vector of $V$. By the fact that $\Fhi(p_0)$ intertwines the actions of
$\mathfrak{so}(n-1)$ on $\Sym^2_0(T_{p_0}^* \bS^{n-1})$ and on $V$, it has the
form
\[
      \Fhi(p_0) ( T ) =
      \big\< T , (e_{+1})^{\otimes n_1} \otimes (e_{+1} \wedge \cdot)
      \otimes (e_{+1} \wedge \cdot) \big\>.
\]
Note that since $v$ is harmonic, this map extends by zero to symmetric
2-tensors proportional to $\sigma(p_0)$. The map $\Fhi(p_0)$ is therefore
defined by the above formula on $\Sym^2(T_{p_0}^* \bS^{n-1})$. Due to the
invariance under the action of $SO(n)$, the density $\Fhi$ has the expression
\[
      \begin{array}{rccl}
            \Fhi: & \Sym^2(T^*\bS^{n-1}) & \to     & V \\
                  & m                    & \mapsto &
            \big\< m ,
            (e_+)^{\otimes n_1} \otimes (e_+ \wedge \cdot )
            \otimes (e_+ \wedge \cdot) \big\> .
      \end{array}
\]
As in the $3$-dimensional case (see Proposition~\ref{thmWeylInvariants3d}), we
can view $\Fhi$ as acting on the space $\cW_{n_1}$ of polynomial Weyl tensors
of degree $n_1$ in the variables $X^0,\ldots,X^n$, and the corresponding
intertwining maps can be written, up to a constant factor, as
\[
      \begin{array}{rccl}
            \Fhi(m): & \cW_{n_1} & \to     & \bR \\
                     & W         & \mapsto &
            \displaystyle{\int_{\bS^{n-1}} \left\<m, W(e_+, \cdot, e_+, \cdot)\right\>
                  \, d\mu^\sigma}.
      \end{array}
\]
Here, the polynomial part of the integrand is evaluated at
$(X^0,X^1,\ldots,X^n) = (1,x^1,\ldots,x^n)$.

We now collect the results of this section in the following proposition.

\begin{proposition}\label{thmClassificationso(n,1)}
      Let $\Fhi : S^2(\bS^{n-1}) \to V$ be an $O_{\uparrow}(n,1)$-intertwining map with $n \geq 4$, where $V$ is an irreducible, finite-dimensional representation of $O_{\uparrow}(n,1)$. Then one of the following cases holds.
      \begin{itemize}
            \item
                  Either $V$, as a $\mathfrak{so}(n,1)$-representation, has highest weight $n_1 \omegabar_1$ with $n_1 \geq 0$. Then $V$ can be identified with the dual of $\cH_{n_1}$. The map $\Fhi$ is a multiple of the map
                  \[ \begin{array}{rccl}
                              \Fhi_h (m) : & \cH_{n_1} & \to     & \bR \\
                                           & P         & \mapsto &
                              \displaystyle{\int_{\bS^{n-1}} P(1, x^1,\ldots,x^n) \tr^\sigma(m)
                                    \, d\mu^\sigma},
                        \end{array} \]
                  for $m \in S^2(\bS^{n-1})$. \item Or $V$, as a $\mathfrak{so}(n,1)$-representation, has highest weight
                  $n_1\omegabar_1 + 2\omegabar_2 $ with $n_1 \geq 0$. Then $V$ can be identified
                  with the dual of $\cW_{n_1}$. The map $\Fhi$ is a multiple of the map
                  \[ \begin{array}{rccl}
                              \Fhi_w (m) : & \cW_{n_1} & \to     & \bR \\
                                           & W         & \mapsto &
                              \displaystyle{\int_{\bS^{n-1}} \left\<m, W(e_+, \cdot, e_+, \cdot) \right\>
                                    \, d\mu^\sigma}
                        \end{array} \]
      \end{itemize}
      Further, these maps are $O_{\uparrow}(n,1)$-intertwining, hence are linear masses at infinity, in the first case for asymptotically hyperbolic metrics of order $k = n - 1 + n_1$, in the second case, for metrics of order $k = n + 1 + n_1$.
\end{proposition}

Note that $\Fhi_h$ depends only on the trace of the mass aspect tensor $m$,
while $\Fhi_w$ depends only on the trace-free part of $m$.

\begin{proof}
      The last remaining point to prove in both cases is the $O_{\uparrow}(n,1)$-intertwining property of the maps $\Fhi_h$ and $\Fhi_w$. To do this, we need to carefully keep track of the points where integration takes place. We start by a calculation for the function $u$ defined in~\eqref{def_u[A]}. The function $t \definedas X^0 \circ p^{-1}$, which is the time component of a point $x$ in the hyperbolic space, is given by
      \[
            t(x) = \frac{1 + |x|^2}{1 - |x|^2} = \frac{1}{\rho(x)} - 1.
      \]
      As a consequence we have
      \[\begin{split}
                  u[A](\xhat)
                   & = \lim_{\lambda \to 1}
                  \frac{\rho(\Abar^{-1}(\lambda \xhat))}{\rho(\lambda \xhat)}      \\
                   & = \lim_{\lambda \to 1}
                  \frac{t(\lambda \xhat) + 1}{t(\Abar^{-1}(\lambda \xhat)) + 1}    \\
                   & = \lim_{\lambda \to 1}
                  \frac{t(\lambda \xhat)}{t(\Abar^{-1}(\lambda \xhat))}            \\
                   & = \lim_{\lambda \to 1} \left(
                  X^0\left(A^{-1} \cdot
                  \frac{p^{-1}(\lambda \xhat)}{X^0(p^{-1}(\lambda \xhat))}\right)
                  \right)^{-1}                                                     \\
                   & = \left( X^0\left(A^{-1} \cdot (1, \xhat)\right) \right)^{-1}
            \end{split} \]
      for any $\xhat \in \bS^{n-1}$. Recall here that $\Abar = p A p^{-1}$ gives the
      action of $A \in O_\uparrow(n, 1)$ on the ball model of hyperbolic space. This
      can be rewritten as
      \[
            A^{-1} (1, \xhat) = \frac{1}{u[A](\xhat)} (1, \Abar^{-1}(\xhat)).
      \]

      From the proof of Corollary~\ref{remGroupAction} we have
      \[
            \Abar_* \sigma = u[A]^2 \sigma,
      \]
      and thus
      \[
            d\mu^{\Abar_* \sigma} = u[A]^{n-1} d\mu^\sigma.
      \]
      From Lemma~\ref{lmGroupActionMassAspect} we know that $A \cdot m = u[A]^{k-2}
            \Abar_* m$. So for $\Fhi_h$ we get
      \[ \begin{split}
                  \Fhi_h(A \cdot m) (P \circ A^{-1})
                   & =
                  \int_{\bS^{n-1}} P(A^{-1}(1, \xhat)) u[A]^{k-2} \tr^\sigma (\Abar_* m)
                  \, d\mu^{\sigma}                                          \\
                   & =
                  \int_{\bS^{n-1}} u[A]^{k-2-n_1} P(1, \Abar^{-1}\xhat) \tr^\sigma (\Abar_* m)
                  \, d\mu^{\sigma}                                          \\
                   & =
                  \int_{\bS^{n-1}} u[A]^{k-n_1-n + 1} P(1, \Abar^{-1}\xhat)
                  \tr^{\Abar_* \sigma} (\Abar_* m) \, d\mu^{\Abar_* \sigma} \\
                   & =
                  \int_{\bS^{n-1}}
                  P(1, \Abar^{-1}\xhat) \tr^{\Abar_* \sigma} (\Abar_* m) \, d\mu^{\Abar_* \sigma}
            \end{split} \]
      since $P$ is a homogeneous polynomial of degree $n_1$, and since $k = n-1+n_1$.
      By a change of variables we conclude that
      \[
            \Fhi_h(A \cdot m) (P \circ A^{-1}) = \Fhi_h(m) (P)
      \]
      which is the desired $O_{\uparrow}(n,1)$-intertwining property of the map
      $\Fhi_h$.

      To prove the intertwining property for $\Fhi_w$ we compute
      \[ \begin{split}
                   & \Fhi_w(A \cdot m)(A_* W) \\
                   & \qquad=
                  \int_{\bS^{n-1}}
                  \left\<u[A]^{k-2} \Abar_* m,
                  (A_* W)(e_+, \cdot, e_+, \cdot)(1, \xhat) \right\>_\sigma
                  \, d\mu^\sigma(\xhat)       \\
                   & \qquad=
                  \int_{\bS^{n-1}} u[A]^{k-2}
                  \left\<\Abar_* m, W(A^*e_+, A^*\cdot,
                  A^*e_+, A^*\cdot)(A^{-1}(1, \xhat)) \right\>_\sigma
                  \, d\mu^\sigma(\xhat).
            \end{split} \]
      We extend $e_+$ to all of $\bR^{n, 1}$ by defining it to be the position vector
      field, that is
      \[
            e_+ \definedas X^\mu \partial_\mu.
      \]
      Then
      \[\begin{split}
                   & \Fhi_w(A \cdot m)(A_* W) \\
                   & \qquad=
                  \int_{\bS^{n-1}} u[A]^{k-2}
                  \left\<\Abar_* m,
                  W(e_+, A^*\cdot, e_+, A^*\cdot)(A^{-1}(1, \xhat)) \right\>_\sigma
                  \, d\mu^\sigma(\xhat).
            \end{split}\]
      The vectors $e_+$ in this expression are evaluated at the point
      $A^{-1}(1,\xhat)$, and since $A^{-1} (1, \xhat) = (u[A](\xhat))^{-1} (1,
            \Abar^{-1}(\xhat))$ we have
      \[
            e_+(A^{-1}(1, \xhat)) = (u[A](\xhat))^{-1} e_+(1, \Abar^{-1}\xhat).
      \]
      The vectors $A^* \cdot$ are also evaluated at the point $A^{-1}(1, \xhat)$,
      where $\cdot$ is the place for a vector belonging to $T_{(1,
                        \Abar^{-1}\xhat)}\bS^{n-1}$. Thus we get
      \[
            (A^* \cdot) (A^{-1}(1, \xhat))
            = u[A]^{-1} (\Abar^* \cdot)(1, \Abar^{-1} \xhat).
      \]
      Since $W$ is homogeneous of degree $n_1$, we find that
      \[ \begin{split}
                   & W(e_+, A^*\cdot, e_+, A^*\cdot)(A^{-1}(1, \xhat)) \\
                   & \qquad=
                  u[A]^{-4-n_1} W(e_+, \Abar^*\cdot, e_+, \Abar^*\cdot)((1, \Abar^{-1}\xhat)).
            \end{split} \]
      We can now continue the computation,
      \[ \begin{split}
                   & \Fhi_w(A \cdot m)(A_* W)                              \\
                   & \qquad=
                  \int_{\bS^{n-1}} u[A]^{k-6-n_1}
                  \left\<\Abar_* m,
                  W(e_+, \Abar^*\cdot, e_+, \Abar^*\cdot)((1, \Abar^{-1}\xhat))
                  \right\>_\sigma \, d\mu^\sigma(\xhat)                    \\
                   & \qquad=
                  \int_{\bS^{n-1}} u[A]^{k-n_1-n-1}
                  \left\<
                  \Abar_* m, W(e_+, \Abar^*\cdot, e_+, \Abar^*\cdot)((1, \Abar^{-1}\xhat))
                  \right\>_{\Abar_* \sigma} \, d\mu^{\Abar_*\sigma}(\xhat) \\
            \end{split} \]
      Note that the metric $\sigma$ is used twice for the scalar product, which
      explains the change in the exponent of $u[A]$ when passing to the last line.
      Since $k = n + 1 + n_1$, a change of variables gives us
      \[
            \Fhi_w(A \cdot m)(A_* W) = \Fhi_w(m)(W),
      \]
      which is the $O_{\uparrow}(n,1)$-intertwining property of the map $\Fhi_w$.
\end{proof}

An immediate corollary of Proposition~\ref{thmClassificationso(n,1)} is that,
given any order $k$ of decay for asymptotically hyperbolic metrics, one can
tell the number of irreducible representations of $O_{\uparrow}(n, 1)$ that are
the ranges of linear masses:
\begin{itemize}
      \item In the case $k < n-1$, there is no such representation, hence no linear mass.
      \item In the case $k = n-1,n$, there is a unique such representation that can be
            identified with $\cH_{n_1}^*$, the dual of $\cH_{n_1}$, where $n_1 = k - n +
                  1$. This corresponds to the linear mass of example~\ref{example_trivialmass}
            for $k = n-1$, respectively with example~\ref{example_AHmassWang} (the standard
            mass) for $k = n$.
      \item In the case $k \geq n+1$, there are two real irreducible representations
            of $O_{\uparrow}(n, 1)$ that are the ranges of linear masses (one of which
            splits into two irreducible complex representations of $SO_\uparrow(n, 1)$ when
            $n=3$, see Proposition~\ref{thmWeylInvariants3d}). Those can be identified
            respectively with the dual of $\cH_{n_1}$ with $n_1 = k - n + 1$, and the dual
            of $\cW_{n_1}$ with $n_1 = k - n - 1$.
\end{itemize}

%% file: weyl.tex
An important fact concerning the (complex) representations that appear in
Proposition \ref{thmClassificationso(n,1)} is that they are of real type,
namely they are the complexification of real irreducible representations of
$O_\uparrow(n, 1)$. See for example \cite[Section~26.3]{FultonHarris} for an
introduction to real and quaternionic representations. The situation is
slightly more complicated for the Weyl masses in dimension $3$ since one has to
consider the sum $V_{n_1\omegabar_1 + (n_1 +4)\omegabar_2} \oplus V_{(n_1
      +4)\omegabar_1 + n_1\omegabar_2}$ to get a real representation.

The real representations whose complexifications give the representations
$\cH_{n_1}$ and $\cW_{n_1}$ are easy to find. Namely, they are the set of real
harmonic polynomials and the set of real Weyl tensors (see
Subsection~\ref{secWeyl}). An interesting consequence of this is that the real
representations carry an invariant quadratic form and we will compute here its
signature, see Subsections \ref{secHarmForm} and \ref{secWeylForm}.

\subsection{The set of harmonic polynomials}
\label{secHarm}

From the Weyl dimension formula (see for example
\cite[Section~7.1.2]{GoodmanWallach}) it follows that the dimension of $\cH_p$
is
\[
  \dim \cH_p = \binom{p + n-2}{p} \frac{2p + n-1}{n-1}.
\]
The calculations for this are straightforward so we omit them. It suffices to
plug $\lambda_1 = p$ and $\lambda_i = 0$ for all $i > 1$ in the Weyl dimension
formulas for $SO(2l)$ and $SO(2l+1)$.

\subsubsection{Invariant quadratic form}
\label{secHarmForm}

The aim of this section is to compute the signature $(n_+(p), n_-(p))$ of the
$O_\uparrow(n, 1)$-invariant quadratic form $q$ on $\cH_p$. There are general
methods to compute this signature, but this representation is simple enough so
that a shorter derivation can be found. See Subsection \ref{secWeylForm} for a
more robust method.

Note that the signature of the invariant quadratic form is not unique, but only
defined up to swapping of $n_+(p)$ and $n_-(p)$ since a given invariant
quadratic form can be multiplied by any non-zero number.

For simplicity we make the choice that $q(X^0) = -1$ and $q(X^i) = 1$ for all
$i = 1, \dots, n$, and we extend $q$ to $\bR[X^0, X^1, \dots, X^n]$ by
requiring that distinct monomials are orthogonal and that
\[
  q\left((X^0)^{k_0} (X^1)^{k_1} \cdots (X^n)^{k_n}\right)
  = q(X^0)^{k_0} q(X^1)^{k_1} \cdots q(X^n)^{k_n}
  = (-1)^{k_0}.
\]
This form is the natural extension of the invariant quadratic form $q$ on
$(\bR^{n, 1})^*$ to $\mathrm{Sym}_* \bR^{n, 1}$ and, as such, it is
$O_\uparrow(n, 1)$-invariant.

The signature $(n_{+, 0}(p), n_{-, 0}(p))$ of $q$ restricted to $\bR_p[X^0,
    X^1, \dots, X^n]$, the space of homogeneous polynomials of degree $p$, can be
computed easily from the fact that monomials are orthogonal,
\[
  \left\lbrace
  \begin{aligned}
    n_{+, 0}(p) & = \# \{\text{monomials with even power of $X^0$}\}, \\
    n_{-, 0}(p) & = \# \{\text{monomials with odd power of $X^0$}\}.
  \end{aligned}
  \right.
\]
Thus,
\[
  \left\lbrace
  \begin{aligned}
    n_{+, 0}(p) & = \sum_{\substack{k \in \{0, 1, \dots, p\} \\ k~\mathrm{even}}}
    \binom{p-k+n-1}{n-1},                                    \\
    n_{-, 0}(p) & = \sum_{\substack{k \in \{0, 1, \dots, p\} \\ k~\mathrm{odd}}}
    \binom{p-k+n-1}{n-1}.
  \end{aligned}
  \right.
\]
This can be condensed into
\[
  \left\lbrace
  \begin{aligned}
    n_{+, 0}(p) + n_{-, 0}(p) & = \sum_{k=0}^p \binom{p-k+n-1}{n-1},        \\
    n_{+, 0}(p) - n_{-, 0}(p) & = \sum_{k=0}^p (-1)^k \binom{p-k+n-1}{n-1}.
  \end{aligned}
  \right.
\]
The space $\cH_p$ is the kernel of the map
\[
  \Box: \bR_p[X^0, X^1, \dots, X^n] \to \bR_{p-2}[X^0, X^1, \dots, X^n]
\]
and a complementary orthogonal subspace to $\cH_p$ is given by the space
$\cK_p$ of polynomials divisible by $-(X^0)^2 + (X^1)^2 + \cdots + (X^n)^2$. It
cannot be expected that $\Box$ is an isometry from $\cK_p$ to $\bR_{p-2}[X^0,
    X^1, \dots, X^n]$. We have however the following result.

\begin{lemma}\label{lmSignaturePoly}
  The signature of $q$ restricted to $\cK_p$ is the same as the signature
  of $q$ on $\bR_{p-2}[X^0, X^1, \dots, X^n]$.
\end{lemma}

From this lemma it follows that the signature of the invariant quadratic form
$q$ on $\cH_k$ is given by
\[
  n_+(p) = n_{+, 0}(p) - n_{+, 0}(p-2),\quad
  n_-(p) = n_{-, 0}(p) - n_{-, 0}(p-2),
\]
and thus
\[
  n_+(p) = \binom{p+n-1}{n-1}, \quad n_-(p) = \binom{p+n-2}{n-1}.
\]

\begin{proof}[Proof of Lemma \ref{lmSignaturePoly}]
  The proof is by a Brauer-style diagrammatic argument. The representation
  $\bR_p[X^0, X^1, \dots, X^n]$ decomposes as a sum of irreducible
  subrepresentations,
  \[
    \bR_p[X^0, X^1, \dots, X^n] =
    \bigoplus_{k=0}^{\left\lfloor p/2\right\rfloor}
    \left(-(X^0)^2 + (X^1)^2 + \cdots + (X^n)^2\right)^k \cH_{p-2k}.
  \]
  It is not complicated to check that the map
  \[
    \begin{array}{rccc}
      m_k : & \cH_r & \to     & \bR_{r+2k}[X^0, X^1, \cdots, X^n]                      \\
            & P     & \mapsto & \left(-(X^0)^2 + (X^1)^2 + \cdots + (X^n)^2\right)^k P
    \end{array}
  \]
  preserves the quadratic form up to a positive constant, that is there exists a
  constant $\lambda_{k, r}$ such that for any $P \in \cH_r$, we have $q(m_k(P)) =
    \lambda_{k, r} q(P)$.

  We represent an arbitrary polynomial $Q \in \bR_r [X^0, X^1, \cdots, X^n]$ as a
  symmetric tensor with $r$ indices,
  \[
    Q = Q_{\mu_1 \cdots \mu_r} X^{\mu_1} \cdots X^{\mu_r}.
  \]
  Multiplication by $\left(-(X^0)^2 + (X^1)^2 + \cdots + (X^n)^2\right)^k$ is
  then the same as tensoring $k$ times with $\eta = \eta_{\mu\nu} X^\mu X^\nu$
  and symmetrizing. That is
  \[ \begin{split}
       & (m_k(Q))_{\mu_1 \cdots \mu_{r+2k}}
      =
      Q_{(\mu_1 \cdots \mu_r} \eta_{\mu_{r+1}\mu_{r+2}} \cdots \eta_{\mu_{r+2k-1}\mu_{r+2k})} \\
       & \qquad
      = \frac{1}{(r+2k)!} \sum_{\sigma \in S_{r+2k}} Q_{\mu_{\sigma(1)} \cdots \mu_{\sigma(r)}}
      \eta_{\mu_{\sigma(r+1)}\mu_{\sigma(r+2)}} \cdots \eta_{\mu_{\sigma(r+2k-1)}\mu_{\sigma(r+2k)}}.
    \end{split} \]

  The form $q$ can be written as
  \[
    q(Q) = Q_{\mu_1 \cdots \mu_r} Q^{\mu_1 \cdots \mu_r},
  \]
  where indices are raised with respect to $\eta$. We have
  \[\begin{split}
      q(m_k(Q))
       & =
      \frac{1}{((r+2k)!)^2} \sum_{\sigma \in S_{r+2k}} Q_{\mu_{\sigma(1)} \cdots \mu_{\sigma(r)}}
      \eta_{\mu_{\sigma(r+1)}\mu_{\sigma(r+2)}} \cdots \eta_{\mu_{\sigma(r+2k-1)}\mu_{\sigma(r+2k)}}          \\
       & \qquad
      \times \sum_{\tau \in S_{r+2k}} Q^{\mu_{\tau(1)} \cdots \mu_{\tau(r)}}
      \eta^{\mu_{\tau(r+1)}\mu_{\tau(r+2)}} \cdots \eta^{\mu_{\tau(r+2k-1)}\mu_{\tau(r+2k)}}                  \\
       & =
      \frac{1}{(r+2k)!} Q^{\mu_1 \cdots \mu_r} \eta^{\mu_{r+1}\mu_{r+2}} \cdots \eta^{\mu_{r+2k-1}\mu_{r+2k}} \\
       & \qquad
      \times \sum_{\sigma \in S_{r+2k}} Q_{\mu_{\sigma(1)} \cdots \mu_{\sigma(r)}}
      \eta_{\mu_{\sigma(r+1)}\mu_{\sigma(r+2)}} \cdots \eta_{\mu_{\sigma(r+2k-1)}\mu_{\sigma(r+2k)}}.
    \end{split} \]
  That is to say that $q(m_k(Q))$ is, up to a factor $\frac{1}{((r+2k)!)^2}$, the
  sum over all possible ways of pairing indices from two copies of $m_k(Q)$. Such
  a pairing can be represented as a Brauer-like diagram:
  \begin{center}
    \begin{tikzpicture}[scale=1]
      \node[bnode] at (0,0) (p0) {};
      \node[bnode] at (1,0) (p1) {};
      \node at (2,0) (p2) {$\cdots$};
      \node[bnode] at (3,0) (p3) {};
      \node[bnode] at (4,0) (p4) {};
      \draw (p0) -- +(0, -0.2) -- +(1, -0.2) -- (p1);
      \draw (p3) -- +(0, -0.2) -- +(1, -0.2) -- (p4);
      \node[bnode] at (5, 0) (p5) {};
      \node[bnode] at (6, 0) (p6) {};
      \node at (7,0) (p7) {$\cdots$};
      \node[bnode] at (8, 0) (p8) {};
      \draw[decorate,decoration={brace, mirror}] (0, -0.4) -- (1, -0.4);
      \node[below] at (0.5, -0.4) {$\eta$};
      \draw[decorate,decoration={brace, mirror}] (3, -0.4) -- (4, -0.4);
      \node[below] at (3.5, -0.4) {$\eta$};
      \draw[decorate,decoration={brace, mirror}] (5, -0.2) -- (8, -0.2);
      \node[below] at (6.5, -0.2) {$Q$};

      \node[bnode] at (0,2) (q0) {};
      \node[bnode] at (1,2) (q1) {};
      \node at (2,2) (q2) {$\cdots$};
      \node[bnode] at (3,2) (q3) {};
      \node[bnode] at (4,2) (q4) {};
      \draw (q0) -- +(0, 0.2) -- +(1, 0.2) -- (q1);
      \draw (q3) -- +(0, 0.2) -- +(1, 0.2) -- (q4);
      \node[bnode] at (5, 2) (q5) {};
      \node[bnode] at (6, 2) (q6) {};
      \node at (7,2) (q7) {$\cdots$};
      \node[bnode] at (8, 2) (q8) {};

      \draw (q0) .. controls +(0, -0.5) and +(0, 0.5) .. (p3);
      \draw (q1) .. controls +(0, -0.5) and +(0, 0.5) .. (p4);
      \draw (q3) .. controls +(0, -0.5) and +(0, 0.5) .. (p0);
      \draw (q4) .. controls +(0, -0.5) and +(0, 0.5) .. (p6);
      \draw (q5) .. controls +(0, -0.5) and +(0, 0.5) .. (p1);
      \draw (q6) .. controls +(0, -0.5) and +(0, 0.5) .. (p8);
      \draw (q8) .. controls +(0, -0.5) and +(0, 0.5) .. (p5);
    \end{tikzpicture}
  \end{center}
  where there are $k$
  \begin{tikzpicture}[scale=0.4]
    \node[bnode2] at (0,0) (p0) {};
    \node[bnode2] at (1,0) (p1) {};
    \draw (p0) -- +(0, -0.2) -- +(1, -0.2) -- (p1);
  \end{tikzpicture}
  and as many
  \begin{tikzpicture}[scale=0.4]
    \node[bnode2] at (0,0) (p0) {};
    \node[bnode2] at (1,0) (p1) {};
    \draw (p0) -- +(0, 0.2) -- +(1, 0.2) -- (p1);
  \end{tikzpicture}
  corresponding to the $\eta$'s we added to $Q$
  and $2r$ isolated points ($r$ above and $r$ below) corresponding to
  the $r$ indices of $Q$. Each diagram contains (at most) three types of
  curves drawn by contractions (curvy lines) and $\eta$'s:
  \begin{enumerate}
    \item Closed curves (loops) corresponding to contractions of $\eta$'s with themselves
          so they give rise to a factor $(n+1)$ to the diagram.
    \item Curves that connect the same $Q$ (either the upper one or the lower one)
          corresponding to a contraction of $Q$ with itself. Since $Q$ is traceless,
          having such a curve in a diagram makes its contribution vanish.
    \item Curves that connect the upper $Q$ and the lower $Q$.
  \end{enumerate}
  In the end, we get the following formula for $q(m_k(Q))$:
  \[
    q(m_k(Q))
    = \frac{q(Q)}{(r+2k)!} \sum_{D \in \cD_0} (n+1)^{\# \mathrm{loops}(D)},
  \]
  where $\cD_0$ denotes the set of diagrams for which there is no curve
  contracting the same $Q$. It follows that $q(m_k(Q)) = \lambda_{k, r} q(Q)$
  with
  \[
    \lambda_{k, r}
    = \frac{1}{(r+2k)!} \sum_{D \in \cD_0} (n+1)^{\# \mathrm{loops}(D)} > 0.
  \]
\end{proof}

\subsection{The set of polynomial Weyl tensors}
\label{secWeyl}

Let $\cE_{p+2}$ be the set of solutions $h_{\mu\nu}$ of the Einstein equations
linearized around the Minkowski metric $\eta$ that are homogeneous polynomials
of degree $p+2$. That is $h_{\mu\nu}$ satisfies
\begin{equation}\label{eqLinEinstein}
  0 =
  - \Box h_{\mu\nu}
  + \partial_\mu \partial^\alpha h_{\alpha\nu}
  + \partial_\nu \partial^\alpha h_{\mu\alpha}
  - \partial_\mu \partial_\nu h^\alpha_\alpha
  - \partial^\alpha \partial^\beta h_{\alpha\beta} \eta_{\mu\nu}
  + \eta_{\mu\nu} \Box h^\alpha_\alpha.
\end{equation}
The following proposition gives an alternative characterization of
the set of polynomial Weyl tensors $\cW_{p}$.

\begin{proposition}\label{propWeyl}
  The set $\cW_{p}$ of polynomial Weyl tensors on $\bR^{n, 1}$ is the
  (dual of) the irreducible representation of $O_\uparrow(n, 1)$ given
  by the following Young tableau
  \begin{center}
    $Y_{p} = $ \gyoung(125:\cdots_3<p+4>,34)
  \end{center}
  This is the image of $\cE_{p+2}$ under the linearized Riemann tensor,
  \[
    \riem: h_{\mu\nu} \mapsto R_{\mu\nu\alpha\beta} \definedas
    -\frac{1}{2} \left( \partial_\mu \partial_\alpha h_{\nu\beta}
    + \partial_\nu \partial_\beta h_{\mu\alpha}
    - \partial_\mu \partial_\beta h_{\nu\alpha}
    - \partial_\nu \partial_\alpha h_{\mu\beta}\right).
  \]
\end{proposition}

We start the proof by showing that given $W \in \cW_{p}$ there exists a
solution $h$ to the linearized Einstein equations such that $W = \riem(h)$.

\begin{lemma}\label{lmSurjectivity}
  Given $W \in \cW_p$ there exists $h \in \cE_{p+2}$ such that $W = \riem(h)$.
\end{lemma}

The proof of this is based on a straightforward refinement of the well-known
Poincar\'e Lemma for cohomology, see for example \cite[Chapter~1]{BottTu}.

\begin{lemma}\label{lmPoincare}
  If $\omega$ is a closed $k$-form, $k \geq 1$, on $\bR^{n, 1}$ with
  coefficients that are homogeneous polynomials of degree $l$, there
  exists a $(k-1)$-form $\xi$ whose coefficients are homogeneous
  polynomials of degree $l+1$ such that $\omega = d \xi$.
\end{lemma}

\begin{proof}
  We introduce the operators $\cI_k: \Lambda^k\bR^{n, 1} \to \Lambda^{k-1}\bR^{n,
      1}$, for all $k>0$, defined as follows. For $\omega = f dX^{\mu_1} \wedge \dots
    \wedge dX^{\mu_k} \in \Lambda^k\bR^{n, 1}$ we set
  \[
    \cI_k(\omega)(p) \definedas
    \sum_{j=1}^k (-1)^{j-1} X^{\mu_j} \int_0^1 t^{k-1} f(tp) dt
    dX^{\mu_1} \wedge \dots \wedge \widehat{dX^{\mu_j}} \wedge \dots \wedge dX^{\mu_k}.
  \]
  and extend $\cI_k$ to all of $\Lambda^k\bR^{n, 1}$ by linearity. We claim that
  for any $\omega \in \Lambda^k\bR^{n, 1}$ we have
  \[
    \omega = d (\cI_k \omega) + \cI_{k+1} (d\omega).
  \]
  Indeed,
  \[\begin{split}
       & d\left(\cI_k \omega\right)(p)                                                                                                                                        \\
       & \qquad =
      \sum_{j=1}^k (-1)^{j-1} \int_0^1 t^{k-1} f(tp) dt~dX^{\mu_j} \wedge dX^{\mu_1} \wedge \dots \wedge \widehat{dX^{\mu_j}} \wedge \dots \wedge dX^{\mu_k}                  \\
       & \qquad\qquad
      + \sum_{j=1}^k (-1)^{j-1} X^{\mu_j} \int_0^1 t^k \partial_\mu f(tp) dt~dX^\mu \wedge dX^{\mu_1} \wedge \dots \wedge \widehat{dX^{\mu_j}} \wedge \dots \wedge dX^{\mu_k} \\
       & \qquad =
      k \int_0^1 t^{k-1} f(tp) dt~dX^{\mu_j} \wedge dX^{\mu_1} \wedge \dots \wedge \widehat{dX^{\mu_j}} \wedge \dots \wedge dX^{\mu_k}                                        \\
       & \qquad\qquad
      + \sum_{j=1}^k (-1)^{j-1} X^{\mu_j} \int_0^1 t^k \partial_\mu f(tp) dt~dX^\mu \wedge dX^{\mu_1} \wedge \dots \wedge \widehat{dX^{\mu_j}} \wedge \dots \wedge dX^{\mu_k},
    \end{split}\]
  and
  \[\begin{split}
       & \cI_{k+1} \left(d\omega\right)(p)                                               \\
       & \qquad =
      \int_0^1 t^k X^\mu \partial_\mu f(tp) dt~dX^{\mu_1} \wedge \dots \wedge dX^{\mu_k} \\
       & \qquad\qquad
      - \sum_{j=1}^k (-1)^{j-1} X^{\mu_j} \int_0^1 t^k \partial_\mu f(tp) dt~dX^\mu \wedge dX^{\mu_1} \wedge \dots \wedge \widehat{dX^{\mu_j}} \wedge \dots \wedge dX^{\mu_k},
    \end{split}\]
  so
  \[\begin{split}
       & d\left(\cI_k \omega\right)(p) + \cI_{k+1} \left(d\omega\right)(p)                                              \\
       & \qquad =
      \int_0^1 \left(k t^{k-1} f(tp) + t^k X^\mu \partial_\mu f(tp)\right) dt~dX^{\mu_1} \wedge \dots \wedge dX^{\mu_k} \\
       & \qquad =
      \int_0^1 \left(k t^{k-1} f(tp) + t^k \frac{d}{dt}(f(tp))\right) dt~dX^{\mu_1} \wedge \dots \wedge dX^{\mu_k}      \\
       & \qquad =
      \left[t^k f(tp) \right]_0^1 dX^{\mu_1} \wedge \dots \wedge dX^{\mu_k}                                             \\
       & \qquad =
      f(p) dX^{\mu_1} \wedge \dots \wedge dX^{\mu_k}.
    \end{split}\]
  In particular, if $\omega$ is closed, that is $d\omega = 0$, we have $\omega =
    d(\cI \omega)$. It suffices to note that if $\omega$ has coefficients that are
  homogeneous polynomials of degree $l$, then $\cI \omega$ has coefficients that
  are homogeneous polynomials of degree $l+1$.
\end{proof}

We now prove Lemma~\ref{lmSurjectivity} in a sequence of claims.

\begin{claim}\label{cl1}
  There exists a tensor $f_{\mu\alpha\beta}$ which is antisymmetric with
  respect to the last two indices such that
  \[
    W_{\mu\nu\alpha\beta}
    = \partial_\mu f_{\nu\alpha\beta} - \partial_\nu f_{\mu\alpha\beta} .
  \]
  Further, the components of $f$ are homogeneous polynomials of degree $p+1$.
\end{claim}

\begin{proof}
  The second Bianchi identity for $W$ can be rewritten
  \[
    d (W_{\mu\nu\alpha\beta} dX^\mu \wedge dX^\nu) = 0.
  \]
  Since $H^2(\bR^{n, 1}) = 0$, the existence of $f$ follows from Lemma
  \ref{lmPoincare}.
\end{proof}

\begin{claim}\label{cl2}
  We can assume that $f$ satisfies
  \[
    f_{\nu\alpha\beta} + f_{\alpha\beta\nu} + f_{\beta\nu\alpha} = 0.
  \]
\end{claim}

\begin{proof}
  Note that in the choice of $f$, we are free to add an arbitrary exact
  form $d\theta$. So we can change $f$ to $\ftil = f + d \theta$ where
  $\theta$ is an arbitrary $\Lambda^2\bR^{n, 1}$-valued $0$-form. In
  component notation,
  \[
    f_{\mu\alpha\beta} \rightsquigarrow
    \ftil_{\mu\alpha\beta} = f_{\mu\alpha\beta} + \partial_\mu \theta_{\alpha\beta}.
  \]
  The condition
  \[
    \ftil_{\nu\alpha\beta} + \ftil_{\alpha\beta\nu} + \ftil_{\beta\nu\alpha} = 0
  \]
  reads
  \[
    0 = f_{\nu\alpha\beta} + f_{\alpha\beta\nu} + f_{\beta\nu\alpha}
    + \partial_\nu \theta_{\alpha\beta}
    + \partial_\alpha \theta_{\beta\nu}
    + \partial_\beta \theta_{\nu\alpha}
  \]
  or
  \[
    0 =
    f_{\nu\alpha\beta} dX^\nu \wedge dX^\alpha \wedge dX^\beta
    + d (\theta_{\alpha\beta} dX^\alpha \wedge dX^\beta).
  \]
  To guarantee the existence of a $\theta$ solving this equation, we have to
  prove that $f_{\nu\alpha\beta} dX^\nu \wedge dX^\alpha \wedge dX^\beta$ is
  closed. To do this, we use the first Bianchi identity for $W$,
  \[\begin{split}
      0 & =
      W_{\mu\nu\alpha\beta} + W_{\alpha\mu\nu\beta} + W_{\nu\alpha\mu\beta} \\
        & =
      \partial_\mu f_{\nu\alpha\beta}
      + \partial_\alpha f_{\mu\nu\beta}
      + \partial_\nu f_{\alpha\mu\beta}
      - \partial_\nu f_{\mu\alpha\beta}
      - \partial_\alpha f_{\nu\mu\beta}
      - \partial_\mu f_{\alpha\nu\beta}                                     \\
        & =
      \partial_\mu (f_{\nu\alpha\beta} - f_{\alpha\nu\beta})
      + \partial_\nu (f_{\alpha\mu\beta} - f_{\mu\alpha\beta})
      + \partial_\alpha (f_{\mu\nu\beta} - f_{\nu\mu\beta}),
    \end{split} \]
  or
  \[
    0 = d (f_{\mu\nu\beta} dX^\mu \wedge dX^\nu).
  \]
  Taking the wedge product with $dX^\beta$, we obtain
  \[
    d (f_{\mu\nu\beta} dX^\mu \wedge dX^\nu \wedge dX^\beta)
    = d (f_{\mu\nu\beta} dX^\mu \wedge dX^\nu) \wedge dX^\beta
    = 0
  \]
  which is the desired identity.
\end{proof}

\begin{claim}\label{cl3}
  There exists a tensor $h_{\alpha\beta}$ such that
  \[
    f_{\mu\alpha\beta} = \partial_\alpha h_{\mu\beta} - \partial_\beta h_{\mu\alpha}.
  \]
  Further the components of $h$ are homogeneous polynomials of degree $p+2$.
\end{claim}

\begin{proof}
  As for Claim~\ref{cl1}, it suffices to see that $f_{\mu\alpha\beta} dX^\alpha\wedge dX^\beta$
  is closed. This fact was already obtained in the proof of Claim~\ref{cl2}.
\end{proof}

\begin{claim}\label{cl4}
  The tensor $h_{\alpha\beta}$ can be assumed to be symmetric.
\end{claim}

\begin{proof}
  The solution $h_{\alpha\beta}$ from Claim~\ref{cl3} is not unique, the
  form $h_{\alpha\beta} dX^\beta$ is only defined up to the addition of
  an exact form $dv_\alpha$. So we can change $h_{\alpha\beta}$ to
  $\htil_{\alpha\beta} = h_{\alpha\beta} + \partial_\beta v_\alpha$
  without changing $f$ (and hence $W$). We get a symmetric $\htil$ if
  and only if the equation
  \[
    0 = h_{\alpha\beta} - h_{\beta\alpha}
    + \partial_\beta v_\alpha - \partial_\alpha v_\beta
  \]
  admits a solution $v_\alpha$. This condition is the same as
  \[
    d(v_\alpha dX^\alpha) = h_{\alpha\beta} dX^\alpha \wedge dX^\beta
  \]
  so it suffices to check that $h_{\alpha\beta} dX^\alpha \wedge dX^\beta$ is
  closed. This follows from Claim~\ref{cl2}. Indeed,
  \[ \begin{split}
      0 & =
      f_{\nu\alpha\beta} + f_{\alpha\beta\nu} + f_{\beta\nu\alpha} \\
        & =
      \partial_\alpha h_{\mu\beta}
      + \partial_\beta h_{\alpha\mu}
      + \partial_\mu h_{\beta\alpha}
      - \partial_\beta h_{\mu\alpha}
      - \partial_\alpha h_{\beta\mu}
      - \partial_\mu h_{\alpha\beta}                               \\
        & =
      \partial_\mu (h_{\beta\alpha} - h_{\alpha\beta})
      + \partial_\alpha (h_{\mu\beta} - h_{\beta\mu})
      + \partial_\beta (h_{\alpha\mu} - h_{\mu\alpha}),
    \end{split} \]
  so
  \[
    0 = -i_{\partial_\mu} d(h_{\alpha\beta} dX^\alpha \wedge dX^\beta).
  \]
\end{proof}

Now $W$ can be written in terms of $h$ and we find that $W = -2\riem(h)$. The
condition \eqref{eqLinEinstein} follows at once by the zero trace condition
imposed on $W$ since this condition is nothing but the vanishing of the
linearized Ricci tensor around the Minkowski metric taken with respect to $h$.

We now want to identify a subrepresentation of the homogeneous polynomial
solutions to the linearized Einstein equations with the property that it
projects onto the set of Weyl tensors. A natural gauge condition in the context
of gravitational waves is de Donder's gauge condition, which requires that
\begin{equation}\label{eqDeDonder}
  h^\mu_\mu = 0,\quad \partial^\mu h_{\mu\nu} = 0.
\end{equation}
The gauge freedom that appears for the linearized Einstein equations
is the action of infinitesimal isometries,
\[
  h_{\mu\nu} \to h_{\mu\nu} + \partial_\mu \xi_\nu + \partial_\nu \xi_\mu,
\]
for an arbitrary 1-form $\xi_\mu$ that is a homogeneous polynomial of degree
$p+3$. The main interest of this gauge condition is that Equation
\eqref{eqLinEinstein} reduces to the wave equation
\[
  \Box h_{\mu\nu} = 0.
\]
The fact that de Donder's gauge condition can be satisfied is the content of
the next lemma.

\begin{lemma}\label{lmDeDonder}
  Given an arbitrary homogeneous polynomial symmetric $2$-tensor
  $h_{\mu\nu}$ of degree $p+2$ solving \eqref{eqLinEinstein}, there
  exists a homogeneous polynomial $1$-form $\xi_\mu$ of degree $p+3$
  such that $\htil_{\mu\nu} \definedas
    h_{\mu\nu} + \partial_\mu \xi_\nu + \partial_\nu \xi_\mu$ satisfies de
  Donder's gauge condition \eqref{eqDeDonder}. In particular, $\htil$
  satisfies the wave equation
  \[
    \Box \htil_{\mu\nu} = 0.
  \]
\end{lemma}

\begin{proof}
  Straightforward calculations show that de Donder's gauge condition for
  $\htil$ imposes the following equations for $\xi_\mu$,
  \begin{subequations}
    \begin{align}
      \partial^\mu \xi_\mu & = - \frac{1}{2} h^\mu_\mu, \label{eqTraceless}                                        \\
      \Box \xi_\nu         & = - \partial^\mu h_{\mu\nu} + \frac{1}{2} \partial_\nu h^\mu_\mu.\label{eqTransverse}
    \end{align}
  \end{subequations}
  We remind the reader that the operator $\Box$ is a surjection from the
  set of homogeneous polynomials of degree $p+3$ to the set of
  homogeneous polynomials of degree $p+1$, so
  Equation~\eqref{eqTransverse} admits a (non-unique) solution
  $\xi_\mu^0$. Upon changing $h_{\mu\nu}$ to
  $h'_{\mu\nu} = h_{\mu\nu} + \partial_\mu \xi^0_\nu + \partial_\nu \xi^0_\mu$,
  we can assume that Equation~\eqref{eqTransverse} reduces to
  \begin{equation} \label{eqTransverse2}\tag{\ref{eqTransverse}'}
    \Box \xi_\nu
    = - \partial^\mu h_{\mu\nu} + \frac{1}{2} \partial_\nu h^\mu_\mu
    = 0.
  \end{equation}
  Taking the trace of Equation \eqref{eqLinEinstein} for $h$ together with
  Condition \eqref{eqTransverse2}, we obtain that
  \[
    \Box h^\mu_\mu = 0.
  \]
  As can be easily seen by taking the d'Alembertian of
  Equation~\eqref{eqTraceless}, this condition is necessary to find a $\xi_\mu$
  solving \eqref{eqTraceless}-\eqref{eqTransverse}. Conversely, since the space
  of harmonic homogeneous polynomials of degree $p+2$ is an irreducible
  representation of $O(n, 1)$, it suffices to check that a single such polynomial
  $-\frac{1}{2} h^\mu_\mu$ can be written as $\partial^\mu \xi_\mu$. For example,
  if
  \[
    -\frac{1}{2} h^\mu_\mu = (X^0-X^1)^{p+2},
  \]
  then we can choose
  \[
    \xi = -\frac{1}{2(p+3)}(X^0 - X^1)^{p+3} (dX^0 + dX^1)
  \]
  so that
  \[
    \Box \xi = 0 \text{ and } \partial^\mu \xi_\mu = -\frac{1}{2} h_\mu^\mu.
  \]
\end{proof}

Note that de Donder's gauge condition does not determine $\xi_\mu$ completely,
yet it serves an important purpose. Namely, any polynomial solution to the
linearized Einstein equation has representatives modulo infinitesimal
diffeomorphisms in the space $\cH_{p+2} \otimes \mathring{\Sym}_2(\bR^{n, 1})$,
that is in the (dual of) the tensor product
\begin{center}
  \gyoung(\ _3\hdts\ ) $\otimes$ \yng(2) .
\end{center}

The next step will be to decompose this tensor product into a sum of
irreducible representations. The decomposition of the tensor product of two
irreducible representations of $O_\uparrow(n, 1)$ is given by a generalization
of the classical Littlewood-Richardson rules obtained in \cite{Littelmann}. We
need some further information. Then we will have to check which irreducible
representation satisfies de Donder's gauge condition (this will ensure that the
corresponding irreducible representation is a set of solutions to the Einstein
equations) and does not belong to the kernel of the operator $\riem$. It should
be noted that checking either de Donder's gauge condition or the non-triviality
of $\riem$ for an arbitrary irreducible representation can be done at the level
of the highest weight vector. Indeed, both de Donder's gauge condition and
$\riem$ can be viewed as intertwining maps for the action of $O_\uparrow(n,
  1)$. And, under such a map, the image of an irreducible representation is
either $0$ or is isomorphic to the irreducible representation itself.

The decomposition of $\cH_{p+2} \otimes \mathring{\Sym}_2(\bR^{n, 1}) \otimes
  \bC$ depends on whether $n = 3$ or $n \geq 4$.

\begin{lemma}[The case $n = 3$]\label{lmTensorProde3}
  Assuming that $n = 3$, the representation $\cH_{p+2} \otimes\mathring{\Sym}_2(\bR^{3, 1})$
  decomposes as the sum
  \[\begin{split}
      \cH_{p+2} \otimes \mathring{\Sym}_2(\bR^{3, 1}) \otimes \bC
       & =
      V_{(p+4) \omegabar_1 + (p+4) \omegabar_2}
      \oplus V_{(p+4) \omegabar_1 + (p+2) \omegabar_2}
      \oplus V_{(p+2)\omegabar_1 + (p+4) \omegabar_2} \\
       & \qquad
      \oplus V_{(p+4) \omegabar_1 + p \omegabar_2}
      \oplus V_{(p+2) \omegabar_1 + (p+2) \omegabar_2}
      \oplus V_{p\omegabar_1 + (p+4) \omegabar_2}     \\
       & \qquad
      \oplus V_{(p+2) \omegabar_1 + p \omegabar_2}
      \oplus V_{p \omegabar_1 + (p+2)\omegabar_2}
      \oplus V_{p\omegabar_1 + p\omegabar_2},
    \end{split} \]
  of irreducible $SO_\uparrow(3, 1)$-representations. Here $V_\omega$ denotes the
  irreducible representation with highest weight $\omega$. The first highest
  weight vectors are
  \begin{equation}\label{eqHighestWeighte3}
    \begin{aligned}
      H_{(p+4) \omegabar_1 + (p+4) \omegabar_2}
       & = (X^0 + X^1)^{p+2} (dX^0 + dX^1)^2,                                                               \\
      H_{(p+2) \omegabar_1 + (p+4) \omegabar_2}
       & = (X^0 + X^1)^{p+1}(X^2 + i X^3) (dX^0 + dX^1)^2                                                   \\
       & \qquad - (X^0 + X^1)^{p+2} (dX^0 + dX^1) (dX^2 + i dX^3),                                          \\
      H_{(p+4) \omegabar_1 + (p+2) \omegabar_2}
       & = (X^0 + X^1)^{p+1}(X^2 - i X^3) (dX^0 + dX^1)^2                                                   \\
       & \qquad - (X^0 + X^1)^{p+2} (dX^0 + dX^1) (dX^2 - i dX^3),                                          \\
      H_{p \omegabar_1 + (p+4) \omegabar_2}
       & = (X^0 + X^1)^p \left[(X^2 + i X^3) (dX^0 + dX^1) - (X^0 + X^1) (dX^2 + i dX^3)\right]^2,          \\
      H_{(p+4) \omegabar_1 + p \omegabar_2}
       & = (X^0 + X^1)^p \left[(X^2 - i X^3) (dX^0 + dX^1) - (X^0 + X^1) (dX^2 - i dX^3)\right]^2,          \\
      h_{(p+2) \omegabar_1 + (p+2) \omegabar_2}
       & = \frac{p+2}{2} (X^0 + X^1)^{p+2} \left((dX^0)^2 - (dX^1)^2 + (dX^2)^2 + (dX^3)^2\right)           \\
       & \qquad - 2 (p+2) (X^0 + X^1)^{p+1} (dX^0 + dX^1) (X^2 dX^2 + X^3 dX^3)                             \\
       & \qquad + \left[(p+1) ((X^2)^2 + (X^3)^2) + (X^0)^2 - (X^1)^2\right] (X^0 + X^1)^p (dX^0 + dX^1)^2.
    \end{aligned}
  \end{equation}
\end{lemma}

It should be noted that we only give the formulas for the first highest weight
vectors. The last three highest weight vectors can be computed but their
expressions are rather long and unenlightening. The proof of
Lemma~\ref{lmTensorProde3} is a simple exercise of the computation of the
Clebsch-Gordan coefficients. For more details on this method, the interested
reader can consult a standard textbook on quantum mechanics such as
\cite[Chapter~X]{CohenTannoudji}. As an example, the vectors $H_{(p+2)
      \omegabar_1 + (p+4) \omegabar_2}$ and $H_{p \omegabar_1 + (p+4) \omegabar_2}$
can be computed to be
\begin{align*}
  H_{(p+2) \omegabar_1 + (p+4) \omegabar_2}
   & = \left(f_1 \otimes 1 - \frac{p+2}{2} 1 \otimes f_1\right)
  H_{(p+4) \omegabar_1 + (p+4) \omegabar_2}                                                         \\
  H_{p \omegabar_1 + (p+4) \omegabar_2}
   & = \left(\frac{1}{p+1} f_1^2 \otimes 1 - f_1 \otimes f_1 + \frac{p+2}{2} 1 \otimes f_1^2\right)
  H_{(p+4) \omegabar_1 + (p+4) \omegabar_2}.
\end{align*}

\begin{lemma}[The case $n \geq 4$]\label{lmTensorProdg3}
  Assuming that $n \geq 4$, the representation
  $\cH_{p+2} \otimes \mathring{\Sym}_2(\bR^{n, 1})$
  decomposes as the sum
  \[ \begin{split}
      \cH_{p+2} \otimes \mathring{\Sym}_2(\bR^{n, 1}) \otimes \bC
       & =
      V_{(p+4) \omegabar_1} \oplus V_{(p+2) \omegabar_1 + \omegabar_2} \oplus V_{p \omegabar_1 + 2 \omegabar_2} \\
       & \qquad
      \oplus V_{(p+2) \omegabar_1} \oplus V_{p \omegabar_1 + \omegabar_2}                                       \\
       & \qquad
      \oplus V_{p \omegabar_1},
    \end{split} \]
  of irreducible $SO_\uparrow(n, 1)$-representations. Here $V_\omega$ denotes the
  irreducible representation with highest weight $\omega$. The highest weight
  vectors are the following,
  \begin{equation}\label{eqHighestWeightg3}
    \begin{aligned}
      H_{(p+4) \omegabar_1}
       & = (Z^{-1})^{p+2} dZ^{-1} \otimes dZ^{-1},                                                   \\
      H_{(p+2) \omegabar_1 + \omegabar_2}
       & = (Z^{-1})^{p+1} Z^{-2} dZ^{-1} \otimes dZ^{-1}
      - \frac12 (Z^{-1})^{p+2}\left(dZ^{-1} \otimes dZ^{-2} + dZ^{-2} \otimes dZ^{-1}\right),        \\
      H_{p \omegabar_1 + 2 \omegabar_2}
       & = (Z^{-1})^{p+2} dZ^{-2} \otimes dZ^{-2} - Z^{-2} (Z^{-1})^{p+1}
      \left(dZ^{-1} \otimes dZ^{-2} + dZ^{-2} \otimes dZ^{-1}\right)                                 \\
       & \qquad + (Z^{-1})^p (Z^{-2})^2 dZ^{-1} \otimes dZ^{-1},                                     \\
      h_{(p+2) \omegabar_1}
       & = \frac12 (Z^{-1})^{p+1} \sum_{j=-l}^{l} Z^j
      \left(dZ^{-1} \otimes dZ^{-j} +  dZ^{-j} \otimes dZ^{-1}\right)                                \\
       & \qquad - \frac{p+1}{2p + n+1} (Z^{-1})^p \sum_{j=-l}^{l} Z^j Z^{-j} dZ^{-1} \otimes dZ^{-1} \\
       & \qquad - \frac{1}{n+1} (Z^{-1})^{p+2} \sum_{j=-l}^{l} dZ^j \otimes dZ^{-j},                 \\
      h_{p \omegabar_1 + \omegabar_2}
       & = \frac12 \sum_{j=-l}^{l} (Z^{-1})^{p+1} Z^j
      \left(dZ^{-2} \otimes dZ^{-j} + dZ^{-j} \otimes dZ^{-2}\right)                                 \\
       & \qquad - \frac12 \sum_{j=-l}^{l} (Z^{-1})^p Z^{-2} Z^j
      \left(dZ^{-1} \otimes dZ^{-j} + dZ^{-j} \otimes dZ^{-1}\right)                                 \\
       & \qquad - \frac{p}{2(2p+n+1)} \left(\sum_{j=-l}^{l} Z^j Z^{-j}\right) \times                 \\
       & \qquad\qquad
      \left[(Z^{-1})^p \left(dZ^{-1} \otimes dZ^{-2} + dZ^{-2} \otimes dZ^{-1}\right)
      - 2 (Z^{-1})^{p-1} Z^{-2} dZ^{-1} \otimes dZ^{-1}\right],                                      \\
      k_{p \omegabar_1}
       & = (Z^{-1})^p \left(\sum_{j=-l}^{l} Z^j dZ^{-j}\right)^{\otimes 2}                           \\
       & \qquad - \frac{1}{2p+n+1} \left(\sum_{a=-l}^{l} Z^a Z^{-a}\right)
       \left[ p (Z^{-1})^{p-1} \sum_{j=-l}^{l} Z^j
      \left(dZ^{-1} \otimes dZ^{-j} + dZ^{-j} \otimes dZ^{-1}\right)\right.                          \\
       & \hphantom{\qquad - \frac{1}{2p+n+1} \left(\sum_{a=-l}^{l} Z^a Z^{-a}\right) \left[\right.}
      \left. + (Z^{-1})^p \sum_{j=-l}^{l} dZ^j \otimes dZ^{-j}\right]                                \\
       & \qquad + \frac{p(p-1)}{(2p+n+1)(2p+n-1)} (Z^{-1})^{p-2}
      \left(\sum_{j=-l}^{l} Z^j Z^{-j}\right)^2 dZ^{-1} \otimes dZ^{-1}.
    \end{aligned}
  \end{equation}
  Here $Z^i$ denotes the coordinates on $\bR^{n, 1}$ associated to the basis
  $(e_j)$ introduced in Subsection~\ref{subsecGeneral}; in the sums above,
  the index $j$ ranges over the basis indices, namely
  $j \in \{-l, \ldots, -1, 0, 1, \ldots, l\}$ when $n+1 = 2l+1$ is odd and
  $j \in \{-l, \ldots, -1, 1, \ldots, l\}$ when $n+1 = 2l$ is even.
\end{lemma}

\begin{proof}
  The proof of this lemma is fairly straightforward, but tedious.
  Therefore, we give only the strategy.
  From \cite{Littelmann} we know that
  \begin{equation}\label{eqTensorProduct}
    \begin{aligned}
      \begin{tikzpicture}[baseline=-.375cm, scale=1]
        \Yboxdim{0.5cm}
        \tgyoung(0cm, -.5cm,\ _3\hdts\ )
        \draw[decorate,decoration={brace}] (0, 0.2) -- (2.5, 0.2);
        \node[above] at (1.25, 0.2) {$p$ boxes};
      \end{tikzpicture}
      \otimes \yng(2)
       & =
      \begin{tikzpicture}[baseline=-.375cm, scale=1]
        \Yboxdim{0.5cm}
        \tgyoung(0cm, -.5cm,\ _3\hdts\ )
        \draw[decorate,decoration={brace}] (0, 0.2) -- (2.5, 0.2);
        \node[above] at (1.25, 0.2) {$p+4$ boxes};
      \end{tikzpicture}
      \oplus
      \begin{tikzpicture}[baseline=-.625cm, scale=1]
        \Yboxdim{0.5cm}
        \tgyoung(0cm, -.5cm,\ _3\hdts\ ,\ )
        \draw[decorate,decoration={brace}] (0, 0.2) -- (2.5, 0.2);
        \node[above] at (1.25, 0.2) {$p+3$ boxes};
      \end{tikzpicture}
      \oplus
      \begin{tikzpicture}[baseline=-.625cm, scale=1]
        \Yboxdim{0.5cm}
        \tgyoung(0cm, -.5cm,\ \ _3\hdts\ ,\ \ )
        \draw[decorate,decoration={brace}] (0, 0.2) -- (3, 0.2);
        \node[above] at (1.25, 0.2) {$p+2$ boxes};
      \end{tikzpicture} \\
       & \qquad\oplus
      \begin{tikzpicture}[baseline=-.375cm, scale=1]
        \Yboxdim{0.5cm}
        \tgyoung(0cm, -.5cm,\ _3\hdts\ )
        \draw[decorate,decoration={brace}] (0, 0.2) -- (2.5, 0.2);
        \node[above] at (1.25, 0.2) {$p+2$ boxes};
      \end{tikzpicture}
      \oplus
      \begin{tikzpicture}[baseline=-.625cm, scale=1]
        \Yboxdim{0.5cm}
        \tgyoung(0cm, -.5cm,\ _3\hdts\ ,\ )
        \draw[decorate,decoration={brace}] (0, 0.2) -- (2.5, 0.2);
        \node[above] at (1.25, 0.2) {$p+1$ boxes};
      \end{tikzpicture}
      \oplus
      \begin{tikzpicture}[baseline=-.375cm, scale=1]
        \Yboxdim{0.5cm}
        \tgyoung(0cm, -.5cm,\ _3\hdts\ )
        \draw[decorate,decoration={brace}] (0, 0.2) -- (2.5, 0.2);
        \node[above] at (1.25, 0.2) {$p$ boxes};
      \end{tikzpicture}
    \end{aligned}
  \end{equation}
  This shows that the highest weight vectors in the tensor product
  $\cH_{p+2} \otimes \mathring{\Sym}_2(\bR^{n, 1}) \otimes \bC$ have the
  highest weights indicated in the statement of the lemma. The first
  three irreducible representations that appear in the right hand side
  of Equation \eqref{eqTensorProduct} have $p+4$ boxes so the
  corresponding highest weight vector in $\cH_{p+2} \otimes
    \mathring{\Sym}_2(\bR^{n, 1})$ can be obtained by using the
  corresponding Young symmetrizer, see
  \cite[Lemma~9.3.2]{GoodmanWallach}. For the fourth and the fifth terms
  in Equation \eqref{eqTensorProduct}, the same procedure leads to
  elements in $\cH_{p+1} \otimes (\bR^{n, 1})^* \otimes \bC$. One then
  has to multiply them by the ``metric'' $\sum_{j=-l}^{l} Z^j dZ^{-j}$. However,
  these elements do not belong to
  $\cH_{p+2} \otimes \mathring{\Sym}_2(\bR^{n, 1}) \otimes \bC$ since
  neither its d'Alembertian nor its trace are zero so we need to subtract
  them by adding terms that are proportional to $\sum_{j=-l}^{l} Z^j Z^{-j}$ and
  to $\sum_{j=-l}^{l} dZ^j \otimes dZ^{-j}$.  All these maps being intertwining
  for $O(n+1, \bC)$, we have then constructed the highest weight vector
  of the corresponding representation. The procedure is similar for the
  sixth term in \eqref{eqTensorProduct}. However, here, we have to
  multiply by $\left(\sum_{j=-l}^{l} Z^j dZ^{-j}\right)^{\otimes 2}$.
\end{proof}

The last step in the proof of Proposition \ref{propWeyl} is to identify $\cW_p$
with the irreducible representation of highest weight $p \omegabar_1 + 2
  \omegabar_2$ if $n \geq 4$ or to $V_{p \omegabar_1 + (p+4) \omegabar_2} \oplus
  V_{(p+4) \omegabar_1 + p \omegabar_2}$ if $n = 3$.

From Lemma~\ref{lmDeDonder}, we know that any $W \in \cW_p$ can be represented
as $W = \riem(h)$ for some $h \in \cE_{p+2}$ satisfying de Donder's gauge
condition. In particular, $h$ belongs to $\cH_{p+2} \otimes
  \mathring{\Sym}_2(\bR^{n, 1})$. This means that $\cW_p$ is the image of the set
of $O_\uparrow(n, 1)$ irreducible subrepresentations of $\cH_{p+2} \otimes
  \mathring{\Sym}_2(\bR^{n, 1})$ that satisfy de Donder's gauge condition. It can
be checked that only the first three vectors in \eqref{eqHighestWeightg3},
which we denoted $H_\omega$ to distinguish them from the others, and similarly
only the first five vectors in \eqref{eqHighestWeighte3} satisfy this condition
(as mentioned earlier, this gauge condition must only be verified for the
highest weight vector).

We now assume $n \geq 4$, the case $n = 3$ being similar. The first two terms
actually belong to the kernel of $\riem$. The simplest way to check this is by
noticing that
\begin{align*}
  H_{(p+4) \omegabar_1}                & = \lie_{\xi_{(p+2) \omegabar_1}} \eta,               \\
  H_{(p+2) \omegabar_1 + \omegabar_2 } & = \lie_{\xi_{(p+2) \omegabar_1 + \omegabar_2}} \eta,
\end{align*}
where
\[
  \begin{aligned}
    \xi_{(p+4) \omegabar_1} & = \frac{1}{2(p+3)} (Z^{-1})^{p+3} e_{+1},                                          \\
    \xi_{(p+2) \omegabar_1 + \omegabar_2 }
                            & = \frac{1}{2(p+2)} \left((Z^{-1})^{p+2} Z^{-2} e_1 - (Z^{-1})^{p+3} e_{+2}\right).
  \end{aligned}
\]
So only the third term in \eqref{eqTensorProduct} remains and it is a simple
matter to check that
\[
  \riem(H_{p\omegabar_1 + 2 \omegabar_2})(e_{-1}, e_{-2}, e_{-1}, e_{-2})
  = (p+2)(p+3) (Z^{-1})^p \neq 0.
\]
This proves that $\cW_p$ is the image of $V_{p \omegabar_1 + 2 \omegabar_2}$
under $\riem$, and thus completing the proof of Proposition~\ref{propWeyl}.

The dimension of $\cW_p$ is
\[
  \dim \cW_p =
  \frac{1}{2} \frac{n+1}{n-1} \binom{p + n}{p+3}
  \frac{(p+1)(p+n+2)(2p+n+3)}{p+n}.
\]
This formula follows once again by the Weyl dimension formula. Note that the
case $n+1=4$ has to be checked separately since the representation $\cW_p = (p,
  p+4) \oplus (p+4, p)$ is not irreducible for $n+1=4$.

\subsubsection{Invariant quadratic form}
\label{secWeylForm}

To compute the signature of the invariant quadratic form on $\cW_p$, we use the
Cartan decomposition of $O_\uparrow(n, 1)$ (see for example
\cite{KnappRepresentationTheory} for a presentation of the Cartan
decomposition). The maximal compact subgroup $K$ of $O_\uparrow(n, 1)$ is
$O(n)$, where it is understood that elements in $O(n)$ leave invariant the time
direction $X^0$. Let $\mathfrak{p}$ denote the subspace of the Lie algebra
$\mathfrak{so}(n, 1)$ generated by the boosts $a_i$, $i=1, \dots, n$.

The subspace $\mathfrak{p}$ is left invariant by the adjoint action of $O(n)$.
As such, $\mathfrak{p}$ is a representation of $O(n)$ which turns out to be
equivalent to the standard one.

From the classical branching rules (see for example
\cite[Chapter~8]{GoodmanWallach}), we know how $\cW_p$ decomposes into
representations of $O(n)$,
\begin{equation}\label{eqBranchingRule}
  \begin{aligned}
    \mathrm{Res}^{O_\uparrow(n, 1)}_{O(n)} \left(\cW_p\right)
     & =
    \begin{tikzpicture}[baseline=-.625cm, scale=1]
      \Yboxdim{0.5cm}
      \tgyoung(0cm, -.5cm,\ \ _3\hdts\ ,\ \ )
      \draw[decorate,decoration={brace}] (0, 0.2) -- (3, 0.2);
      \node[above] at (1.25, 0.2) {$p+2$ boxes};
    \end{tikzpicture}
    \oplus \cdots \oplus
    \begin{tikzpicture}[baseline=-.625cm, scale=1]
      \Yboxdim{0.5cm}
      \tgyoung(0cm, -.5cm,\ \ \ ,\ \ )
    \end{tikzpicture}
    \oplus
    \begin{tikzpicture}[baseline=-.625cm, scale=1]
      \Yboxdim{0.5cm}
      \tgyoung(0cm, -.5cm,\ \ ,\ \ )
    \end{tikzpicture} \\
     & \qquad\oplus
    \begin{tikzpicture}[baseline=-.625cm, scale=1]
      \Yboxdim{0.5cm}
      \tgyoung(0cm, -.5cm,\ \ _3\hdts\ ,\ )
      \draw[decorate,decoration={brace}] (0, 0.2) -- (3, 0.2);
      \node[above] at (1.25, 0.2) {$p+2$ boxes};
    \end{tikzpicture}
    \oplus \cdots \oplus
    \begin{tikzpicture}[baseline=-.625cm, scale=1]
      \Yboxdim{0.5cm}
      \tgyoung(0cm, -.5cm,\ \ \ ,\ )
    \end{tikzpicture}
    \oplus
    \begin{tikzpicture}[baseline=-.625cm, scale=1]
      \Yboxdim{0.5cm}
      \tgyoung(0cm, -.5cm,\ \ ,\ )
    \end{tikzpicture} \\
     & \qquad \oplus
    \begin{tikzpicture}[baseline=-.375cm, scale=1]
      \Yboxdim{0.5cm}
      \tgyoung(0cm, -.5cm,\ _3\hdts\ )
      \draw[decorate,decoration={brace}] (0, 0.2) -- (2.5, 0.2);
      \node[above] at (1.25, 0.2) {$p+2$ boxes};
    \end{tikzpicture}
    \oplus \cdots \oplus
    \begin{tikzpicture}[baseline=-.375cm, scale=1]
      \Yboxdim{0.5cm}
      \tgyoung(0cm, -.5cm,\ \ \ )
    \end{tikzpicture}
    \oplus
    \begin{tikzpicture}[baseline=-.375cm, scale=1]
      \Yboxdim{0.5cm}
      \tgyoung(0cm, -.5cm,\ \ )
    \end{tikzpicture}
  \end{aligned}
\end{equation}
An important fact to notice is that the restriction is multiplicity
free and this decomposition holds true considering either the real
representation $\cW_p$ or the complex representation $\cW_p \otimes
  \bC$.

We now study the action of $\mathfrak{p}$ closer. The mapping $a \otimes W
  \mapsto a \cdot W$ from $\mathfrak{p} \otimes \cW_p$ to $\cW_p$ is
$O(n)$-equivariant meaning that for any $R \in O(n)$, we have $Ad(R)(a) \cdot
  (R_* W) = R_* (a \cdot W)$ (we remind the reader that $O_\uparrow(n, 1)$ acts
on $\cW_p$ by push-forward). As a consequence, elements belonging to an
irreducible representation $V_0$ of $O(n)$ in the right-hand side of
\eqref{eqBranchingRule} are mapped by $\mathfrak{p}$ into the subspace of
$\mathrm{Res}^{O_\uparrow(n,1)}_{O(n)}$ which is the direct sum of irreducible
representations that appear both in the right-hand side of
\eqref{eqBranchingRule} and in the decomposition of the tensor product $V_0
  \otimes \bR^n$ into irreducible representations of $O(n)$. As an example,
\[
  \mathfrak{p} \quad \text{maps} \quad
  \Yboxdim{0.25cm}
  \yng(3,2)
  \quad \text{to} \quad
  \yng(4,2) \oplus \yng(2,2) \oplus \yng(3,1)
  \quad .
\]
What is important to notice is that $O(n)$-representations with an odd number
of boxes are mapped to representations with an even number of boxes and vice
versa by $\mathfrak{p}$.

The invariant quadratic form $q$ on $\cW_p$ restricts to each element in the
right-hand side of Equation \eqref{eqBranchingRule}. Being real representation
of $O(n)$, each of them carry a unique (up to normalization) $O(n)$-invariant
quadratic form which can be chosen to be positive definite. We conclude that
the restriction of $q$ to any element in the right-hand side of
\eqref{eqBranchingRule} is either positive or negative definite.

The Lie group with Lie algebra $\mathfrak{g} = \mathfrak{o}(n) \oplus i
  \mathfrak{p}$ is $SO(n+1)$, the compact form of $SO_\uparrow(n, 1)$. The
representation $\cW_p \otimes \bC$ is also of real type for $SO(n+1)$. A real
form of $\cW_p \otimes \bC$ for $SO(n+1)$ can be constructed as follows. Let us
abuse notation and assume that a Young diagram represents the real
representation inside $\cW_p$ that is isomorphic to the image of the Young
projector.

Then a real form for $SO(n+1)$ is given by
\[
  \widetilde{\cW}_p = \Yboxdim{0.25cm}\left[
    \yng(2) \oplus \yng(2,2) \oplus \yng(3,1) \oplus \cdots\right]
  \oplus i \left[\yng(2,1) \oplus \yng(3) \oplus \yng(3,2) \oplus \cdots\right],
\]
namely, we put an $i$ in front of all subrepresentations with an odd number of
boxes. Since $SO(n+1)$ is a compact Lie group, $\widetilde{\cW}_p$ carries a
unique invariant quadratic form $q$ which is positive definite.

Looking at the invariant quadratic form $q$ on $\cW_p$ this means that $q$ is
positive definite on the $O(n)$-irreducible representations with an even number
of boxes and negative definite on the representations with an odd number of
boxes.

If we let $(n_+(p), n_-(p))$ denote the signature of $q$ on $\cW_p$, we have
$n_+(p) + n_-(p) = \dim \cW_p$ and $n_+(p) - n_-(p)$ can be computed by
induction on $p$. Indeed,
\begin{align*}
   & n_+(p) - n_-(p) = n_+(p-1) - n_-(p-1) \\
   & \qquad + (-1)^p \left[\dim
    \begin{tikzpicture}[baseline=-.375cm, scale=1]
      \Yboxdim{0.5cm}
      \tgyoung(0cm, -.5cm,\ _3\hdts\ )
      \draw[decorate,decoration={brace}] (0, 0.2) -- (2.5, 0.2);
      \node[above] at (1.25, 0.2) {$p+2$ boxes};
    \end{tikzpicture}
    - \dim
    \begin{tikzpicture}[baseline=-.625cm, scale=1]
      \Yboxdim{0.5cm}
      \tgyoung(0cm, -.5cm,\ \ _3\hdts\ ,\ )
      \draw[decorate,decoration={brace}] (0, 0.2) -- (3, 0.2);
      \node[above] at (1.25, 0.2) {$p+2$ boxes};
    \end{tikzpicture}
    + \dim
    \begin{tikzpicture}[baseline=-.625cm, scale=1]
      \Yboxdim{0.5cm}
      \tgyoung(0cm, -.5cm,\ \ _3\hdts\ ,\ \ )
      \draw[decorate,decoration={brace}] (0, 0.2) -- (3, 0.2);
      \node[above] at (1.25, 0.2) {$p+2$ boxes};
    \end{tikzpicture}
    \right].
\end{align*}
For $p = 0$, we find
\[
  n_+(0) - n_-(0) = \frac{1}{12} (n+2) (n-1)(n-2)(n-3).
\]
The following formula can be obtained by induction,
\[
  n_+(p) - n_-(p) = (-1)^p \frac{p+1}{2} (p+n+2) \binom{p+n-1}{p+3}.
\]
At this point it should be noted that it is much more convenient to change the
convention and assume that $q$ is positive definite if $p$ is even and negative
definite if $p$ is odd. That has the effect of removing the annoying factor
$(-1)^p$. With this new convention, we get
\[
  \left\lbrace
  \begin{aligned}
    n_+(p) & = \frac12 (n^2 + (n+1)p + 3)\frac{(p+1)(p+n+2)}{(n-1)(p+n)} \binom{p+n}{p+3}, \\
    n_-(p) & = \frac12 (np + 4n + p)\frac{(p+1)(p+n+2)}{(n-1)(p+n)} \binom{p+n}{p+3}.
  \end{aligned}
  \right.
\]